\documentclass[mathpazo]{cmr}

\usepackage{bm}
\usepackage[pdftex]{graphicx}
\usepackage{here}
\usepackage{subcaption}
\usepackage{caption}
\usepackage{setspace}
\graphicspath{{./picture/}}

\begin{document}
\title{A structure-preserving numerical method\\for the fourth-order geometric evolution equations for planar curves}


\author[ E. Miyazaki, T. Kemmochi, T. Sogabe, and S.-L. Zhang]{ Eiji Miyazaki\affil{1},
                                      Tomoya Kemmochi\affil{1}, Tomohiro Sogabe\affil{1}and  Shao-Liang Zhang\affil{1} }
\address{\affilnum{1}\ 
                      Graduate School of Engineering,
                      Nagoya University,
                      Furo-cho, Chikusa-ku, Nagoya, Aichi 464-8601 Japan \\
        }
\emails{{\tt zhang@na.nuap.nagoya-u.ac.jp} (E.Miyazaki),
       {\tt chan@email}  (A.~Chan),
       {\tt zhao@email}  (A.~Zhao)  }


\begin{abstract}
        For fourth-order geometric evolution equations for planar curves with the dissipation of the bending energy, including the Willmore and the Helfrich flows, we consider a numerical approach. In this study, we construct a structure-preserving method based on a discrete variational derivative method. Furthermore, to prevent the vertex concentration that may lead to numerical instability, we discretely introduce Deckelnick's tangential velocity. Here, a modification term is introduced in the process of adding tangential velocity. This modified term enables the method to reproduce the equations' properties while preventing vertex concentration. Numerical experiments demonstrate that the proposed approach captures the equations' properties with high accuracy and avoids the concentration of vertices.
\end{abstract}

\ams{37M15, 65M06}    
\keywords{Geometric evolution equation, Willmore flow, Helfrich flow, Numerical calculation, Structure-preserving, Discrete variational derivative method, Tangential velocity}    

\maketitle



\section{Introduction}
In this paper, we design a numerical method for the Willmore and the Helfrich flows for planar curves. The Willmore flow is a geometric evolution equation that models the behavior of elastic bodies \cite{20}. It is a gradient flow with regard to the bending energy
\begin{equation}
  \label{dansei}
B(t)=\frac{1}{2}\int_{C(t)}(k-c_0)^2ds,
\end{equation}
where $C(t)$ is a closed curve, $s$ is the arc-length parameter of $C(t)$, $k$ is the curvature of $C(t)$, and $c_0$ is a given constant \cite{17,18}. The Willmore flow is given by 
\begin{equation}
  \label{will}
  \frac{\partial\boldsymbol{X}}{\partial{t}}=   -\boldsymbol{\delta} B,\quad
  \boldsymbol{\delta} B= -\left\{\frac{\partial^2 k}{\partial s^2}+\frac{1}{2}(k-c_0)^3\right\}\boldsymbol{N},
\end{equation}
where $\boldsymbol{\delta}B$ is the gradient of $B$, $\boldsymbol{X}$ is a point on the closed curve $C(t)$, and $\boldsymbol{N}$ is the unit outward normal vector of $C(t)$.

The Helfrich flow is a gradient flow for the bending energy under the constraint that the length and the enclosed area of the curve are conserved. The Helfrich flow is given by 
\begin{equation}
  \begin{split}
  \label{helf1}
 &\frac{\partial\boldsymbol{X}}{\partial{t}}=\left\{\frac{\partial^2 k}{\partial s^2}+\frac{1}{2}(k-c_0)^3+k\lambda +\mu \right\}\boldsymbol{N},\\
  &\begin{pmatrix}
    \lambda  \\
    \mu
    \end{pmatrix}
    =\frac{1}{\langle k\rangle^2-\langle  k^2\rangle}
    \begin{pmatrix}
      1&-\langle k\rangle \\
      -\langle k\rangle &-\langle k^2\rangle 
      \end{pmatrix}
      \begin{pmatrix}
        \langle k\boldsymbol{\delta} B\rangle + c_0^2 \langle k^2\rangle /2\\
        \langle \boldsymbol{\delta} B\rangle + c_0^2 \langle k\rangle /2
        \end{pmatrix},
    \end{split}
\end{equation}
where 
\begin{equation}
  \langle  F \rangle = \frac{1}{L}\int_{C(t)} F ds, \quad L=\int_{C(t)}ds.
\end{equation} Helfrich proposed an optimization problem that mathematically models the shape of red blood cells \cite{helfrich,23}. The Helfrich flow was proposed to solve this optimization problem \cite{21,22}. Note that the Willmore and the Helfrich flows are fourth-order nonlinear evolution equations. 

Some numerical approaches have been investigated for \eqref{will} and \eqref{helf1}. There are two primary problems in the numerical computation of these flows.
\begin{itemize}
  \item Numerical computation using general-purpose approaches (e.g. the Runge--Kutta method) can become unstable if the time step size is too large. 
  \item When one approximates curves by polygonal curves, the vertices may be concentrated as the time step proceeds. The concentration of vertices may make numerical computations unstable.
\end{itemize}
Some numerical methods regarding the above problems for \eqref{will} and \eqref{helf1} are presented.
In \cite{14}, a linear numerical method for \eqref{will} and \eqref{helf1} are proposed. It is based on the finite element method and approximates a closed curve by a closed polygonal curve. Additionally, the method implicitly involves a tangential velocity through a mass lumped inner product. In \cite{19}, a semi-implicit numerical method for \eqref{will} is proposed. The method also uses a polygonal curve and introduces tangential velocity that makes the distribution of vertices asymptotically uniform, which was derived in \cite{24,25}.
These tangential velocities' impact in \cite{14,19} can prevent the polygonal curve's vertices from concentrating.  In \cite{16}, a structure-preserving numerical method for \eqref{will} is proposed. The method approximates the curve by B-spline curves \cite{B-spline}. It is based on the discrete partial derivative method \cite{DPDM} and discretely reproduces the dissipation of the bending energy. The method is a nonlinear scheme.  A method that discretely reproduces a property of an equation is called the structure-preserving numerical method. The numerical solutions obtained using these methods are not only physically correct but also known to be less likely to diverge even when relatively large time step sizes are used \cite{1}.

This study aims to construct a structure-preserving numerical method with tangential velocity for \eqref{will} and \eqref{helf1}. That is, for \eqref{will}, we construct a numerical method that discretely reproduces the bending energy's dissipation, and for \eqref{helf1}, we construct a numerical method that discretely reproduces the bending energy's dissipation, conservation of the length, and the enclosed area.

Throughout this paper, the polygonal curve is used to
discretize the curve. Additionally,
the variational derivative of the bending energy, the length, and the area are discretized based on the discrete variational derivative method (DVDM) \cite{7}. Then, these variational derivatives and the Lagrange multipliers are used to build a structure-preserving method. Additionally, to prevent vertex concentration, we discretely introduce Deckelnick's tangential velocity \cite{12}. 
Here, a modification term is introduced to make our method structure-preserving. 

The proposed method enables us to conduct numerical computations with relatively large time step sizes while avoiding the polygonal curve's vertex concentration. Furthermore, since the energy properties of equation \eqref{will} and \eqref{helf1} are discretely reproduced, it is anticipated that the solutions obtained by using our method are physically correct. The method is a nonlinear scheme. Therefore, the cost of numerical computations tends to be high. 

The rest of this paper is organized as follows. In Section 2, we collect the bending energy's dissipation of \eqref{will}, the bending energy's dissipation, and the conservation of the length and the area of \eqref{helf1}. In Section 3, we construct structure-preserving numerical methods based on the DVDM. Additionally, discrete tangential velocity is introduced. We present numerical examples in Section 4, and finally, give concluding remarks in Section 5.





\section{Properties of the target equations}

In this study, we construct a structure-preserving numerical method by reproducing the process of deriving the properties of the gradient flows in a continuous system. Specifically, we discretize variational derivatives based on DVDM, and discretely reproduce the approach of determining the Lagrange multipliers in the continuous system. In this section, we will collect the variational structure and Lagrange multipliers of target equations, and derive the properties of these equations.

\subsection{Willmore flow}
The Willmore flow is given by
\begin{equation}
        \label{will_hen}
        \frac{\partial\boldsymbol{X}}{\partial{t}}=   -\boldsymbol{\delta} B,\quad
        \boldsymbol{\delta} B= -\left\{\frac{\partial^2 k}{\partial s^2}+\frac{1}{2}(k-c_0)^3\right\}\boldsymbol{N},
\end{equation}
where $\boldsymbol{\delta}B$ represents the gradient of $B$, $\boldsymbol{X}$ represents a point on a closed curve $C(t)$, $k$ represents the curvature, $c_0$ represents a given constant and $\boldsymbol{N}$ represents the unit outward normal vector.
The Willmore flow has the bending energy's dissipation:
\begin{equation}
        \begin{split}
        \frac{dB}{dt}
        &=\int_{C(t)} \boldsymbol{\delta} B\cdot\frac{\partial \boldsymbol{X}}{\partial t}ds\\
      &=-\int_{C(t)}\left| \boldsymbol{\delta} B\right|^2ds\leq0.
        \end{split}
      \end{equation}

\subsection{Helfrich flow}
The Helfrich flow is given by
\begin{equation}
  \begin{split}
  \label{helf}
 &\frac{\partial\boldsymbol{X}}{\partial{t}}=\left\{\frac{\partial^2 k}{\partial s^2}+\frac{1}{2}(k-c_0)^3+k\lambda +\mu\right\}\boldsymbol{N},\\
  &\begin{pmatrix}
    \lambda  \\
    \mu
    \end{pmatrix}
    =\frac{1}{\langle k\rangle^2-\langle  k^2\rangle}
    \begin{pmatrix}
      1&-\langle k\rangle \\
      -\langle k\rangle &-\langle k^2\rangle 
      \end{pmatrix}
      \begin{pmatrix}
        \langle k\boldsymbol{\delta} B\rangle + c_0^2 \langle k^2\rangle /2\\
        \langle \boldsymbol{\delta} B\rangle + c_0^2 \langle k\rangle /2
        \end{pmatrix},
    \end{split}
\end{equation}
where 
\begin{equation}
  \langle  F \rangle = \frac{1}{L}\int_{C(t)} F ds, \quad L=\int_{C(t)}ds.
\end{equation}
The Helfrich flow can be rewritten in the following form using the variational derivative of the bending energy $\boldsymbol{\delta} B$, the variational derivative of the length of $C(t)$ (denoted by $\boldsymbol{\delta} L$) and the enclosed erea of $C(t)$ (denoted by $\boldsymbol{\delta} A$):
\begin{equation}
        \label{hell_hen}
        \frac{\partial\boldsymbol{X}}{\partial{t}}=-\boldsymbol{\delta} B+\lambda\boldsymbol{\delta} L+\mu\boldsymbol{\delta} A ,\quad \boldsymbol{\delta} L=k\boldsymbol{N},\quad \boldsymbol{\delta} A=\boldsymbol{N}.
\end{equation}
Hereafter, for simplicity, we denote the inner product over a curve $C(t)$ by 
\begin{equation}
        \label{inner product}
        (\boldsymbol{X},\boldsymbol{Y})=\int_{C(t)}\boldsymbol{X}\cdot \boldsymbol{Y}ds. 
\end{equation}
We obtain the expression of time derivatives of the length and the area as follows:
\begin{align}
        \label{dldt}
        \frac{d L}{d t}&=\int_{C(t)}\boldsymbol{\delta} L\cdot\frac{\partial  \boldsymbol{X}}{\partial t}ds=\left(\boldsymbol{\delta} L,\frac{\partial  \boldsymbol{X}}{\partial t}\right),\\
        \label{dadt}
         \frac{d A}{d t}&=\int_{C(t)}\boldsymbol{\delta} A\cdot\frac{\partial  \boldsymbol{X}}{\partial t}ds=\left(\boldsymbol{\delta} A,\frac{\partial  \boldsymbol{X}}{\partial t}\right).
\end{align}
Since the length and the area of $C(t)$ are conserved, we have $\frac{d L}{d t}=\frac{d A}{d t}=0$.
By substituting \eqref{hell_hen} into \eqref{dldt} and \eqref{dadt}, we obtain 
\begin{align}
        \label{lambda}
        &\left(\boldsymbol{\delta} L,-\boldsymbol{\delta} B+\lambda\boldsymbol{\delta} L+\mu\boldsymbol{\delta} A\right)=0,\\
        \label{mu}
         &\left(\boldsymbol{\delta} A,-\boldsymbol{\delta} B+\lambda\boldsymbol{\delta} L+\mu\boldsymbol{\delta} A\right)=0.
      \end{align} 
Solving \eqref{lambda} and \eqref{mu} as equations for $\lambda$ and $\mu$, we can obtain $\lambda$ and $\mu$ that conserve the curve length and the area of $C(t)$.

Additionally, The Helfrich flow also has the dissipation of the bending energy:
\begin{equation}
        \begin{split}
          \label{dbdt}
        \frac{dB}{dt}
        &=\int_{C(t)} \boldsymbol{\delta} B\cdot\frac{\partial \boldsymbol{X}}{\partial t}ds\\
        &=\int_{C(t)} \boldsymbol{\delta} B\cdot\left(  -\boldsymbol{\delta} B+\lambda\boldsymbol{\delta} L+\mu\boldsymbol{\delta} A \right)ds\\
      &=-\mathrm{det}\left(\begin{matrix}(\boldsymbol{\delta} B,\boldsymbol{\delta} B)&(\boldsymbol{\delta} B,\boldsymbol{\delta} L)&(\boldsymbol{\delta} B,\boldsymbol{\delta} A)\\
        (\boldsymbol{\delta} L,\boldsymbol{\delta} B)&(\boldsymbol{\delta} L,\boldsymbol{\delta} L)&(\boldsymbol{\delta} L,\boldsymbol{\delta} A)\\
        (\boldsymbol{\delta} A,\boldsymbol{\delta} B)&(\boldsymbol{\delta} A,\boldsymbol{\delta} L)&(\boldsymbol{\delta} A,\boldsymbol{\delta} A)
      \end{matrix} \right)
      /\mathrm{det}\left(\begin{matrix}
        (\boldsymbol{\delta} L,\boldsymbol{\delta} L)&(\boldsymbol{\delta} L,\boldsymbol{\delta} A)\\
        (\boldsymbol{\delta} A,\boldsymbol{\delta} L)&(\boldsymbol{\delta} A,\boldsymbol{\delta} A)
      \end{matrix}\right)\\
      &\leq0.
        \end{split}
      \end{equation}
In the process of above transformation, we solved \eqref{lambda} and \eqref{mu} as a linear system for $\lambda$ and $\mu$, and then substituted them.





\section{Construction of structure-preserving\\ numerical method for each equation}
In this section, various quantities related to closed curves are discretized based on the approximation of the closed curve by a closed polygonal curve.
Let us discretize each variational derivative according to the concept of DVDM and determine the Lagrange multipliers to preserve the equation's structure. The tangential velocity is also introduced for stabilization. Note that in the process of introducing the tangential velocity, we add a modification term to preserve the discrete bending energy's dissipation even after the tangential velocity is added. 

\subsection{Discretization of the curve}
Let us discretize a curve and geometric quantities according to \cite{3}.
We represent a time-evolving polygonal curve in the plane by $\Gamma(t)$ and the domain of inside $\Gamma(t)$ by $\Omega(t)$. Note that in this subsection, $t$ is excluded because we fix the time $t$ and consider various quantities defined on the polygonal curve $\Gamma(t)$. Let $N$ be the number of edges of $\Gamma$ and $\boldsymbol{X}_i, \; i=1,\dots,N$ be the $i$th vertex. The index $i$ shall increase counterclockwise. 

Next, let $\Gamma_i$ be the line segment connecting $\boldsymbol{X}_{i-1}$ and $\boldsymbol{X}_{i}$. We represent the length of $\Gamma_i$ by $r_i=|\boldsymbol{X}_i-\boldsymbol{X}_{i-1}|$ .
Then, the unit tangent vector $\boldsymbol{t}_i$ on $\Gamma_i$ is defined by
\begin{equation}
  \label{t}
\boldsymbol{t}_i=\left(\boldsymbol{X}_i-\boldsymbol{X}_{i-1}\right)/{r_i}.
\end{equation}
Additionally, the unit tangent vector $\boldsymbol{T}_i$ on the vertex $\boldsymbol{X}_i$ is defined by
\begin{equation}
  \label{T}
  \boldsymbol{T}_i=\frac{\boldsymbol{t}_{i}+\boldsymbol{t}_{i+1}}{|\boldsymbol{t}_{i}+\boldsymbol{t}_{i+1}|}.
\end{equation}

Then, the curvature at the $i$th vertex $\boldsymbol{X}_i$ is defined by 
\begin{equation}
  \label{k_i}
  k_i=\frac{\mathrm{det}\left[\delta_i^{\left<1\right>}\boldsymbol{X}_i,\delta_i^{\left<2\right>}\boldsymbol{X}_i\right]}{|\delta_i^{\left<1\right>}\boldsymbol{X}_i|^3},
\end{equation}
where
\begin{equation}
  \delta_i^{\left<1\right>}\boldsymbol{X}_i=\frac{\boldsymbol{X}_{i+1}-\boldsymbol{X}_{i-1}}{2\Delta u},\quad \delta_i^{\left<2\right>}\boldsymbol{X}_i=\frac{\boldsymbol{X}_{i+1}-2\boldsymbol{X}_{i}+\boldsymbol{X}_{i-1}}{\Delta u},\quad\Delta u=1/N.
\end{equation}
Here, $\mathrm{det}[\cdot,\cdot]$ means the following operations for two vectors $\boldsymbol{A}$ and $\boldsymbol{B}$:
\begin{equation}
  \boldsymbol{A}=\left(\begin{matrix}
    a_1\\
    a_2
    \end{matrix}\right),\quad
    \boldsymbol{B}=\left(\begin{matrix}
    b_1\\
      b_2
      \end{matrix}\right)
    \implies \mathrm{det}[\boldsymbol{A},\boldsymbol{B}]=\mathrm{det}\left[\begin{matrix}
    a_1&b_1  \\
    a_2&b_2
    \end{matrix}
    \right].
\end{equation}
As above, the curvature is discretized by the central differencing of the derivatives.
Then, using discretized curvature \eqref{k_i}, the discrete bending energy is defined by 
\begin{equation}
  B_d(\Gamma)=\frac{1}{2}\sum_{i=1}^N(k_i-c_0)^2\hat{r}_i,
\end{equation}
where $\hat{r}_i=(r_i+r_{i+1})/2$.
The length $L_d$ and the enclosed area $A_d$ of $\Gamma $ are defined by
\begin{align}
  L_d(\Gamma)=\sum_{i=1}^Nr_i,\quad A_d(\Gamma)=\frac{1}{2}\sum_{i=1}^N \mathrm{det}\left[\boldsymbol{X}_{i-1},\boldsymbol{X}_{i}\right].
\end{align}
\subsection{Discrete tangential velocity}
In this section, discrete tangential velocity is defined. It is known that tangential velocity in the continuous system does not influence the shape of the curve \cite{gage}. Thus, those are used as a means of examining the equations \cite{12}. In a discrete system, if the number of vertices $N$ is sufficiently large, the impact on the shape of the polygonal curve is considered small. Numerical computation of time-evolving curves usually leads to concentration of vertices, which induces numerical instability if the tangent velocity is not taken properly. Various numerical methods have demonstrated the efficiency of introducing discrete tangential velocities\cite{24,25,3,13}. 

In this study, the tangent velocity of Deckelnick\cite{12} is used in the discrete system to prevent vertex concentration.
Deckelnick's tangential velocity is given by 
\begin{equation}
  \label{w_c}
  w(u,t)=-\frac{\partial}{\partial u}\left(\frac{1}{|\partial_u\boldsymbol{X}|}\right).
\end{equation}
We discretize \eqref{w_c} as follows
\begin{equation}
  \label{w}
  w_i(t)=-\alpha\delta_i^{-}\left(\frac{1}{|\delta_i^{+}\boldsymbol{X}_i(t)|}\right),
\end{equation}
where \[\quad\delta_i^{+}\boldsymbol{X}_{i}=\frac{\boldsymbol{X}_{i+1}-\boldsymbol{X}_{i}}{\Delta u},\quad\delta_i^{-}\boldsymbol{X}_{i}=\frac{\boldsymbol{X}_{i}-\boldsymbol{X}_{i-1}}{\Delta u},\]
and $w_i(t)$ represents the discretized tangential velocity at the vertex $\boldsymbol{X}_i(t)$, and $\alpha$ is a positive constant. 
 The tangential velocity for a vertex $\boldsymbol{X}_i(t)$ reduces the difference in distance to its two neighboring vertices $\boldsymbol{X}_{i-1}(t)$ and $\boldsymbol{X}_{i+1}(t)$. The larger the value of $\alpha$, the greater the impact of tangential velocity. We must modify the value of $\alpha$ according to the equations and conditions of the numerical computation, including the curves' initial arrangement.

 \subsection{Discretization of time and construction of structure-preserving schemes}
 We present the $n$-step time by $t^{(n)}=n\Delta t$, where $\Delta t$ represents the time step size.
 The quantity at time $t^{(n)}$ is described by the superscript $n$. In this subsection, \eqref{will_hen} and \eqref{hell_hen} are discretized temporally to preserve the structure.
 According to the concept of DVDM, the variational derivative in the discrete system is defined as follows:
\begin{definition}
  \label{definition:henbun}
  If the time-difference of the energy $E^{(n)}$ in the discrete system can be transformed into the following form,
  \begin{equation}
        \frac{E^{(n+1)}-E^{(n)}}{\Delta t}=\sum_{i=1}^N\boldsymbol{\delta} _dE_i^{(n+1)}\cdot \frac{\boldsymbol{X}_i^{(n+1)}-\boldsymbol{X}_i^{(n)}}{\Delta t}{\hat{r}_i}^{(n+1/2)},\quad \hat{r}_i^{\left(n+1/2\right)}=\frac{\hat{r}_i^{(n+1)}+\hat{r}_i^{(n)}}{2},
  \end{equation}
 then, we call $\boldsymbol{\delta} _dE_i^{(n+1)}$ discrete variational derivative of the energy $E^{(n)}$.
\end{definition}
\subsubsection{The Willmore flow}
We discretize the Wilmore flow \eqref{will_hen} with tangential velocity as follows:
\begin{equation}
  \begin{split}
  \label{sc}
  &\frac{\boldsymbol{X}_i^{(n+1)}-\boldsymbol{X}^{(n)}_i}{\Delta t}=-\boldsymbol{\delta}_d B_i^{(n+1)}+w_i^{(n)}\boldsymbol{T}_i^{(n)}+\gamma^{(n+1)}\boldsymbol{\delta}_d B_i^{(n+1)},\\
  &\qquad \qquad\qquad\qquad\qquad\qquad\qquad\qquad\:\:\:\:\:\:i=1,\dots,N,\quad n=0,1,\dots,
  \end{split}
\end{equation}
where $\boldsymbol{\delta} _dB_i^{(n+1)}$ represents the bending energy's discrete variational derivative, which is obtained according to definition \ref{definition:henbun}.
Detailed derivation and the concrete form of $\boldsymbol{\delta} _dB_i^{(n+1)}$ is presented in the Appendix. In \eqref{sc}, a modification term $\gamma^{(n+1)}\boldsymbol{\delta}_d B_i^{(n+1)}$ is added to preserve the bending energy's discrete dissipation. This is because, in a continuous system, the variations are parallel to the normal vector $\boldsymbol{N }(u,t)$ and orthogonal to the tangent vector $\boldsymbol{T}(u,t)$, but the discrete variation $\boldsymbol{\delta}_d B_i^{(n+1)}$ and the tangent vector $\boldsymbol{T}_i$ obtained from the formula (\ref{T}) are not necessarily orthogonal. Thus, the modification term $\gamma^{(n+1)}\boldsymbol{\delta}_d B_i^{(n+1)}$ is added to make the scheme structure-preserving. 

The Lagrange multiplier $\gamma^{(n+1)}$ is then determined so that the scheme conserves the equation's energy structure.

For simplicity, we represent the discrete inner product as 
\begin{equation}
  \label{naiseki}
(\boldsymbol{X},\boldsymbol{Y})_{d,n}=\sum_{i=1}^N\boldsymbol{X}_i\cdot\boldsymbol{Y}_i\hat{r}_i^{\left(n+1/2\right)},
\end{equation}
where 
\begin{equation}
  \boldsymbol{X}=(\boldsymbol{X}_1,\boldsymbol{X}_2,\dots,\boldsymbol{X}_N),\quad   \boldsymbol{Y}=(\boldsymbol{Y}_1,\boldsymbol{Y}_2,\dots,\boldsymbol{Y}_N)\in \mathbb{R}^{2\times N} 
\end{equation}
represent polygonal curves with vertices $\boldsymbol{X}_i$ and $\boldsymbol{Y}_i$.
By the definition of the discrete variational derivative, we obtain
\begin{equation}
  \begin{split}
    &\frac{B_d^{(n+1)}-B_d^{(n)}}{\Delta t}=\left( \boldsymbol{\delta}_d B^{(n+1)},\frac{\boldsymbol{X}^{(n+1)}-\boldsymbol{X}^{(n)}}{{\Delta t}}\right)_{d,n}\\
    &=\left( \boldsymbol{\delta}_d B^{(n+1)},- \boldsymbol{\delta}_d B^{(n+1)}+w^{(n)}\boldsymbol{T}^{(n)}+\gamma^{(n+1)}\boldsymbol{\delta}_d B^{(n+1)}\right)_{d,n}\\
    &=-\left( \boldsymbol{\delta}_d B^{(n+1)},\boldsymbol{\delta}_d B^{(n+1)}\right)_{d,n}+\left(\boldsymbol{\delta}_d B^{(n+1)},w^{(n)}\boldsymbol{T}^{(n)}+\gamma^{(n+1)}\boldsymbol{\delta}_d B^{(n+1)}\right)_{d,n}.
  \end{split}
\end{equation}
Therefore, determining $\gamma^{(n+1)}$ by 
\begin{equation}
  \left(\boldsymbol{\delta}_d B^{(n+1)},w^{(n)}\boldsymbol{T}^{(n)}+\gamma^{(n+1)}\boldsymbol{\delta}_d B^{(n+1)}\right)_{d,n}=0,
\end{equation}
we obtain
\begin{equation}
  \label{B-B}
  \frac{B^{(n+1)}-B^{(n)}}{\Delta t}=-\left( \boldsymbol{\delta}_d B^{(n+1)},\boldsymbol{\delta}_d B^{(n+1)}\right)_{d,n}\leq0.
\end{equation}
The dissipative nature of the discrete bending energy is verified from \eqref{B-B}. The resulting $\gamma^{(n+1)}$ is anticipated to be quite small. This is because the directions of the discrete variational derivative and the tangent vector are considered almost orthogonal. Thus, it is reasonable to introduce a new variable $\gamma$ in the above way.

From the above, we propose a structure-preserving scheme for the Wilmore flow as follows: for $n=0,1,\dots,$ 
\begin{equation}
  \begin{split}
    \label{will_sk}
    &\frac{\boldsymbol{X}_i^{(n+1)}-\boldsymbol{X}^{(n)}_i}{\Delta t}=- \boldsymbol{\delta}_d B_i^{(n+1)}+w_i^{(n)}\boldsymbol{T}_i^{(n)}-\gamma^{(n+1)} \boldsymbol{\delta}_d B_i^{(n+1)},\quad i=1,2,...N,\\
    & \left(\boldsymbol{\delta}_d B^{(n+1)},w^{(n)}\boldsymbol{T}^{(n)}+\gamma^{(n+1)}\boldsymbol{\delta}_d B^{(n+1)}\right)_{d,n}=0.
  \end{split}
\end{equation}
Note that the above scheme are nonlinear equations for the unknowns $\boldsymbol{X}_i^{(n+1)}(i=1,\dots,N)$ and $\gamma^{(n+1)}$.

\subsubsection{The Helfrich flow}
The Helfrich flow \eqref{hell_hen} is discretized into the following form, and we introduce tangential velocities:
\begin{equation}
  \begin{split}
    \label{hel_d}
&\frac{\boldsymbol{X}_i^{(n+1)}-\boldsymbol{X}_i^{(n)}}{\Delta t}\\
&=-\boldsymbol{\delta} _dB_i^{(n+1)}+\lambda ^{(n+1)}\boldsymbol{\delta} _dL_i^{(n+1)}+\mu ^{(n+1)}\boldsymbol{\delta}_d A_i ^{(n+1)}+w_i^{(n)}\boldsymbol{T}_i^{(n)}+\gamma^{(n+1)}\:\boldsymbol{\delta}_d B_i^{(n+1)},\\
&\qquad \qquad\qquad\qquad\qquad\qquad\qquad\qquad\qquad\:\:\:\:\:\:i=1,\dots,N,\quad n=0,1,\dots,
\end{split}
\end{equation}
where $\boldsymbol{\delta} _dB_i^{(n+1)}$, $\boldsymbol{\delta} _d L_i^{(n+1)}$, and $\boldsymbol{\delta} _dA_i^{(n+1)}$ are, respectively, the discrete variational derivatives of the bending energy $B_d$, the length $L_d$, and the enclosed area $A_d$. The detailed deformations and their respective concrete forms are presented in the Appendix. Tangential velocities and their modification terms are also added, as in the Willmore flow's scheme. 

The Lagrange multipliers $\lambda^{(n+1)},\mu^{(n+1)}$, and $\gamma^{(n+1)}$  are determined so that the scheme reproduces the equation's energy structure. We first determine the Lagrange multipliers $\lambda^{(n+1)}$ and $\mu^{(n+1)}$ for the conservation of the length $L_d$ and the enclosed area $A_d$. By the definition of the discrete variational derivatives, the time-differences of the length and the area are expressed as 
\begin{align}
  \label{lambda_d}
&\frac{L_d^{(n+1)}-L_d^{(n)}}{\Delta t}=\left(\boldsymbol{\delta} _dL^{(n+1)},\frac{\boldsymbol{X}^{(n+1)}-\boldsymbol{X}^{(n)}}{\Delta t}\right)
_{d,n},\\
\label{mu_d}
&\frac{A_d^{(n+1)}-A_d^{(n)}}{\Delta t}=\left(\boldsymbol{\delta} _dA^{(n+1)},\frac{\boldsymbol{X}^{(n+1)}-\boldsymbol{X}^{(n)}}{\Delta t}\right)
_{d,n}.
\end{align}
By substituting \eqref{hel_d} for time-difference of vertices $\boldsymbol{X}$ in  \eqref{lambda_d} and \eqref{mu_d}, we have
\footnotesize
\begin{multline}
  \label{lam_d2} 
  \frac{L_d^{(n+1)}-L_d^{(n)}}{\Delta t}\\
  =(\boldsymbol{\delta} _dL^{(n+1)},-\boldsymbol{\delta} _dB_i^{(n+1)}+\lambda ^{(n+1)}\boldsymbol{\delta} _dL^{(n+1)}+\mu ^{(n+1)}\boldsymbol{\delta}_d A ^{(n+1)}+w^{(n)}\boldsymbol{T}^{(n)}+\gamma^{(n+1)}\:\boldsymbol{\delta}_d B^{(n+1)})_{d,n},
\end{multline}
\begin{multline}
  \label{mu_d2} 
  \frac{A_d^{(n+1)}-A_d^{(n)}}{\Delta t}\\
 = (\boldsymbol{\delta} _dA^{(n+1)},-\boldsymbol{\delta} _dB_i^{(n+1)}+\lambda ^{(n+1)}\boldsymbol{\delta} _dL^{(n+1)}+\mu ^{(n+1)}\boldsymbol{\delta}_d A ^{(n+1)}+w^{(n)}\boldsymbol{T}^{(n)}+\gamma^{(n+1)}\:\boldsymbol{\delta}_d B^{(n+1)})_{d,n}.
\end{multline}
\normalsize
We can reproduce the conservation of the length and the area by letting \eqref{lam_d2} and \eqref{mu_d2} equal to 0.

Next, the Lagrange multiplier $\gamma^{(n+1)}$ for the bending energy's dissipative nature is determined. By substituting \eqref{hel_d} for the time-difference of vertices $\boldsymbol{X}$ and transforming, we obtain 
\footnotesize
\begin{equation}
  \begin{split}
    \label{gamma_d}
    &\frac{B_d^{(n+1)}-B_d^{(n)}}{\Delta t}=\left( \boldsymbol{\delta}_d B^{(n+1)},\frac{\boldsymbol{X}^{(n+1)}-\boldsymbol{X}^{(n)}}{{\Delta t}}\right)_{d,n}\\
    =&\left( \boldsymbol{\delta}_d B^{(n+1)},-\boldsymbol{\delta} _dB_i^{(n+1)}+\lambda ^{(n+1)}\boldsymbol{\delta} _dL_i^{(n+1)}+\mu ^{(n+1)}\boldsymbol{\delta}_d A_i ^{(n+1)}+w_i^{(n)}\boldsymbol{T}_i^{(n)}+\gamma^{(n+1)}\:\boldsymbol{\delta}_d B_i^{(n+1)}\right)_{d,n}\\
    =&\left( \boldsymbol{\delta}_d B^{(n+1)},-\boldsymbol{\delta} _dB_i^{(n+1)}+\lambda ^{(n+1)}\boldsymbol{\delta} _dL_i^{(n+1)}+\mu ^{(n+1)}\boldsymbol{\delta}_d A_i ^{(n+1)}+w_i^{(n)}\boldsymbol{T}_i^{(n)}+\gamma^{(n+1)}\:\boldsymbol{\delta}_d B_i^{(n+1)}\right)_{d,n}\\
    &+\:\mathrm{det}\left(\begin{matrix}(\boldsymbol{\delta} B_d,\boldsymbol{\delta} B_d)&(\boldsymbol{\delta} _dB,\boldsymbol{\delta} _dL)&(\boldsymbol{\delta} _dB,\boldsymbol{\delta} _dA)\\
      (\boldsymbol{\delta} _dL,\boldsymbol{\delta} _dB)&(\boldsymbol{\delta} _dL,\boldsymbol{\delta} _dL)&(\boldsymbol{\delta} _dL,\boldsymbol{\delta} _dA)\\
      (\boldsymbol{\delta} _dA,\boldsymbol{\delta} _dB)&(\boldsymbol{\delta} _dA,\boldsymbol{\delta} _dL)&(\boldsymbol{\delta} _dA,\boldsymbol{\delta} _dA)
    \end{matrix} \right)
  /\mathrm{det}\left(\begin{matrix}
      (\boldsymbol{\delta} _dL,\boldsymbol{\delta} _dL)&(\boldsymbol{\delta} _dL,\boldsymbol{\delta} _dA)\\
      (\boldsymbol{\delta} _dA,\boldsymbol{\delta} _dL)&(\boldsymbol{\delta} _dA,\boldsymbol{\delta} _dA)
    \end{matrix}\right)\\
    &-\:\mathrm{det}\left(\begin{matrix}(\boldsymbol{\delta} B_d,\boldsymbol{\delta} B_d)&(\boldsymbol{\delta} _dB,\boldsymbol{\delta} _dL)&(\boldsymbol{\delta} _dB,\boldsymbol{\delta} _dA)\\
      (\boldsymbol{\delta} _dL,\boldsymbol{\delta} _dB)&(\boldsymbol{\delta} _dL,\boldsymbol{\delta} _dL)&(\boldsymbol{\delta} _dL,\boldsymbol{\delta} _dA)\\
      (\boldsymbol{\delta} _dA,\boldsymbol{\delta} _dB)&(\boldsymbol{\delta} _dA,\boldsymbol{\delta} _dL)&(\boldsymbol{\delta} _dA,\boldsymbol{\delta} _dA)
    \end{matrix} \right)
  /\mathrm{det}\left(\begin{matrix}
      (\boldsymbol{\delta} _dL,\boldsymbol{\delta} _dL)&(\boldsymbol{\delta} _dL,\boldsymbol{\delta} _dA)\\
      (\boldsymbol{\delta} _dA,\boldsymbol{\delta} _dL)&(\boldsymbol{\delta} _dA,\boldsymbol{\delta} _dA)
    \end{matrix}\right).
\end{split}
\end{equation}
\normalsize
Here, for simplicity, the following substitutions are made.
\footnotesize
  \begin{equation}
    D(\boldsymbol{\delta}_d B,\boldsymbol{\delta} _dL,\boldsymbol{\delta} _dA)=
    \mathrm{det}\left(\begin{matrix}(\boldsymbol{\delta} B_d,\boldsymbol{\delta} B_d)&(\boldsymbol{\delta} _dB,\boldsymbol{\delta} _dL)&(\boldsymbol{\delta} _dB,\boldsymbol{\delta} _dA)\\
      (\boldsymbol{\delta} _dL,\boldsymbol{\delta} _dB)&(\boldsymbol{\delta} _dL,\boldsymbol{\delta} _dL)&(\boldsymbol{\delta} _dL,\boldsymbol{\delta} _dA)\\
      (\boldsymbol{\delta} _dA,\boldsymbol{\delta} _dB)&(\boldsymbol{\delta} _dA,\boldsymbol{\delta} _dL)&(\boldsymbol{\delta} _dA,\boldsymbol{\delta} _dA)
    \end{matrix} \right)
  /\mathrm{det}\left(\begin{matrix}
      (\boldsymbol{\delta} _dL,\boldsymbol{\delta} _dL)&(\boldsymbol{\delta} _dL,\boldsymbol{\delta} _dA)\\
      (\boldsymbol{\delta} _dA,\boldsymbol{\delta} _dL)&(\boldsymbol{\delta} _dA,\boldsymbol{\delta} _dA)
    \end{matrix}\right).
  \end{equation}
\normalsize
Determining the Lagrange multiplier $\gamma^{(n+1)}$ from 
\footnotesize
\begin{equation}
  \begin{split}
    \label{helf_gamma}
  &( \boldsymbol{\delta}_d B^{(n+1)},-\boldsymbol{\delta} _dB^{(n+1)}+\lambda ^{(n+1)}\boldsymbol{\delta} _dL^{(n+1)}+\mu ^{(n+1)}\boldsymbol{\delta}_d A^{(n+1)}
  +w^{(n)}\boldsymbol{T}^{(n)}+\gamma^{(n+1)}\:\boldsymbol{\delta}_d B^{(n+1)})_{d,n}\\
  &+ D(\boldsymbol{\delta}_d B,\boldsymbol{\delta} _dL,\boldsymbol{\delta} _dA)=0,
\end{split}
\end{equation}
\normalsize
we obtain the discrete bending energy dissipation.
From the above, we propose a structure-preserving scheme for the Helfrich flow as follows: for $n=0,1,\dots,$ 
\footnotesize
\begin{equation}
  \begin{split}
    \label{hell_ski}
    &\frac{\boldsymbol{X}_i^{(n+1)}-\boldsymbol{X}_i^{(n)}}{\Delta t}=-\boldsymbol{\delta} _dB_i^{(n+1)}+\lambda ^{(n+1)}\boldsymbol{\delta} _dL_i^{(n+1)}+\mu ^{(n+1)}\boldsymbol{\delta}_d A_i ^{(n+1)}+w_i^{(n)}\boldsymbol{T}_i^{(n)}+\gamma^{(n+1)}\:\boldsymbol{\delta}_d B_i^{(n+1)},\\
    & i=1,2,...N,\\
    &   (\boldsymbol{\delta} _dL^{(n+1)},-\boldsymbol{\delta} _dB^{(n+1)}+\lambda ^{(n+1)}\boldsymbol{\delta} _dL^{(n+1)}+\mu ^{(n+1)}\boldsymbol{\delta}_d A ^{(n+1)}+w^{(n)}\boldsymbol{T}^{(n)}+\gamma^{(n+1)}\:\boldsymbol{\delta}_d B^{(n+1)})_{d,n}=0,\\
    &(\boldsymbol{\delta} _dA^{(n+1)},-\boldsymbol{\delta} _dB^{(n+1)}+\lambda ^{(n+1)}\boldsymbol{\delta} _dL^{(n+1)}+\mu ^{(n+1)}\boldsymbol{\delta}_d A ^{(n+1)}+w^{(n)}\boldsymbol{T}^{(n)}+\gamma^{(n+1)}\:\boldsymbol{\delta}_d B^{(n+1)})_{d,n}=0,\\
    &( \boldsymbol{\delta}_d B^{(n+1)},-\boldsymbol{\delta} _dB^{(n+1)}+\lambda ^{(n+1)}\boldsymbol{\delta} _dL^{(n+1)}+\mu ^{(n+1)}\boldsymbol{\delta}_d A ^{(n+1)}+w^{(n)}\boldsymbol{T}^{(n)}+\gamma^{(n+1)}\:\boldsymbol{\delta}_d B^{(n+1)})_{d,n}\\
    &+D(\boldsymbol{\delta}_d B,\boldsymbol{\delta} _dL,\boldsymbol{\delta} _dA)=0.
 \end{split}
\end{equation}
\normalsize
Note that the above scheme are nonlinear equations for the unknowns $\boldsymbol{X}_i^{(n+1)}(i=1,\dots,N)$ , $\lambda^{(n+1)},\mu^{(n+1)}$, and $\gamma^{(n+1)}$.

\begin{remark} 
  By determining $\gamma^{(n+1)}$ in \eqref{helf_gamma}, the difference quotient of the bending energy can be expressed as 
  \footnotesize
  \begin{equation}
    \begin{split}
      &\frac{B_d^{(n+1)}-B_d^{(n)}}{\Delta t}\\
       &=-\:\mathrm{det}\left(\begin{matrix}(\boldsymbol{\delta} B_d,\boldsymbol{\delta} B_d)_{n,d}&(\boldsymbol{\delta} _dB,\boldsymbol{\delta} _dL)_{n,d}&(\boldsymbol{\delta} _dB,\boldsymbol{\delta} _dA)_{n,d}\\
        (\boldsymbol{\delta} _dL,\boldsymbol{\delta} _dB)_{n,d}&(\boldsymbol{\delta} _dL,\boldsymbol{\delta} _dL)_{n,d}&(\boldsymbol{\delta} _dL,\boldsymbol{\delta} _dA)_{n,d}\\
        (\boldsymbol{\delta} _dA,\boldsymbol{\delta} _dB)_{n,d}&(\boldsymbol{\delta} _dA,\boldsymbol{\delta} _dL)_{n,d}&(\boldsymbol{\delta} _dA,\boldsymbol{\delta} _dA)_{n,d}
      \end{matrix} \right)
    /\mathrm{det}\left(\begin{matrix}
        (\boldsymbol{\delta} _dL,\boldsymbol{\delta} _dL)_{n,d}&(\boldsymbol{\delta} _dL,\boldsymbol{\delta} _dA)_{n,d}\\
        (\boldsymbol{\delta} _dA,\boldsymbol{\delta} _dL)_{n,d}&(\boldsymbol{\delta} _dA,\boldsymbol{\delta} _dA)_{n,d}
      \end{matrix}\right),
  \end{split}
  \end{equation}
  \normalsize
  which corresponds to \eqref{dbdt}. Here, superscripts $(n+1)$ for $\boldsymbol{\delta} B_d$, $\boldsymbol{\delta} L_d$, and $\boldsymbol{\delta} A_d$ are omitted. Therefore, the derivation of dissipativity in the continuous system can be reproduced in the discrete system by determining the $\gamma$ by \eqref{helf_gamma}.
  Moreover, explicitly obtaining $\gamma^{(n+1)}$ from the expression\eqref{hell_ski}, we have
  \scriptsize
  \begin{equation}
    \label{yo_gamma}
    \gamma^{(n+1)}\mathrm{det}\left(\begin{matrix}(\boldsymbol{\delta}_d B,\boldsymbol{\delta}_d B)_{n,d}&(\boldsymbol{\delta} _dB,\boldsymbol{\delta} _dL)_{n,d}&(\boldsymbol{\delta} _dB,\boldsymbol{\delta} _dA)_{n,d}\\
      (\boldsymbol{\delta} _dL,\boldsymbol{\delta} _dB)_{n,d}&(\boldsymbol{\delta} _dL,\boldsymbol{\delta} _dL)_{n,d}&(\boldsymbol{\delta} _dL,\boldsymbol{\delta} _dA)_{n,d}\\
      (\boldsymbol{\delta} _dA,\boldsymbol{\delta} _dB)_{n,d}&(\boldsymbol{\delta} _dA,\boldsymbol{\delta} _dL)_{n,d}&(\boldsymbol{\delta} _dA,\boldsymbol{\delta} _dA)_{n,d}
    \end{matrix} \right)
  =
  -w^{(n)}\mathrm{det}
  \left(\begin{matrix}(\boldsymbol{\delta}_d B,\boldsymbol{T})_{n,d}&(\boldsymbol{\delta} _dB,\boldsymbol{T})_{n,d}&(\boldsymbol{\delta} _dB,\boldsymbol{T})_{n,d}\\
      (\boldsymbol{\delta} _dL,\boldsymbol{T})_{n,d}&(\boldsymbol{\delta} _dL,\boldsymbol{T})_{n,d}&(\boldsymbol{\delta} _dL,\boldsymbol{T})_{n,d}\\
      (\boldsymbol{\delta} _dA,\boldsymbol{T} )_{n,d}&(\boldsymbol{\delta} _dA,\boldsymbol{T})_{n,d}&(\boldsymbol{\delta} _dA,\boldsymbol{T} )_{n,d}
    \end{matrix} \right).
  \end{equation}
  \normalsize
   By \eqref{yo_gamma}, the value of $\gamma$ is considered to approach zero in the continuous limit.
\end{remark}




\section{Numerical experiments}
In this section, we provide numerical experiments of the schemes proposed in Section 3 for the Willmore and the Helfrich flows.
 Nonlinear equations \eqref{will_sk} and \eqref{hell_ski} are solved by Scipy.optimize.root, which is a common nonlinear solver in Python with a relative residual 
tolerance of $10^{-6}$.

\subsection{Willmore flow}
In this subsection, the findings of numerical computations by the scheme \eqref{will_sk} are presented.
\subsubsection{Examination. 1}
Numerical experiments are carried out for different values of the parameter $\alpha$ on tangential velocity to verify the effect of the value of $\alpha$.
The initial curve is set to  
\begin{equation}
  \begin{split}
    \label{kubire}
    x_1(t)&=0.5a_1(t),\qquad x_2(t)=0.54a_3(t),\qquad t\in[0,1],\\
    a_1(t)&=1.8\cos (2\pi t),\qquad a_2(t)=0.2+\sin(\pi t)\sin(6\pi t)\sin(2a_1(t)),\\
    a_3(t)&=0.5\sin(2\pi t)+\sin(a_1(t))+a_2(t)\sin(2\pi t),
  \end{split}
\end{equation}
which is taken from \cite{13}.
We set 30 vertices $\boldsymbol{X}_i^0=(X_i^0,Y_i^0)^\mathrm{T}$ to
\begin{equation}
  \label{syoki}
  X_i^0=x_1(i/N),\quad Y_i^0=x_2(i/N),\quad i=1,2,\dots,30.
\end{equation}

\begin{figure}[H]
  \begin{minipage}{0.5\hsize}
   \begin{center}
    \includegraphics[width=50mm]{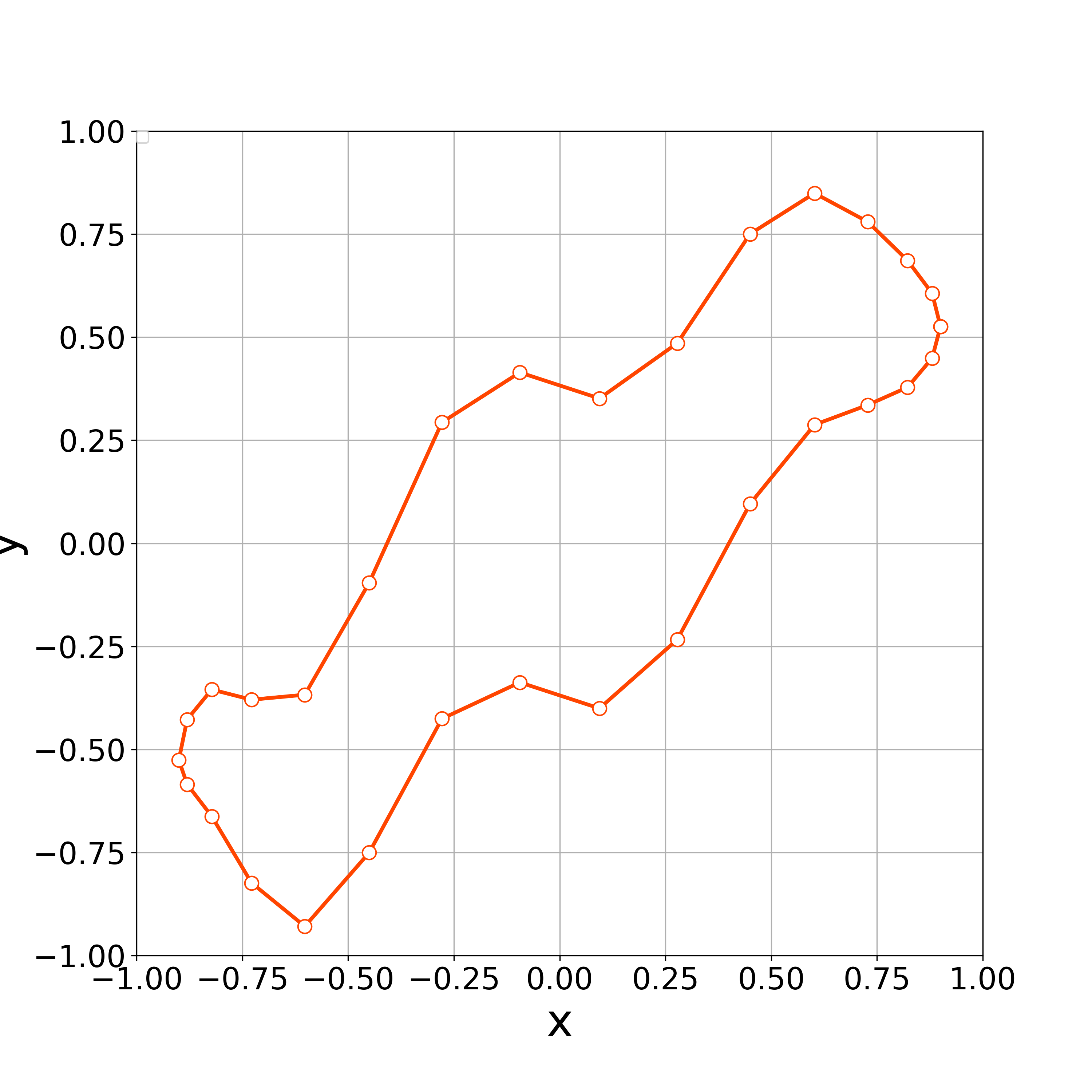}
   \end{center}
   \caption{Initial polygonal curve by \eqref{syoki}.}
   \label{fig:one}
  \end{minipage}
  \begin{minipage}{0.5\hsize}
   \begin{center}
    \includegraphics[width=50mm]{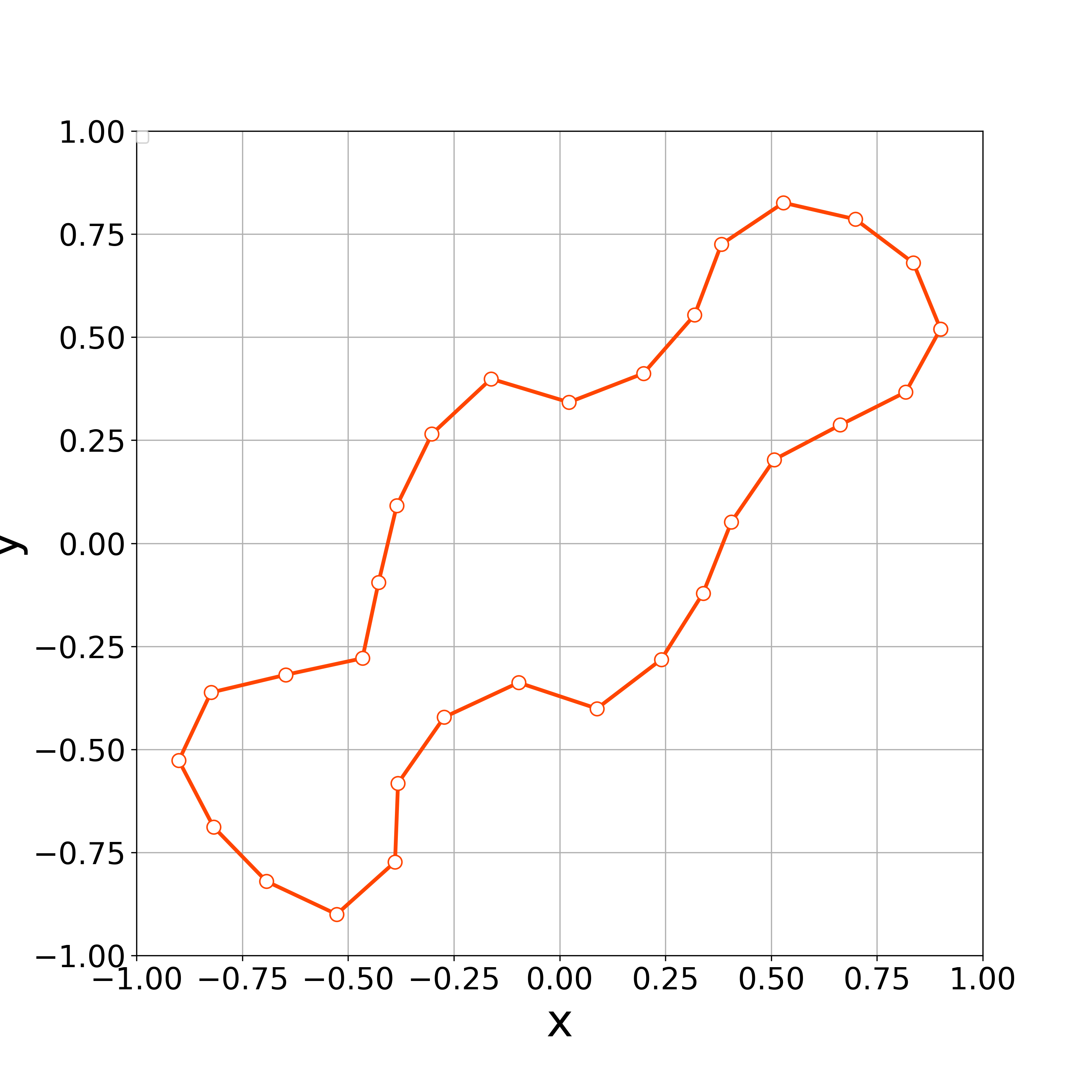}
   \end{center}
   \caption{Initial polygonal curve after relocation.}
   \label{fig:two}
  \end{minipage}
 \end{figure}

 The vertices' arrangement in Fig.\ \ref{fig:one} is not uniform, and numerical computations are prone to instability. Thus, we adjust them by 
 conducting a numerical computation with $\boldsymbol{\delta}_dB_i\;(i=1,2,...,N)$ equal to $\boldsymbol{0}$ in \eqref{will_sk}. The vertices were rearranged under the following conditions
 \begin{equation}
  \alpha=5,\quad\Delta t=10^{-4},\quad \text{up to}\:\; T=0.1.
 \end{equation}
Fig.\ \ref{fig:two} shows a initial polygonal curve after relocation.
Numerical experiments are carried out with the initial curve illustrated in Fig.\ \ref{fig:two}, time step size  $\Delta t=10^{-4}$ and the constant $c_0=2$. We set $\alpha=$0,10,50, and 200.

\begin{figure}[H]
  \begin{minipage}{0.5\hsize}
   \begin{center}
    \includegraphics[width=55mm]{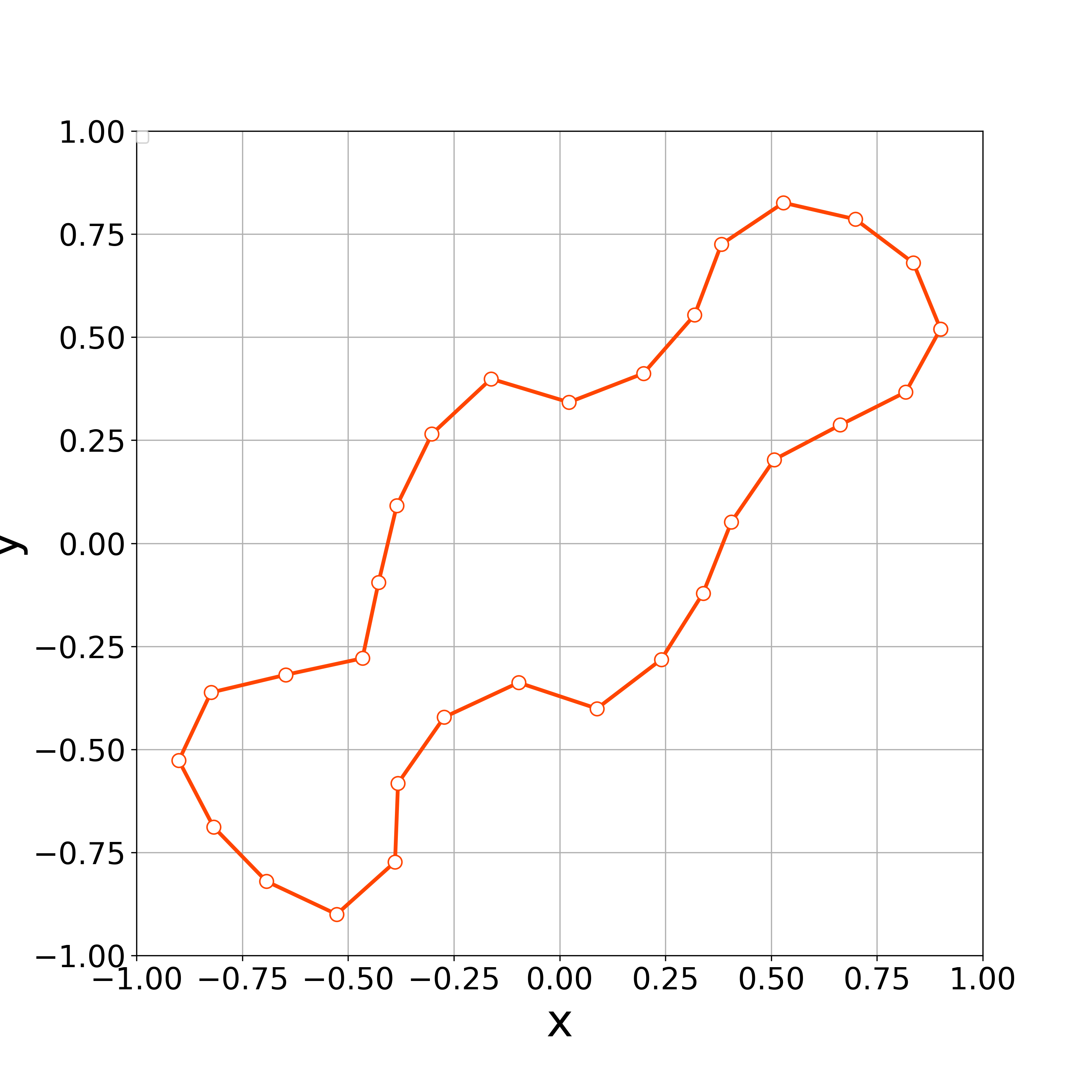}
   \end{center}
   \subcaption{$n=0.$}

  \end{minipage} 
  \begin{minipage}{0.5\hsize}
   \begin{center}
    \includegraphics[width=55mm]{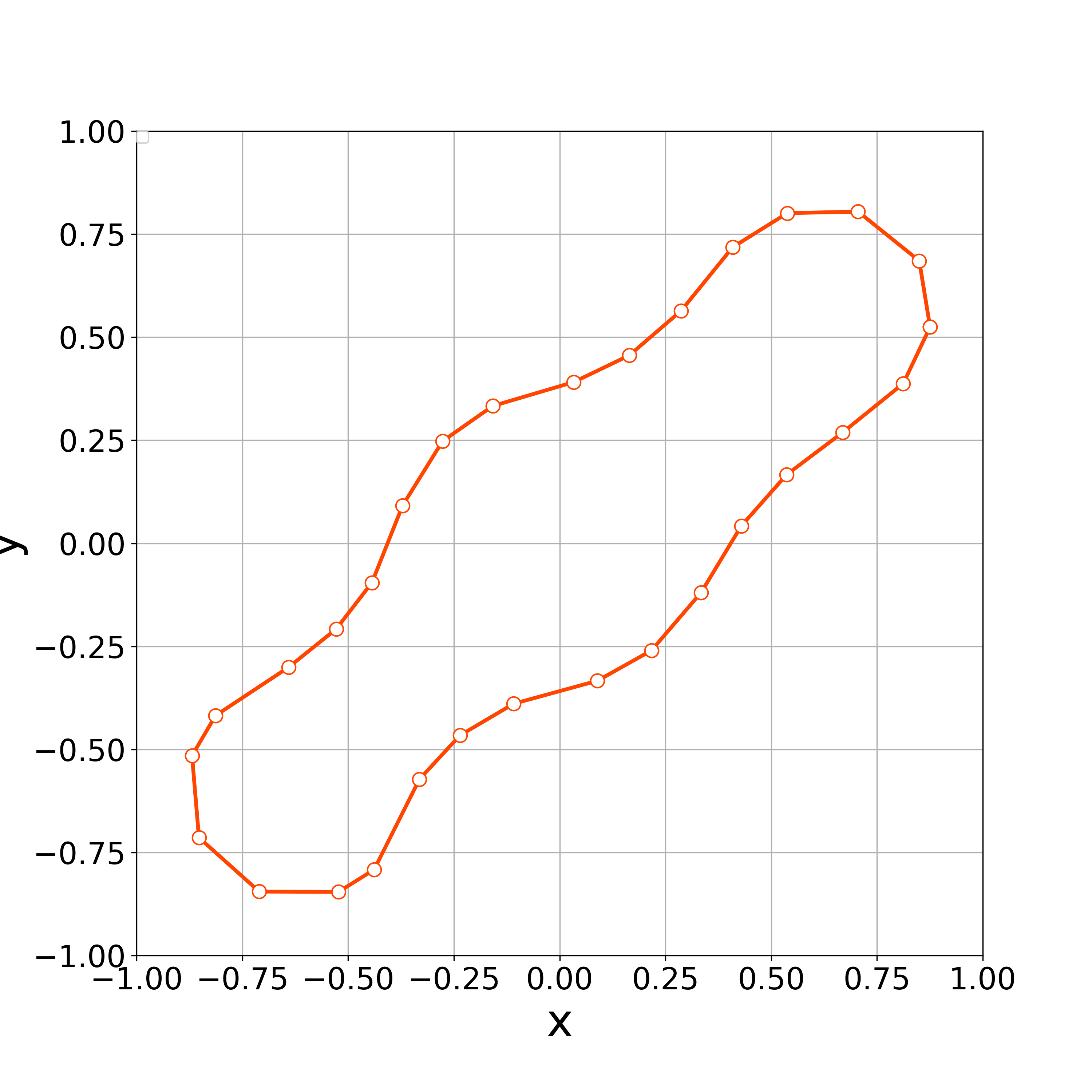}
   \end{center}
   \subcaption{$n=5.$}

  \end{minipage} 

  \begin{minipage}{0.5\hsize}
   \begin{center}
    \includegraphics[width=55mm]{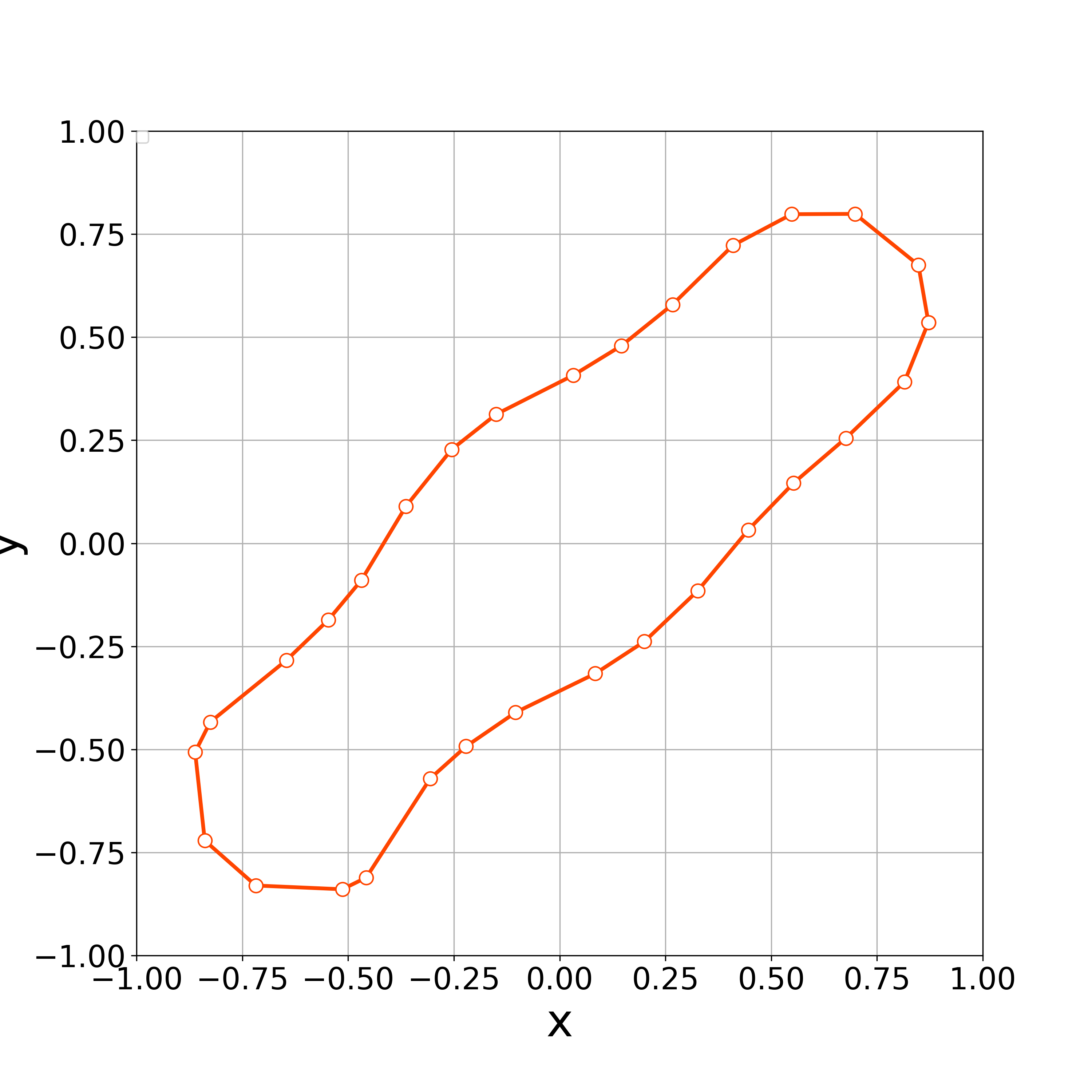}
   \end{center}
   \subcaption{$n=10$}

  \end{minipage}
  \begin{minipage}{0.5\hsize}
   \begin{center}
    \includegraphics[width=55mm]{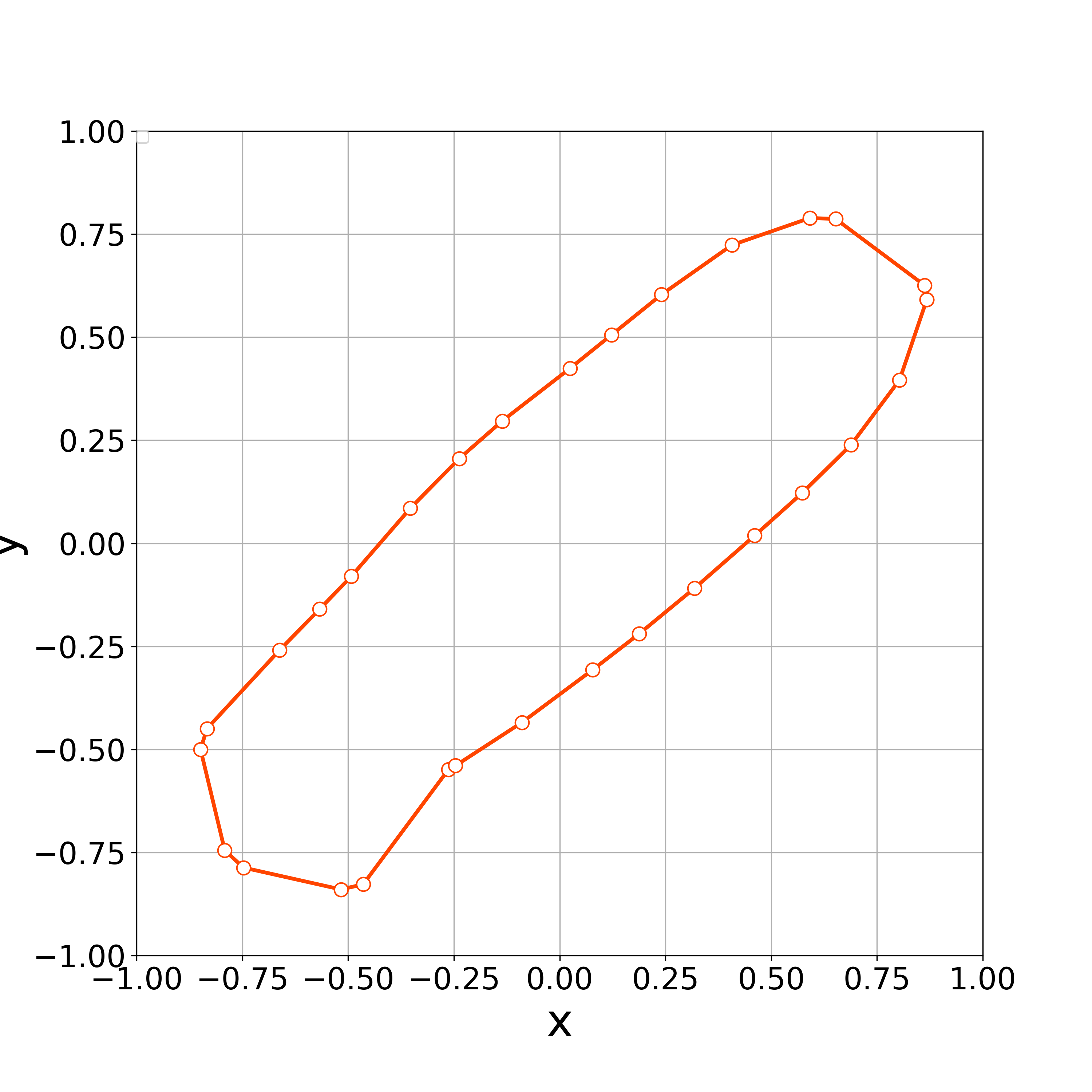}
   \end{center}
   \subcaption{$n=22$}

  \end{minipage}
  \caption{The $n$th step of the polygonal curve without tangential velocity.}
  \label{fig:38}
 \end{figure}

When tangential velocity is absent, i.e., $\alpha=0$, the concentration of vertices occurs at the step 22. Fig.\  \ref{fig:38} demonstrates that the vertices are concentrated at $(x,y)=(-0.25,-0.5),(0.8,0.6)$.

\begin{figure}[H]
  \begin{minipage}{0.5\hsize}
   \begin{center}
    \includegraphics[width=55mm]{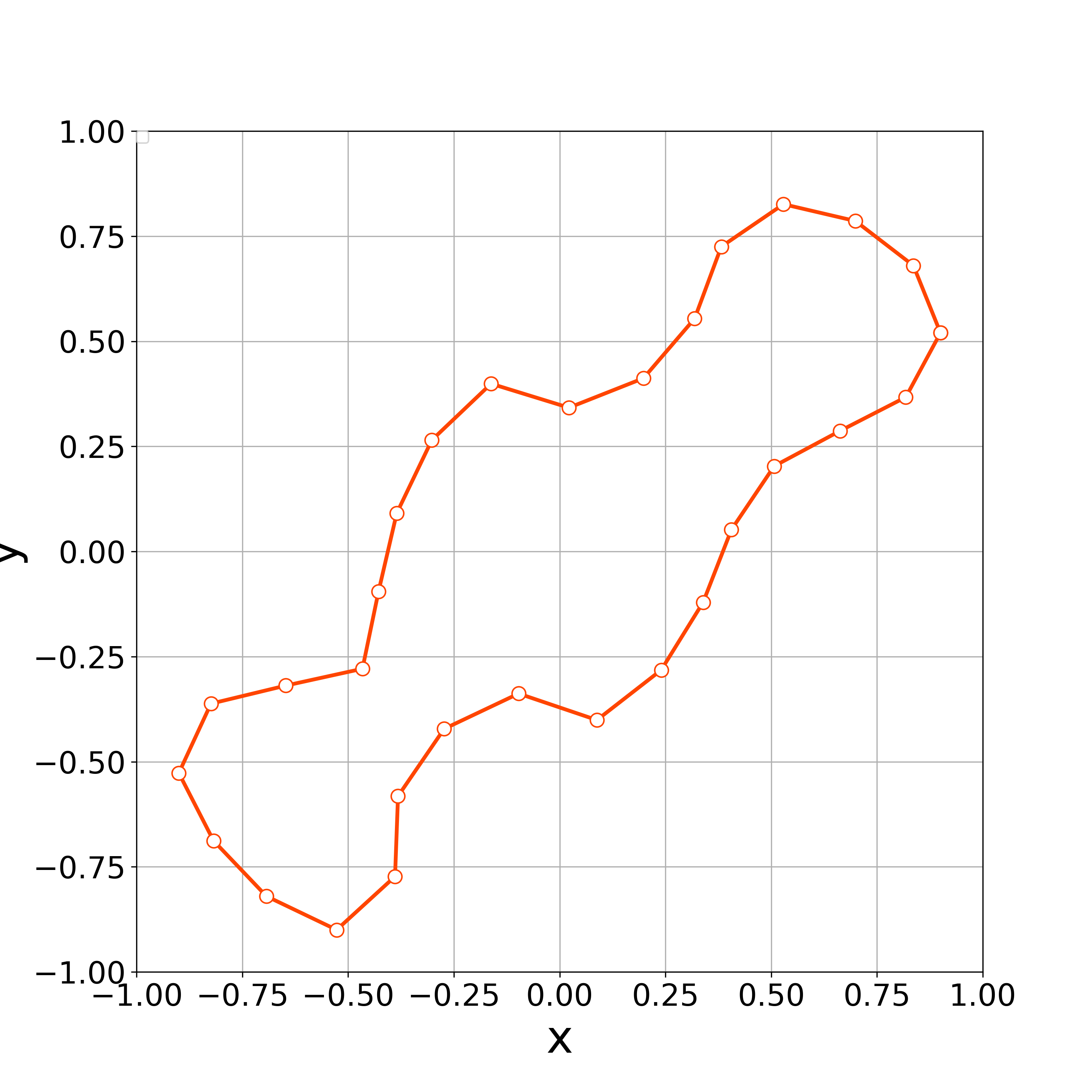}
   \end{center}
   \subcaption{$n=0.$}
   \label{fig:1}
  \end{minipage} 
  \begin{minipage}{0.5\hsize}
   \begin{center}
    \includegraphics[width=55mm]{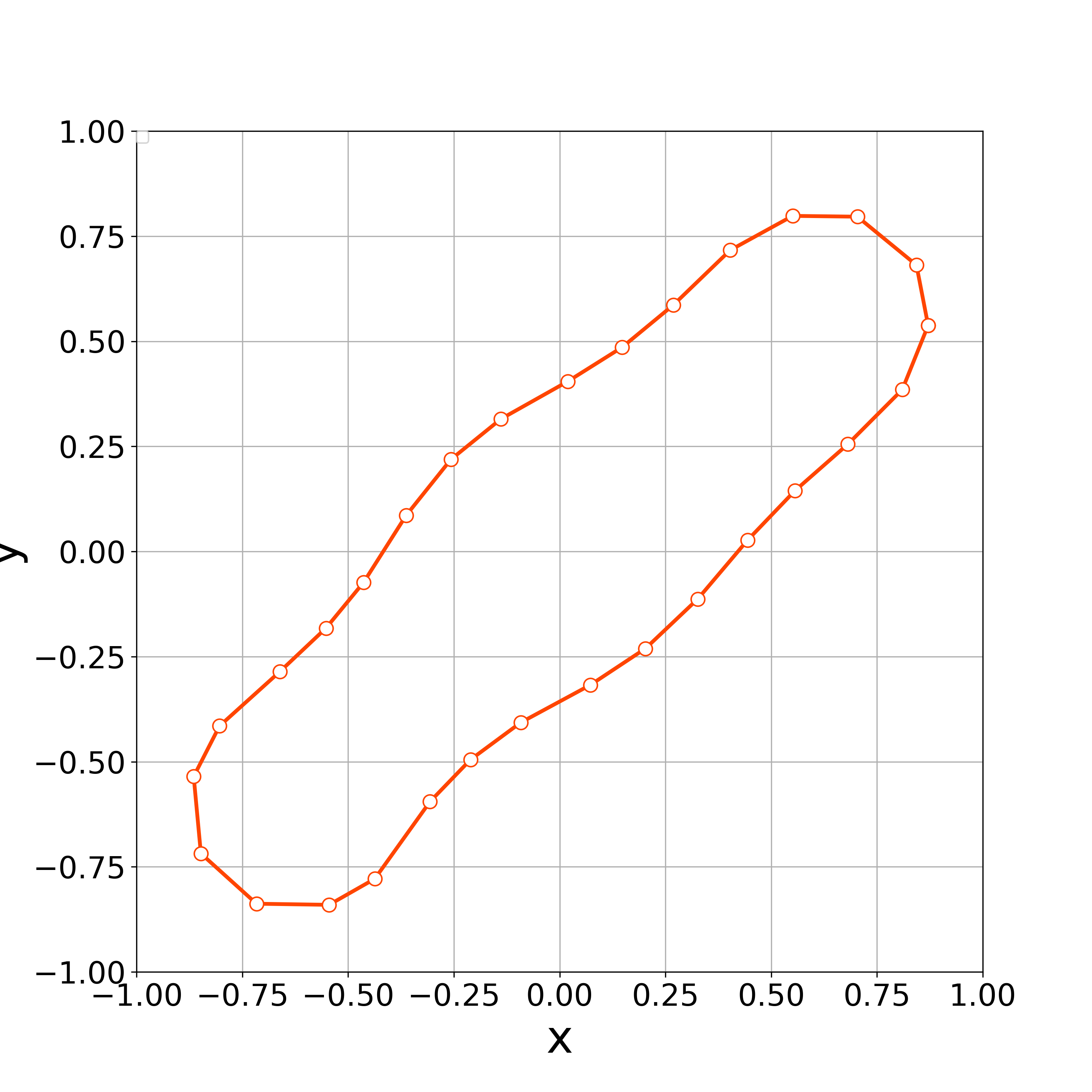}
   \end{center}
   \subcaption{$n=10.$}
   \label{fig:2}
  \end{minipage} 

  \begin{minipage}{0.5\hsize}
   \begin{center}
    \includegraphics[width=55mm]{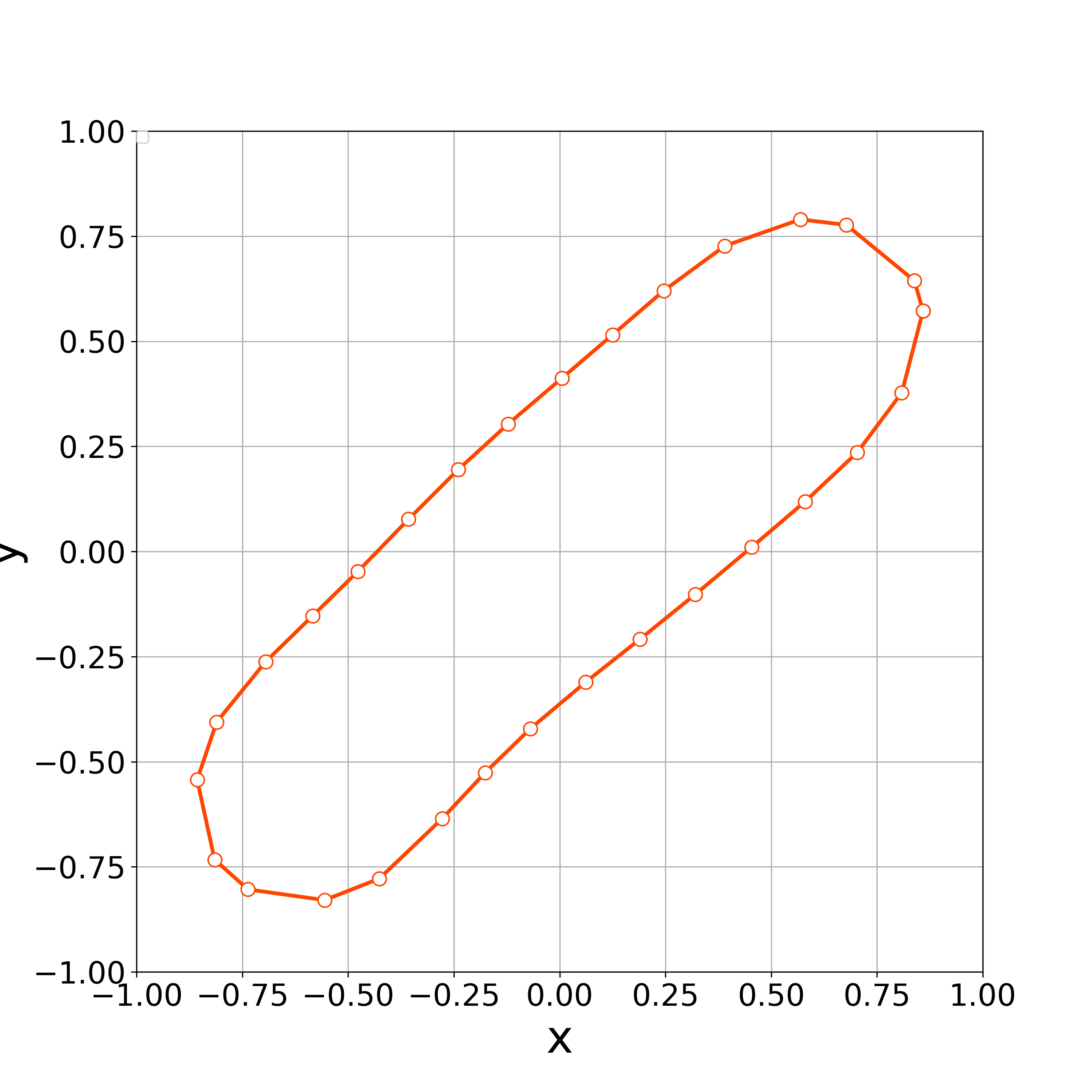}
   \end{center}
   \subcaption{$n=20$}
   \label{fig:3}
  \end{minipage}
  \begin{minipage}{0.5\hsize}
   \begin{center}
    \includegraphics[width=55mm]{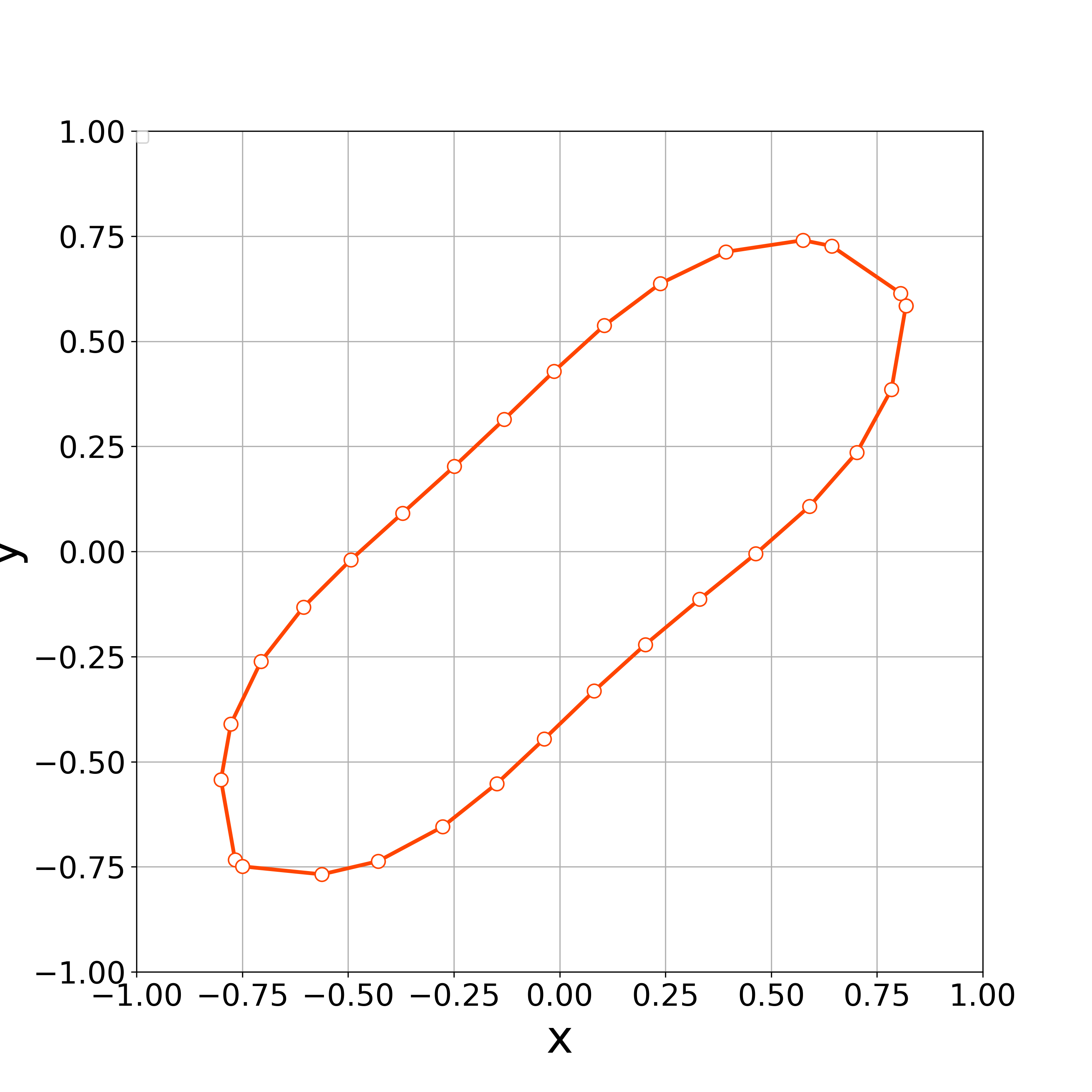}
   \end{center}
   \subcaption{$n=30$}

  \end{minipage}
  \caption{The $n$th step of the polygonal curve with $\alpha=10$.}
  \label{fig:4}
 \end{figure}

When $\alpha=10$, the concentration of vertices occurs the step 30. Fig.\  \ref{fig:4} demonstrates that the vertices are concentrated at $(x,y)=(-0.75,-0.75),(0.8,0.6)$. This result suggests that the value of $\alpha$ is too small to get stable numerical results.

\begin{figure}[H]
  \begin{minipage}{0.5\hsize}
   \begin{center}
    \includegraphics[width=55mm]{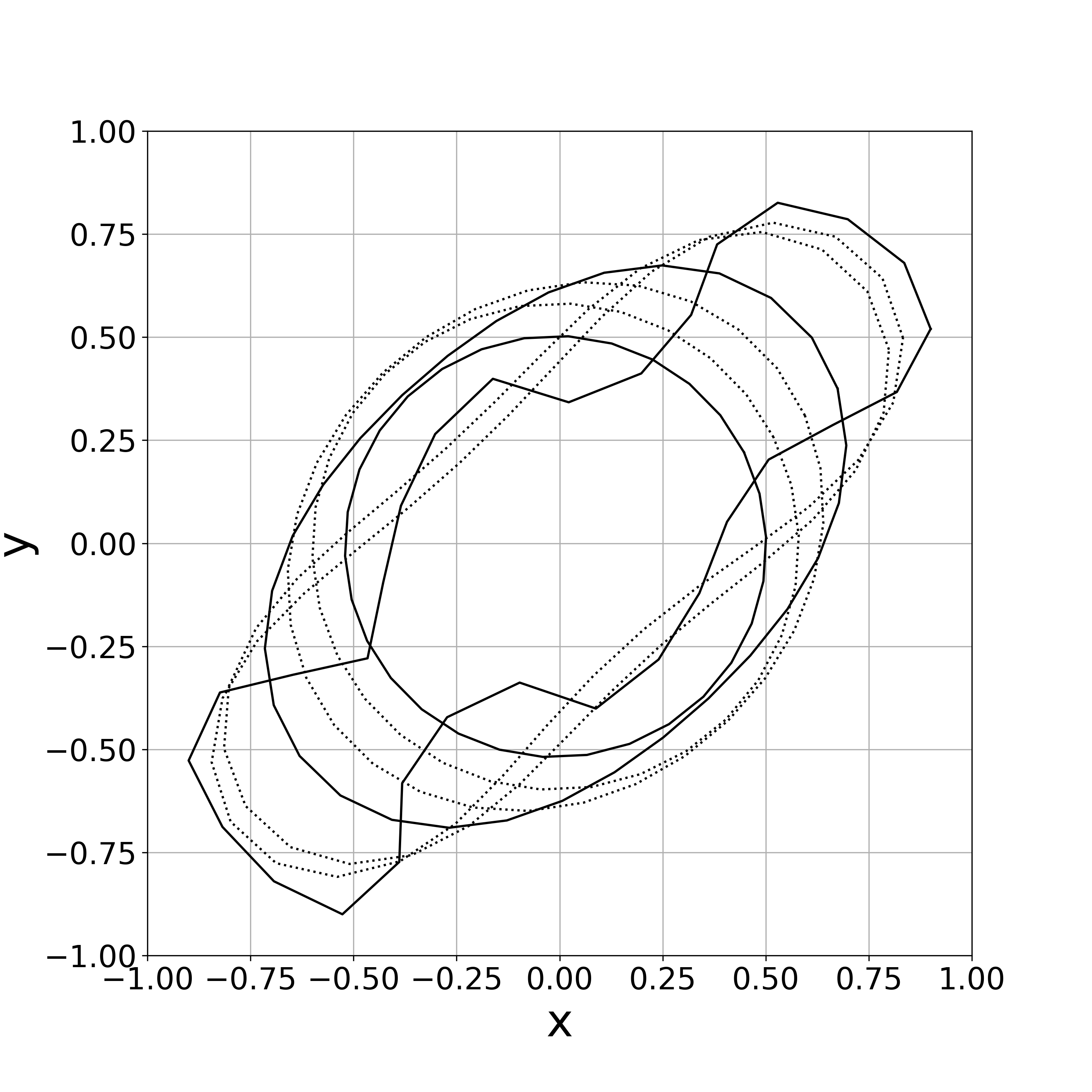}
   \end{center}
   \caption{Numerical solutions with $\alpha=50$\\ (initial polygonal curve and the numerical \\solutions at $n=10$, 50, 100, 300, 500, 1000\\ and 3000th steps).}
   \label{fig:5}
  \end{minipage}
  \begin{minipage}{0.5\hsize}
   \begin{center}
    \includegraphics[width=55mm]{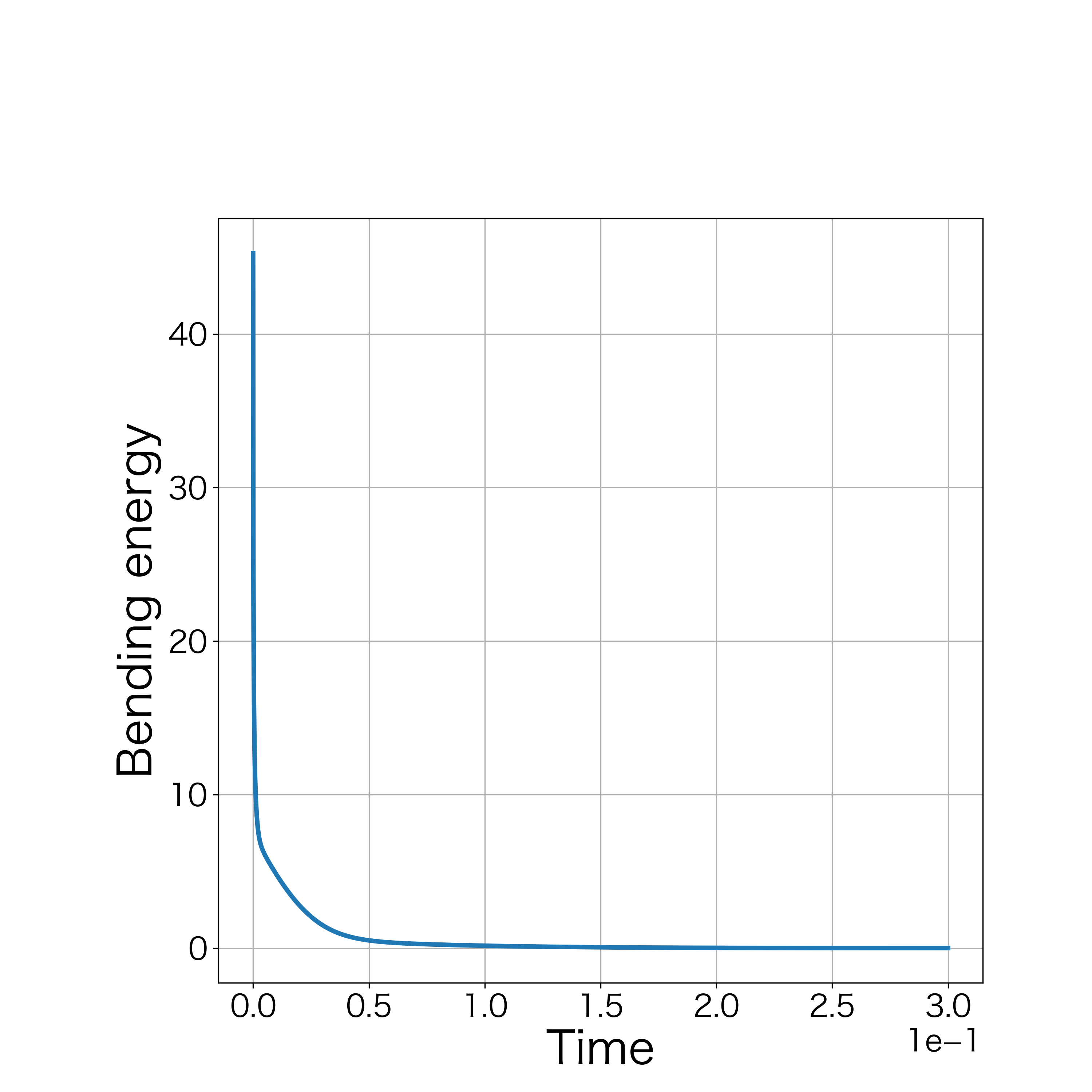}
   \end{center}
   \caption{The time evolution of the bending energy \\($\alpha=50$).}
   \label{fig:6}
  \end{minipage}
 \end{figure}

\begin{figure}[H]
   \begin{center}
    \includegraphics[width=55mm]{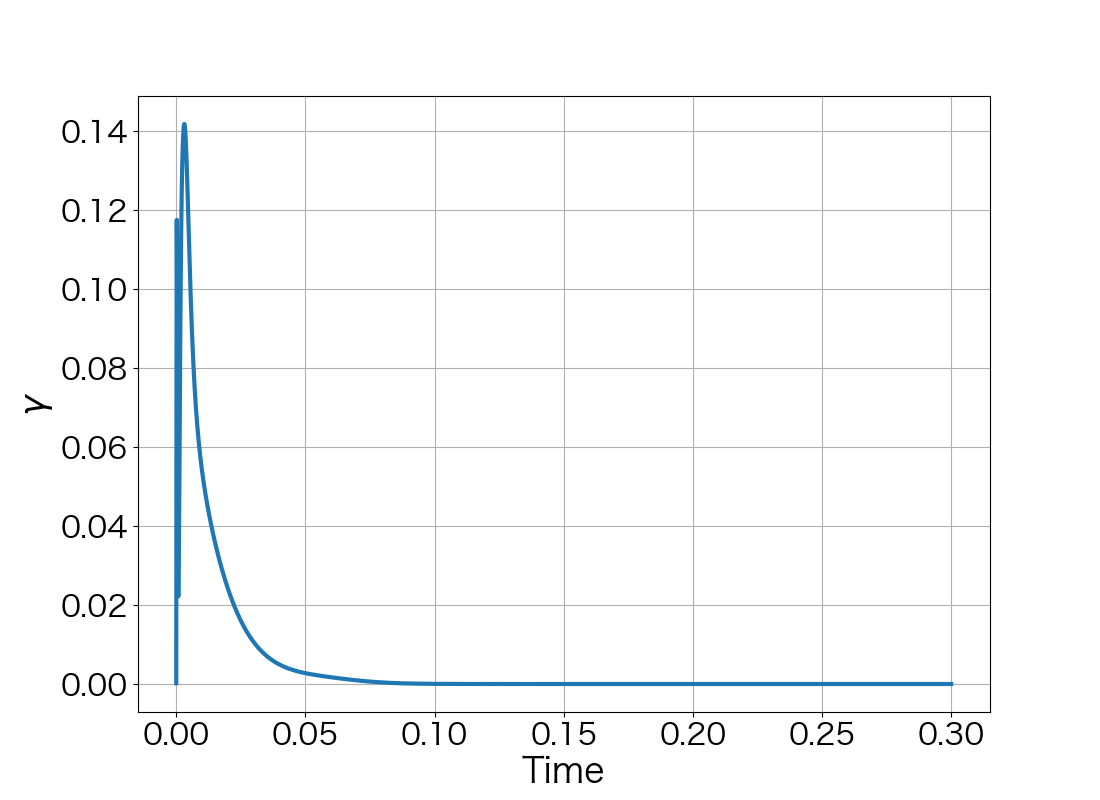}
   \end{center}
   \caption{Value of $\gamma$ at each step ($\alpha=50$).}
   \label{fig:will_gamma}
 \end{figure}

 When $\alpha=50$, the numerical computation demonstrates that the nonlinear solver converges up to the 3000th step. Fig.\ \ref{fig:5} demonstrates the initial polygonal curve and the numerical solutions at $n=10$, 50, 100, 300, 500, 1000, and 3000 steps, and we can see that the solution is approaching a regular circle with a radius $1/c_0=0.5$. Fig.\ \ref{fig:6} demonstrates the time evolution of the bending energy, and we can see the discrete version of energy dissipative nature. Furthermore, the value of $\gamma$ at each step approached 0 as the time step progressed (Fig.\ \ref{fig:will_gamma}).


\begin{figure}[H]
  \begin{minipage}{0.5\hsize}
   \begin{center}
    \includegraphics[width=55mm]{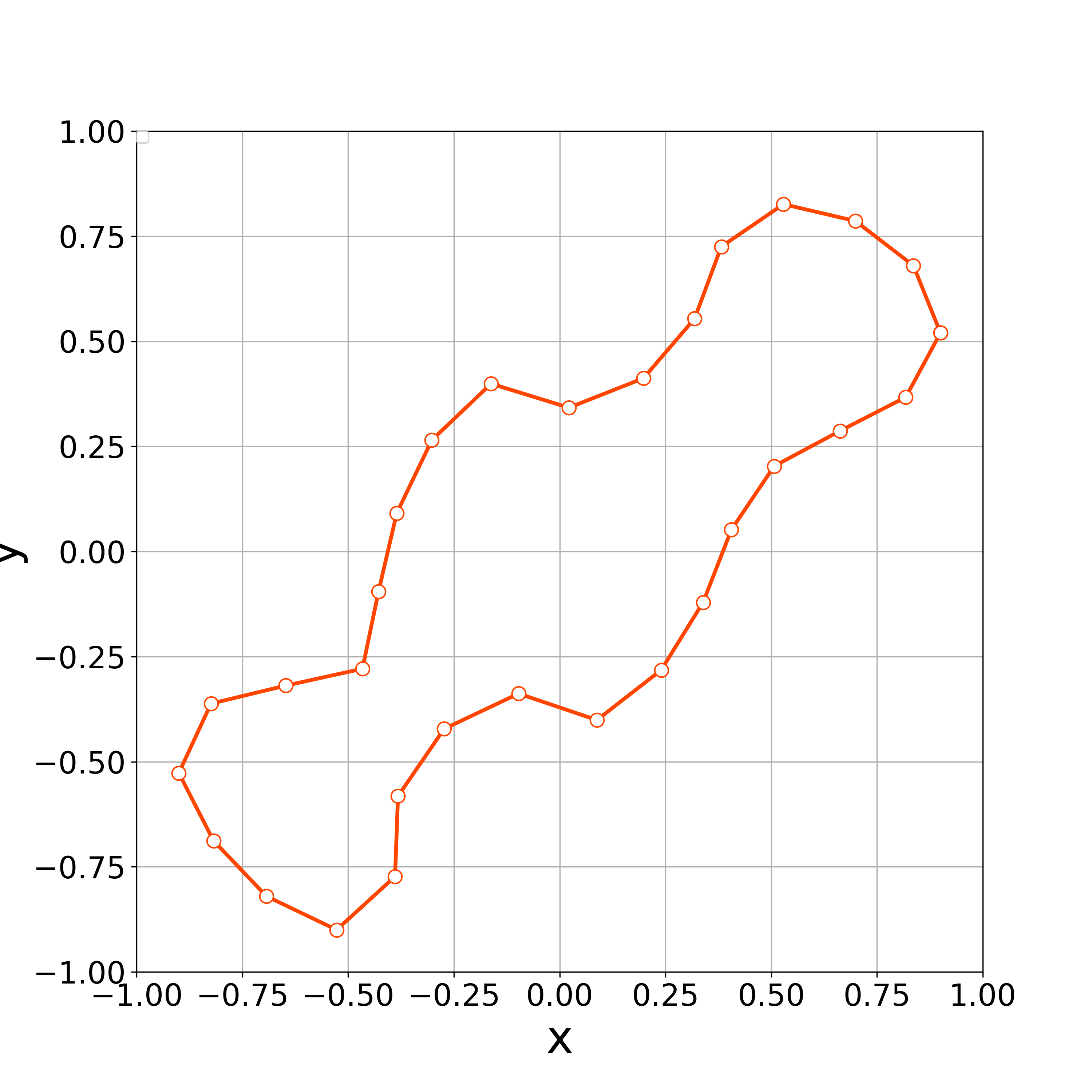}
   \end{center}
   \subcaption{$n=0.$}

  \end{minipage} 
  \begin{minipage}{0.5\hsize}
   \begin{center}
    \includegraphics[width=55mm]{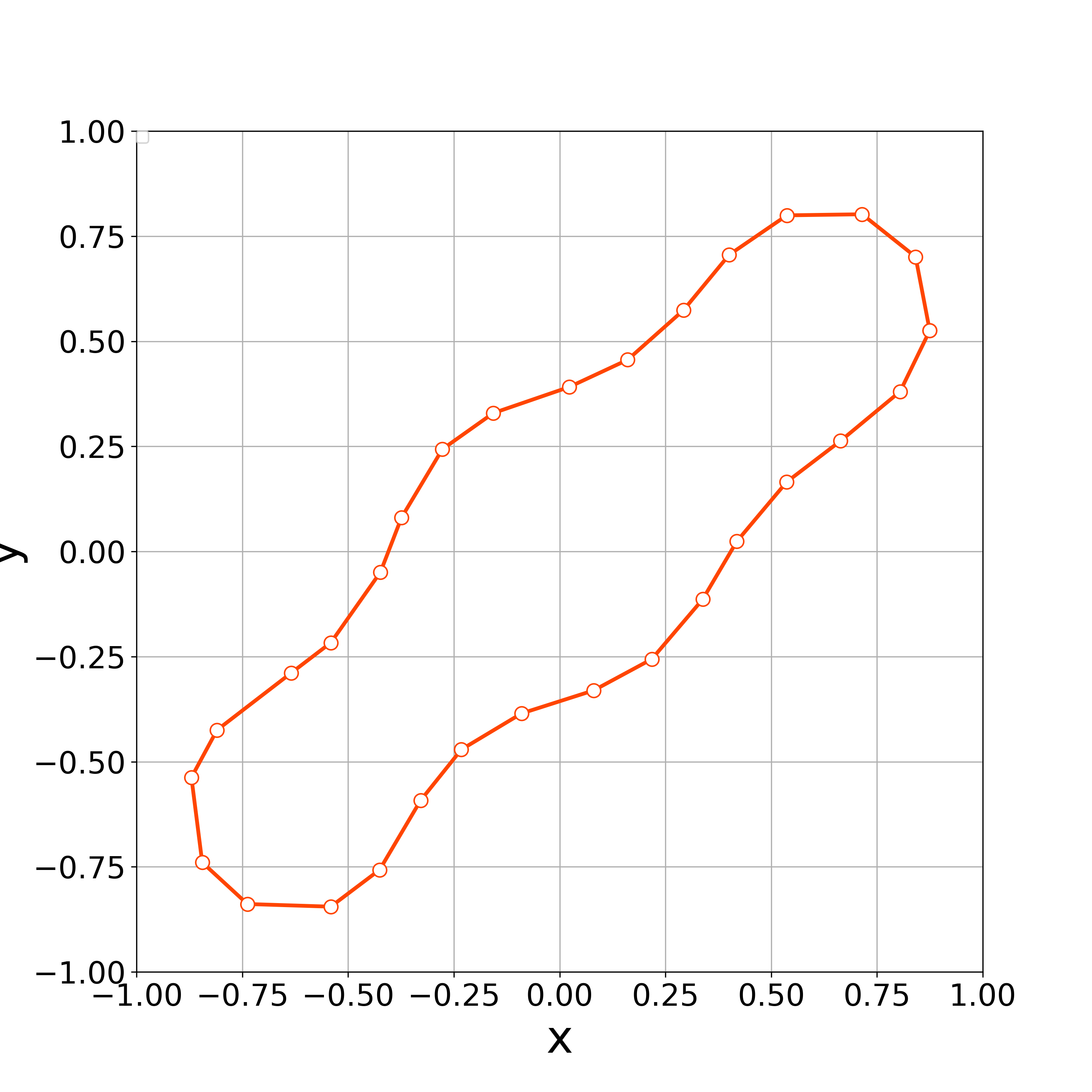}
   \end{center}
   \subcaption{$n=5.$}

  \end{minipage} 

  \begin{minipage}{0.5\hsize}
   \begin{center}
    \includegraphics[width=55mm]{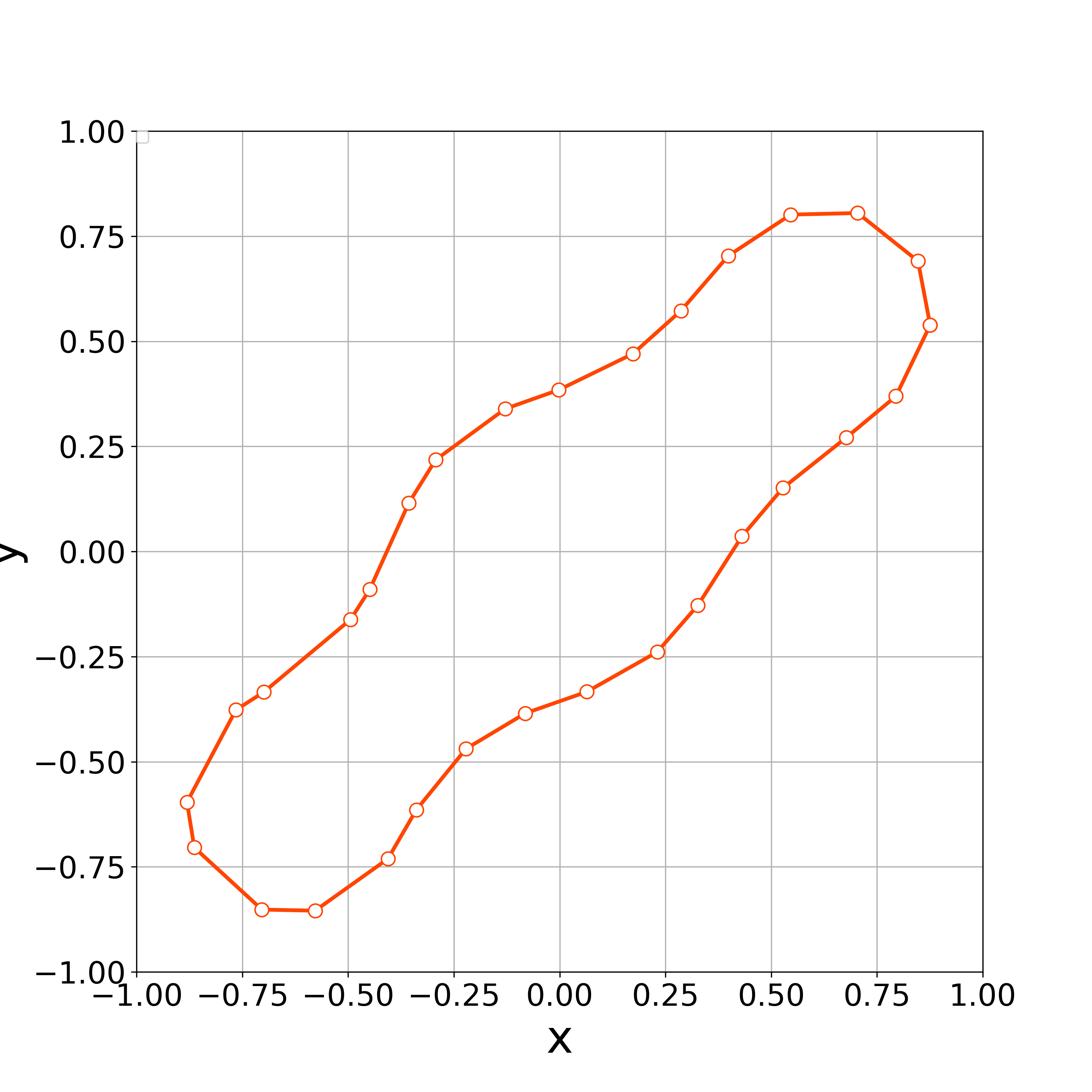}
   \end{center}
   \subcaption{$n=6.$}

  \end{minipage}
  \begin{minipage}{0.5\hsize}
   \begin{center}
    \includegraphics[width=55mm]{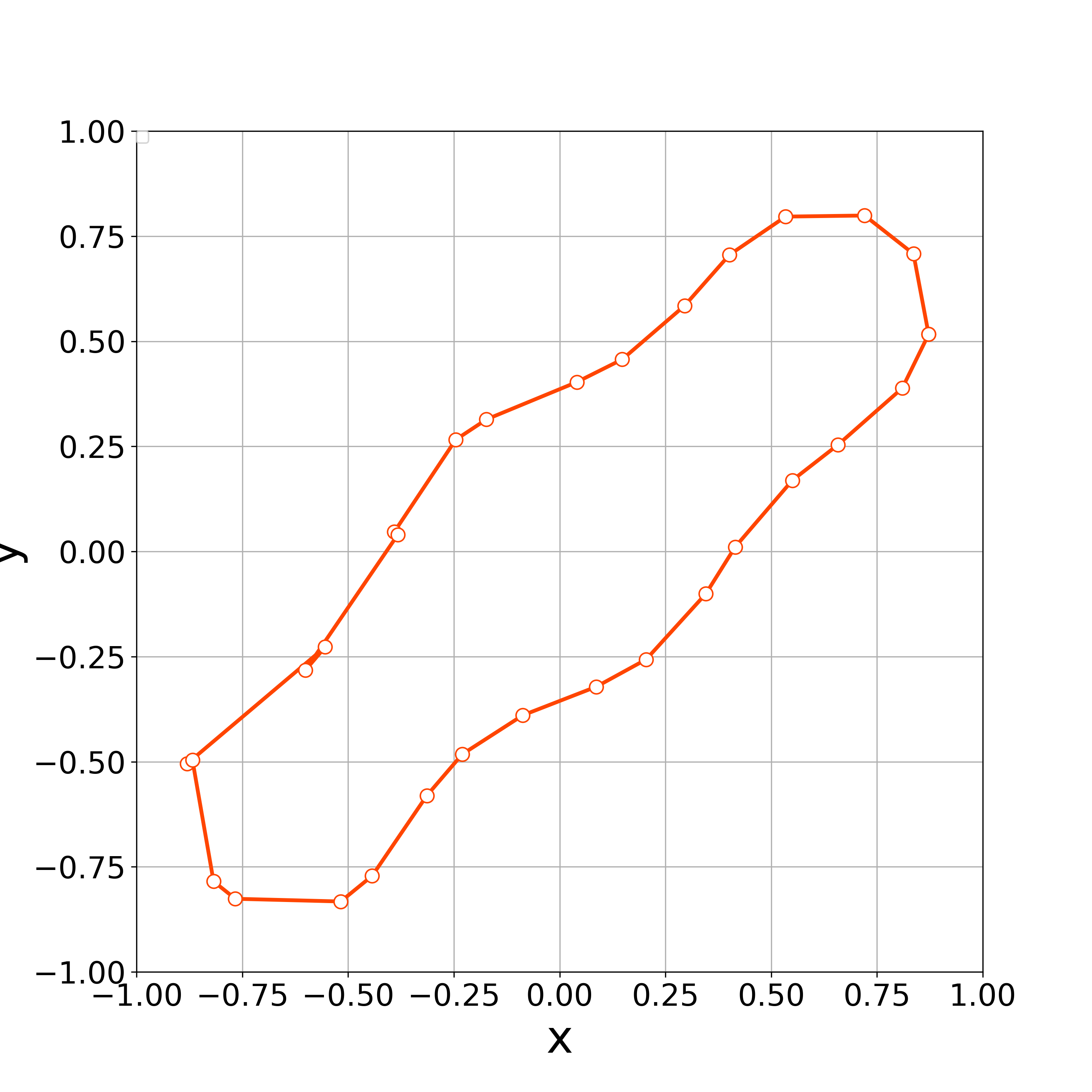}
   \end{center}
   \subcaption{$n=7.$}
  \end{minipage}
  \caption{The $n$th step of the polygonal curve with $\alpha=200$.}
  \label{fig:10}
 \end{figure}

 When $\alpha=200$, the vertices started to oscillate more tangentially from the fifth step, and the nonlinear solver stopped converging in the seventh step (Fig. \ref{fig:10}).

 The above results suggest that the value of $\alpha$ is an important parameter in numerical calculations. If the value of $\alpha$ is too small, the concentration of vertices may occur.
 On the other hand, if the value of $\alpha$ is too big, vertices may oscillate, which causes the concentration of vertices again. Thus, we must choose the value of $\alpha$ appropriately.


\subsubsection{Examination. 2}
\begin{figure}[H]
  \begin{center}
   \includegraphics[width=55mm]{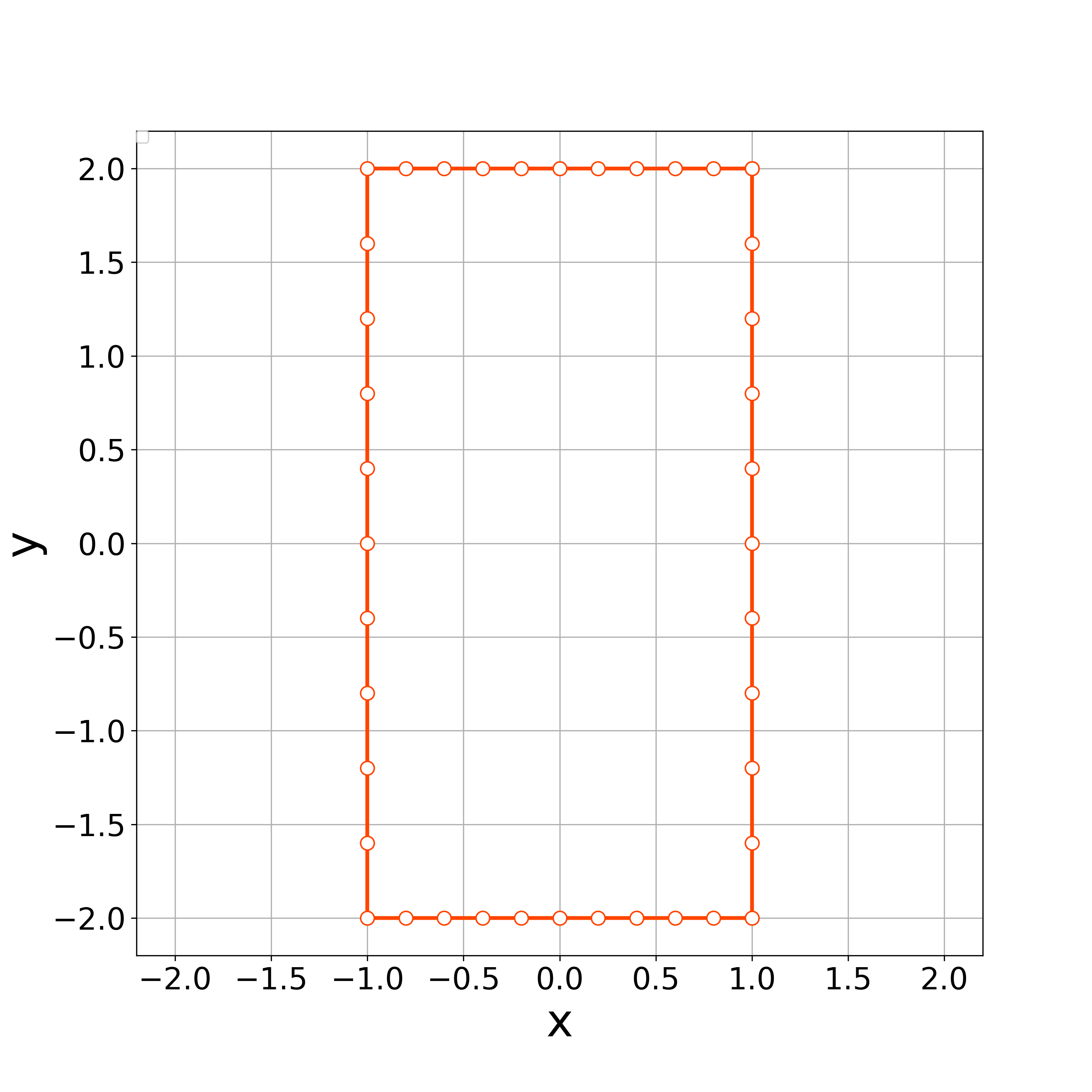}
  \end{center}
  \caption{Rectangular initial arrangement.}
  \label{fig:26}
\end{figure}

\begin{figure}[H]
  \begin{minipage}[b]{0.48\columnwidth}
    \centering
    \includegraphics[width=\columnwidth]{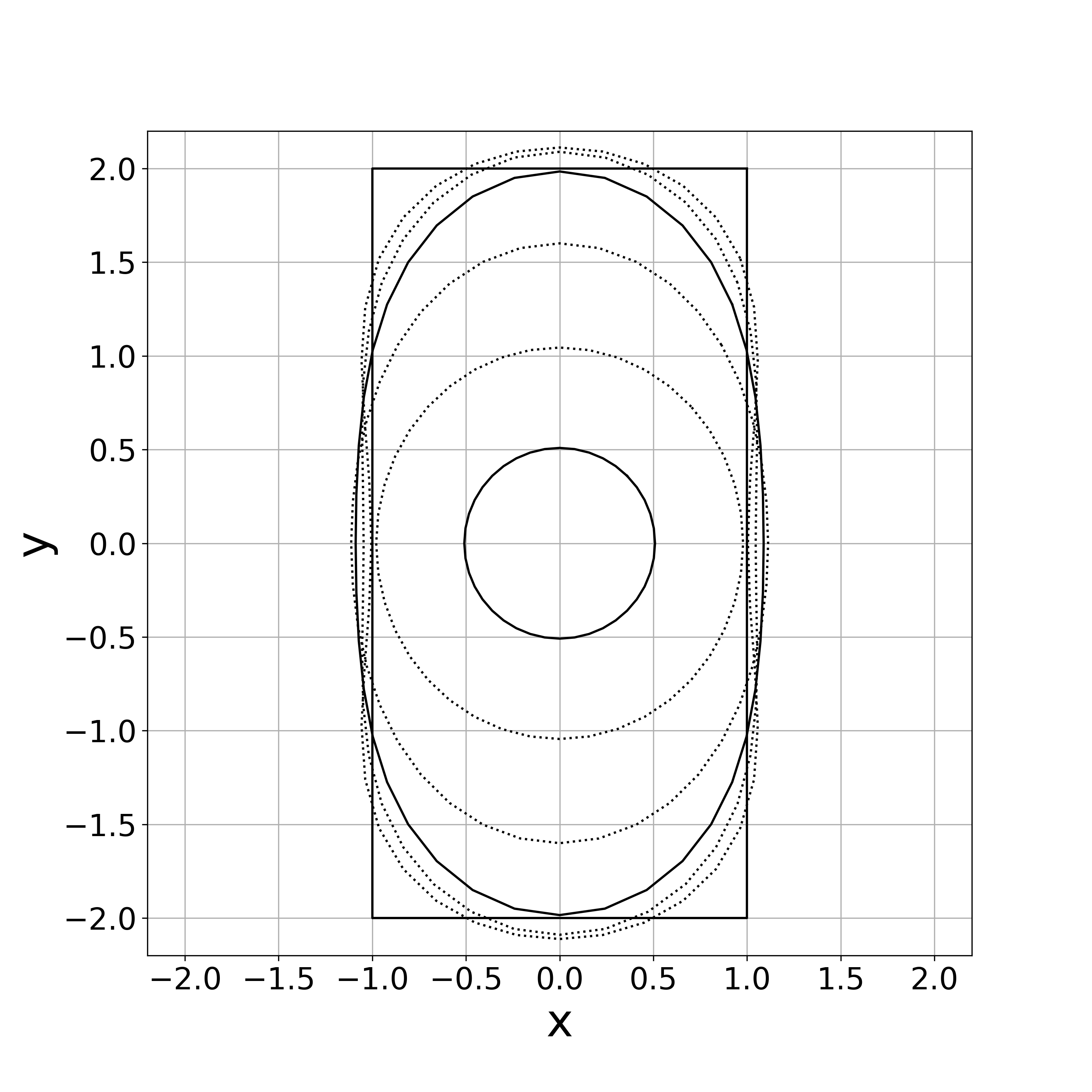}
    \caption{Numerical solutions with $\alpha=50$ (initial rectangular curve and the numerical solutions at $n=50$, 500, 1000, 2500, 5000, and 10000th steps).}
    \label{fig:33}
  \end{minipage}
  \hspace{0.04\columnwidth} 
  \begin{minipage}[b]{0.48\columnwidth}
    \centering
    \mbox{\raisebox{14mm}{\includegraphics[width=\columnwidth]{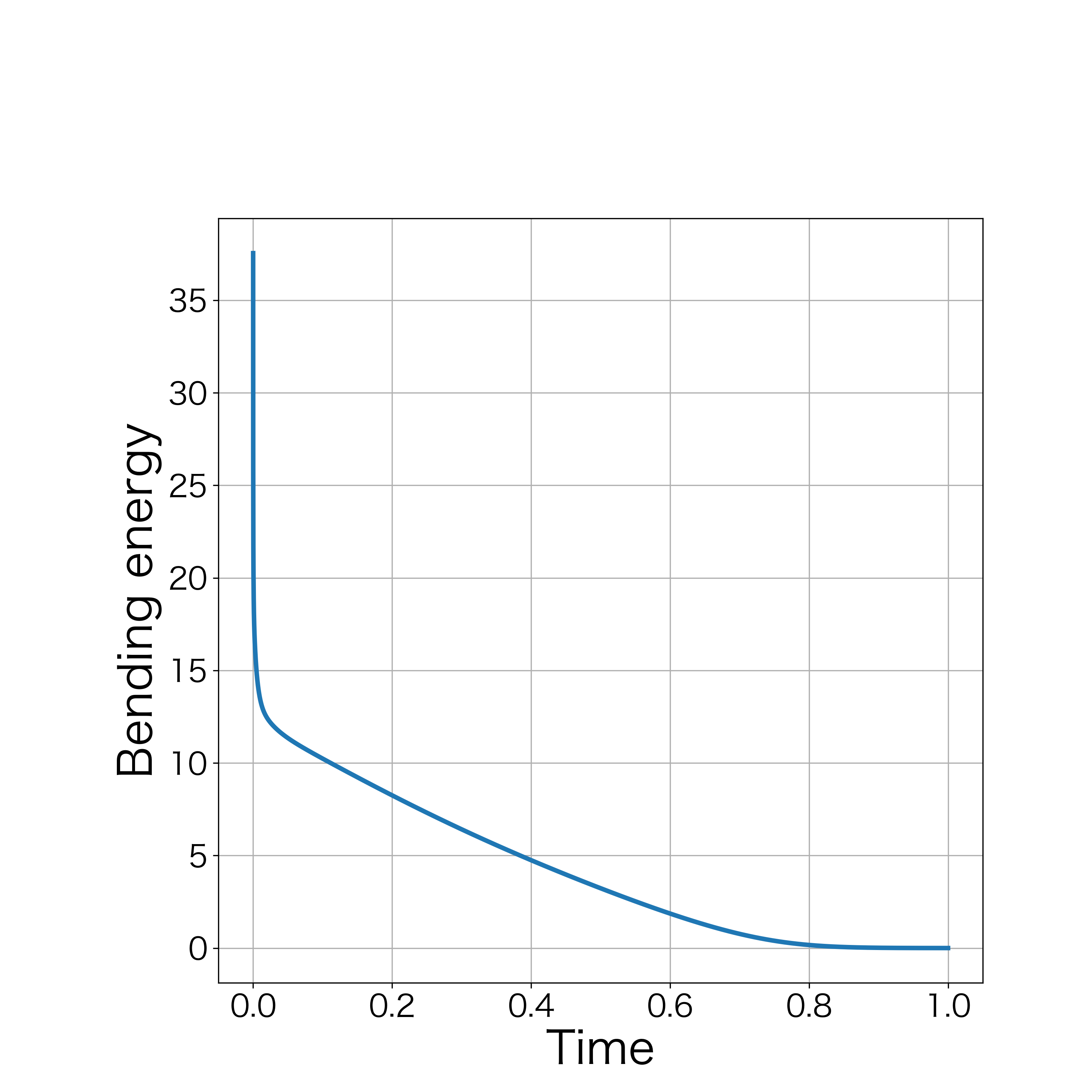}}}
     \vspace*{-1.0cm}
    \caption{The time evolution of the bending energy (rectangular curve).}
    \label{fig:34}
  \end{minipage}
\end{figure}

Numerical experiments are conducted with $N=40$ vertices and the rectangular initial curve shown in Fig.\ \ref{fig:26} . The time step size $\Delta t=10^{-4}$, the constant $c_0=2$, and the tangential velocity parameter $\alpha=30$ are used.
Fig.\ \ref{fig:33} demonstrates the initial polygonal curve and the numerical solutions at $n=50$, 500, 1000, 2500, 5000, and 10000 steps, and we can see that the solution is approaching a regular circle with radius $1/c_0=0.5$. Fig.\ \ref{fig:34} demonstrates the time evolution of the bending energy, and we can observe the discrete version of the energy dissipative nature. 
 


\subsection{Helfrich flow}
In this subsection, the findings of numerical computations using the scheme \eqref{hell_ski} are presented.
\subsubsection{Examination. 3}
Numerical experiments are conducted for various values of the parameter $\alpha$ on tangential velocity to comfirm the impact of $\alpha$ with the initial curve illustrated in Fig.\ \ref{fig:two}, time step size  $\Delta t=10^{-4}$ and the constant $c_0=2$. We set $\alpha=$0,10,100, and 200.

\begin{figure}[H]
  \begin{minipage}{0.5\hsize}
   \begin{center}
    \includegraphics[width=55mm]{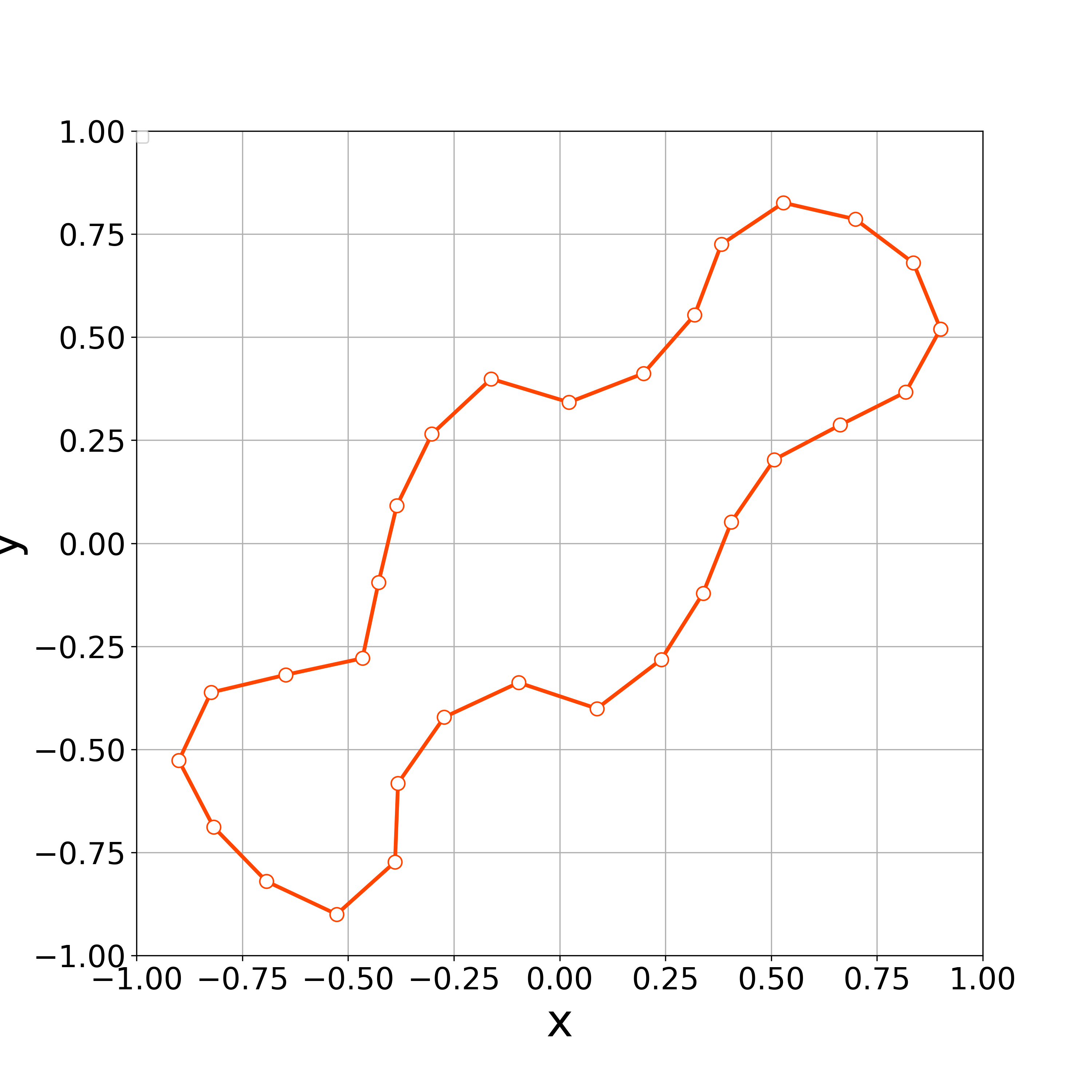}
   \end{center}
   \subcaption{$n=0.$}

  \end{minipage} 
  \begin{minipage}{0.5\hsize}
   \begin{center}
    \includegraphics[width=55mm]{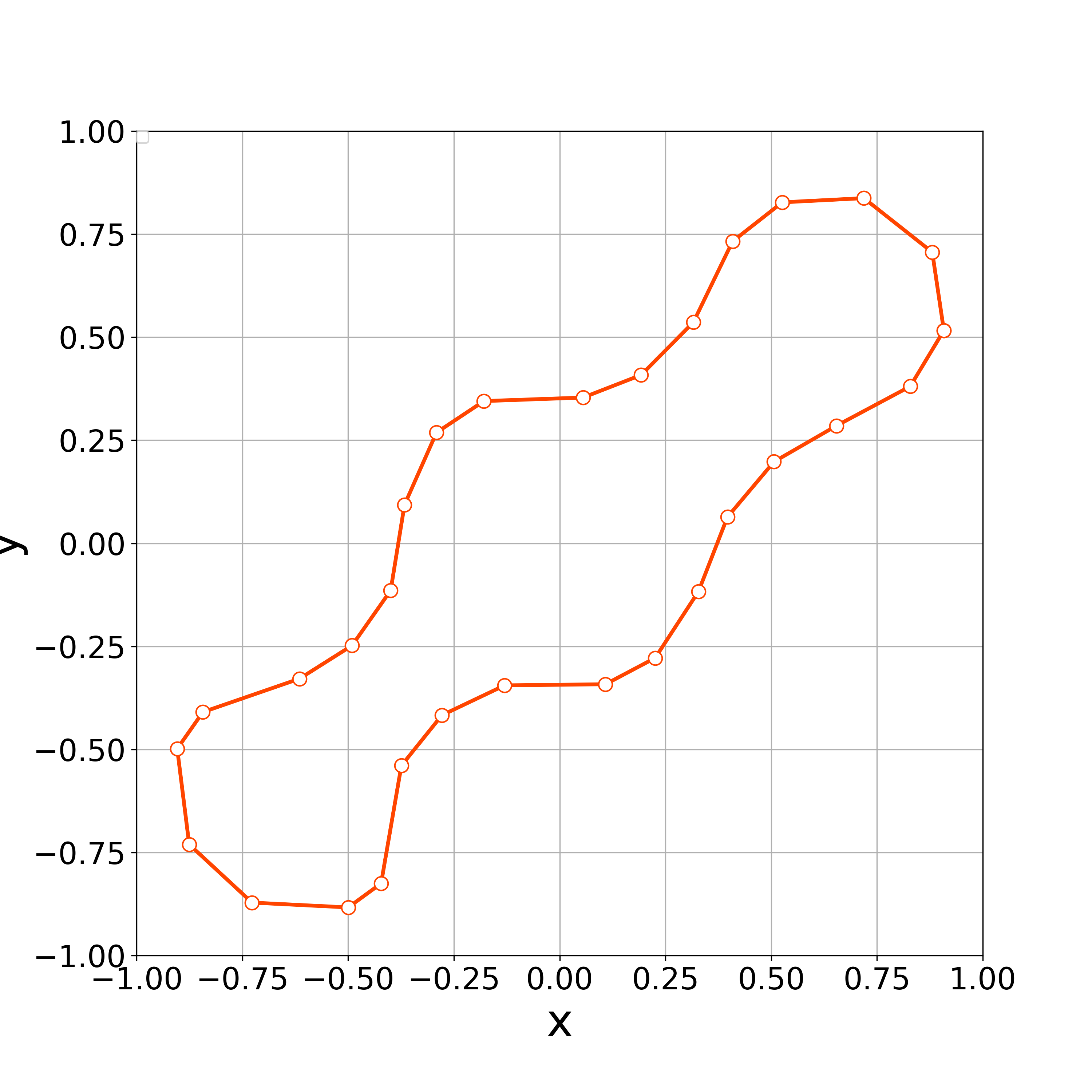}
   \end{center}
   \subcaption{$n=5.$}

  \end{minipage} 

  \begin{minipage}{0.5\hsize}
   \begin{center}
    \includegraphics[width=55mm]{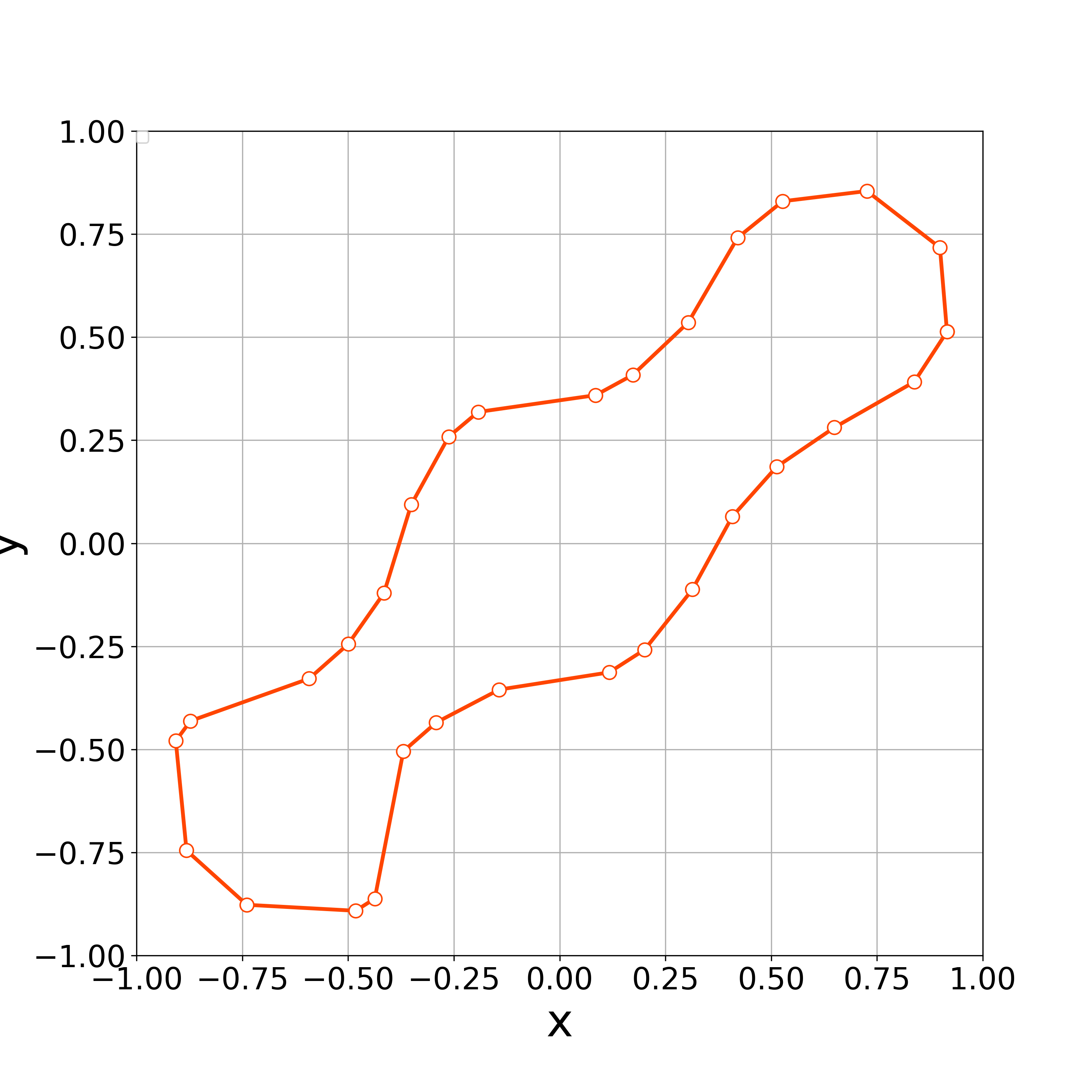}
   \end{center}
   \subcaption{$n=10.$}

  \end{minipage}
  \begin{minipage}{0.5\hsize}
   \begin{center}
    \includegraphics[width=55mm]{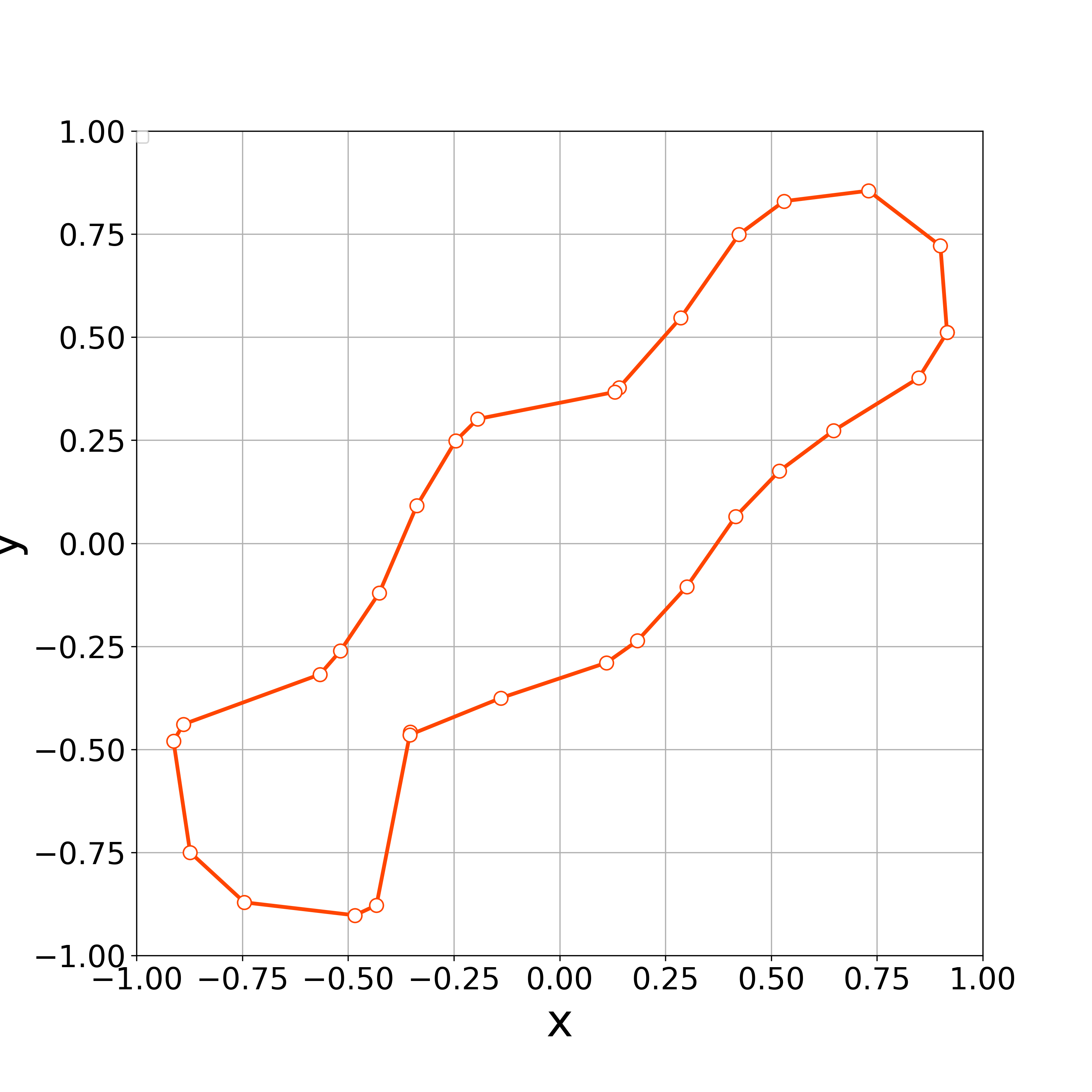}
   \end{center}
   \subcaption{$n=15.$}
  \end{minipage}
  \caption{The $n$th step of the polygonal curve without tangential velocity.}
  \label{fig:42}
 \end{figure}

When tangential velocity is absent, i.e., $\alpha=0$, the nonlinear solver stopped converging at the step 15 because of the concentration of the vertices. Fig.\ \ref{fig:42} demonstrates that the vertices are concentrated at $(x,y)=(-0.4,-0.45)$.

\begin{figure}[H]
  \begin{minipage}{0.5\hsize}
   \begin{center}
    \includegraphics[width=55mm]{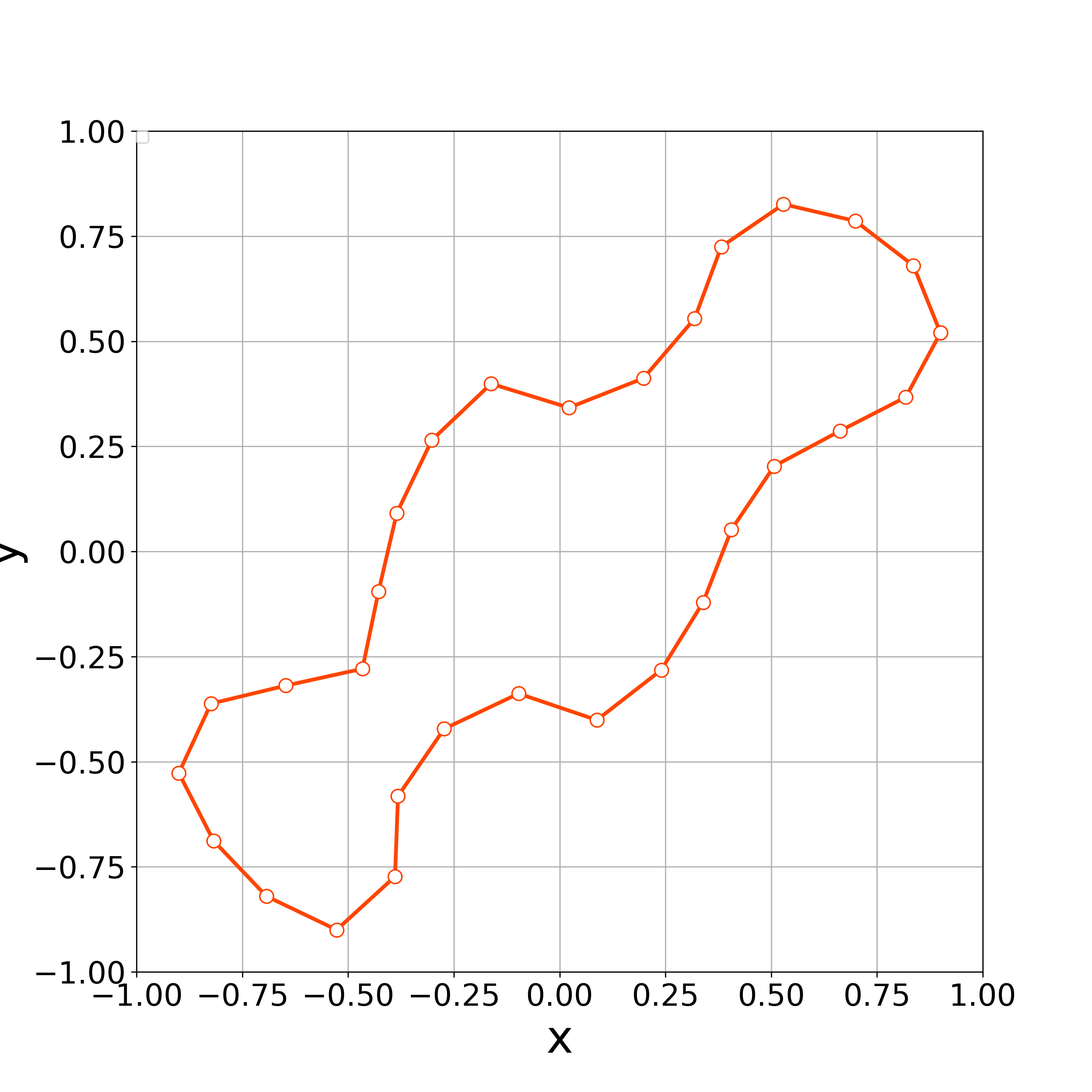}
   \end{center}
   \subcaption{$n=0.$}

  \end{minipage} 
  \begin{minipage}{0.5\hsize}
   \begin{center}
    \includegraphics[width=55mm]{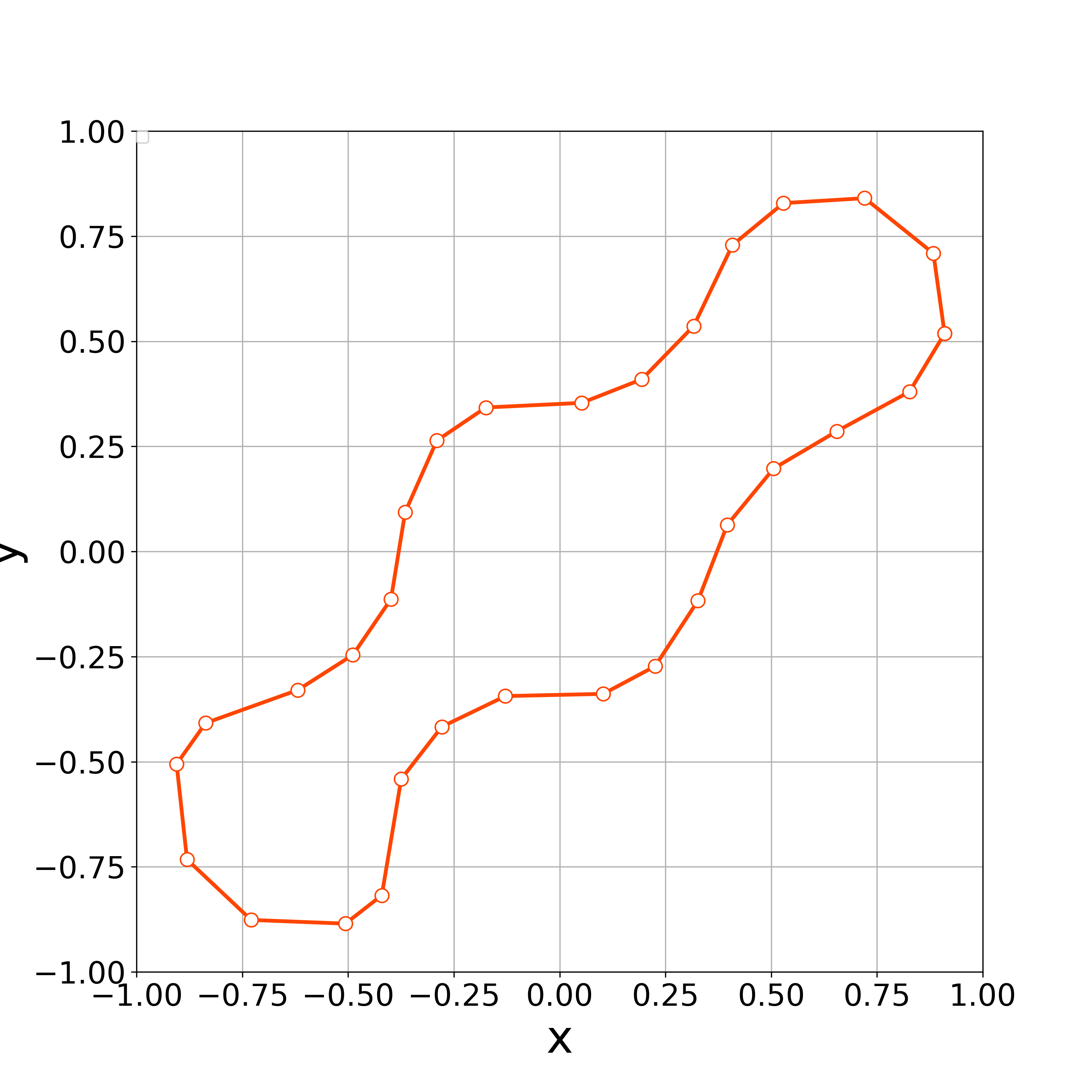}
   \end{center}
   \subcaption{$n=5.$}

  \end{minipage} 

  \begin{minipage}{0.5\hsize}
   \begin{center}
    \includegraphics[width=55mm]{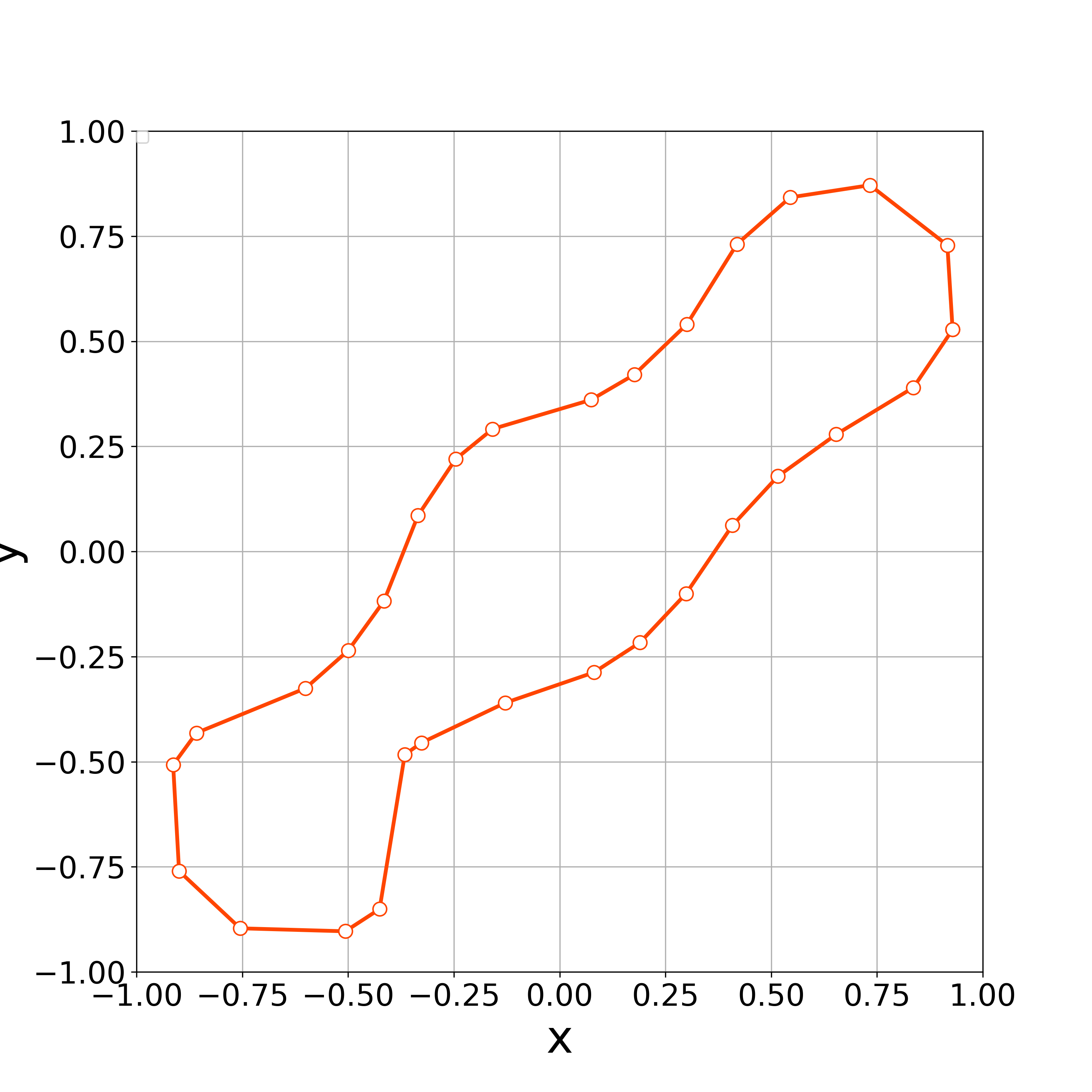}
   \end{center}
   \subcaption{$n=12.$}

  \end{minipage}
  \begin{minipage}{0.5\hsize}
   \begin{center}
    \includegraphics[width=55mm]{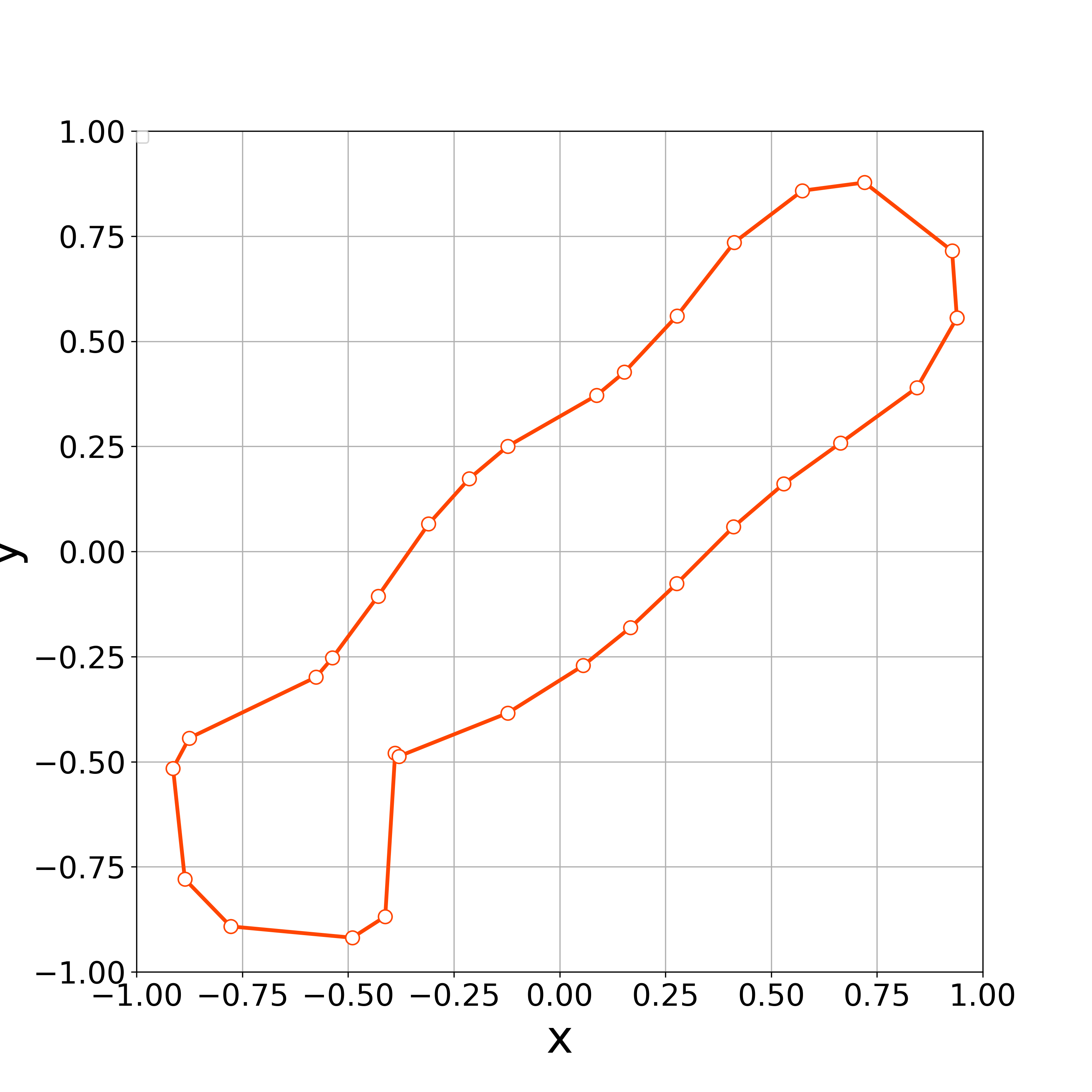}
   \end{center}
   \subcaption{$n=16.$}
  \end{minipage}
  \caption{The $n$th step of the polygonal curve with $\alpha=10$.}
  \label{fig:14}
 \end{figure}

When $\alpha=10$, the nonlinear solver stopped converging at the step 16 because of the concentration of vertices. Fig.\ \ref{fig:14} demonstrates that the vertices are concentrated at $(x,y)=(-0.4,-0.5)$.

\begin{figure}[H]
  \begin{minipage}{0.5\hsize}
   \begin{center}
    \includegraphics[width=55mm]{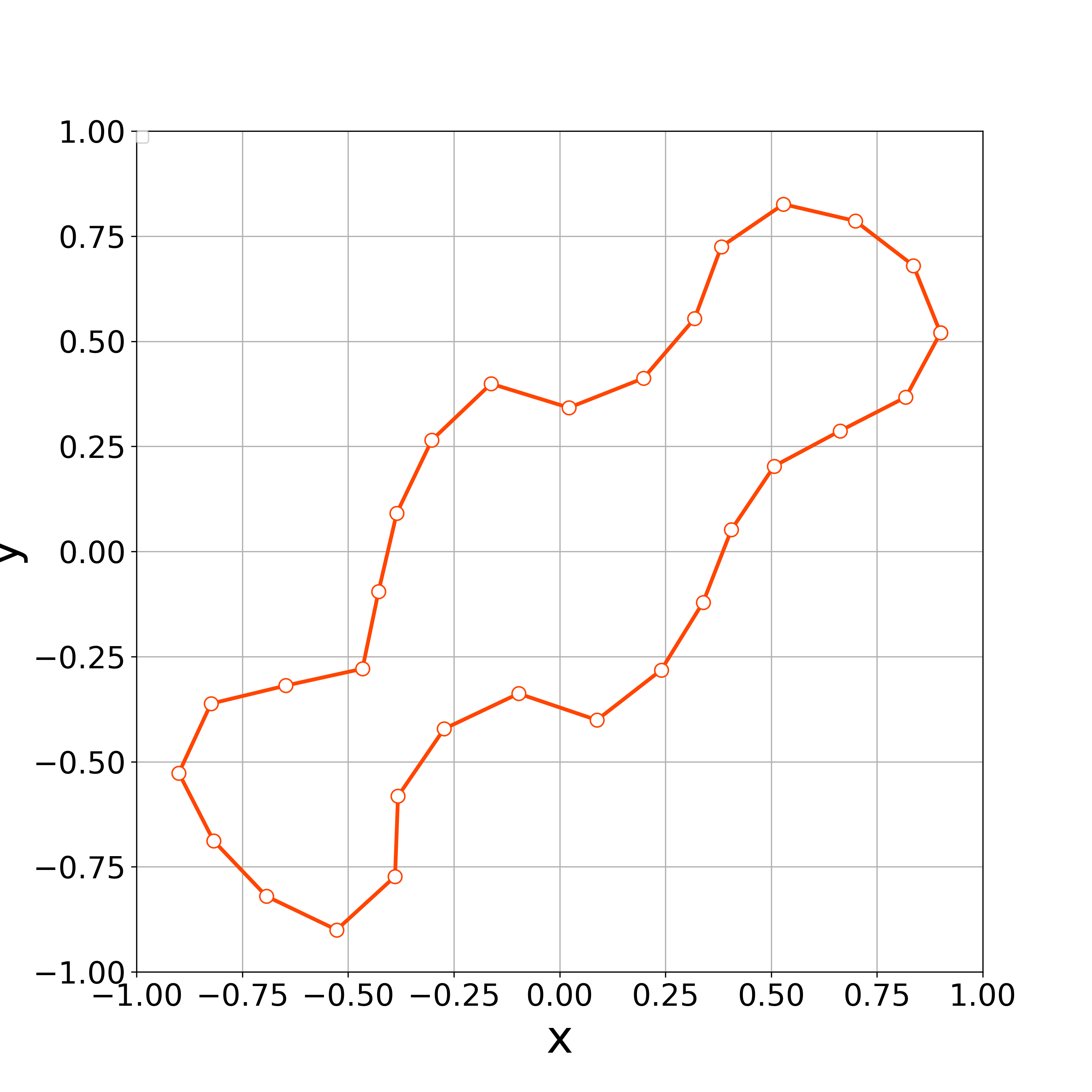}
   \end{center}
   \subcaption{$n=0.$}

  \end{minipage} 
  \begin{minipage}{0.5\hsize}
   \begin{center}
    \includegraphics[width=55mm]{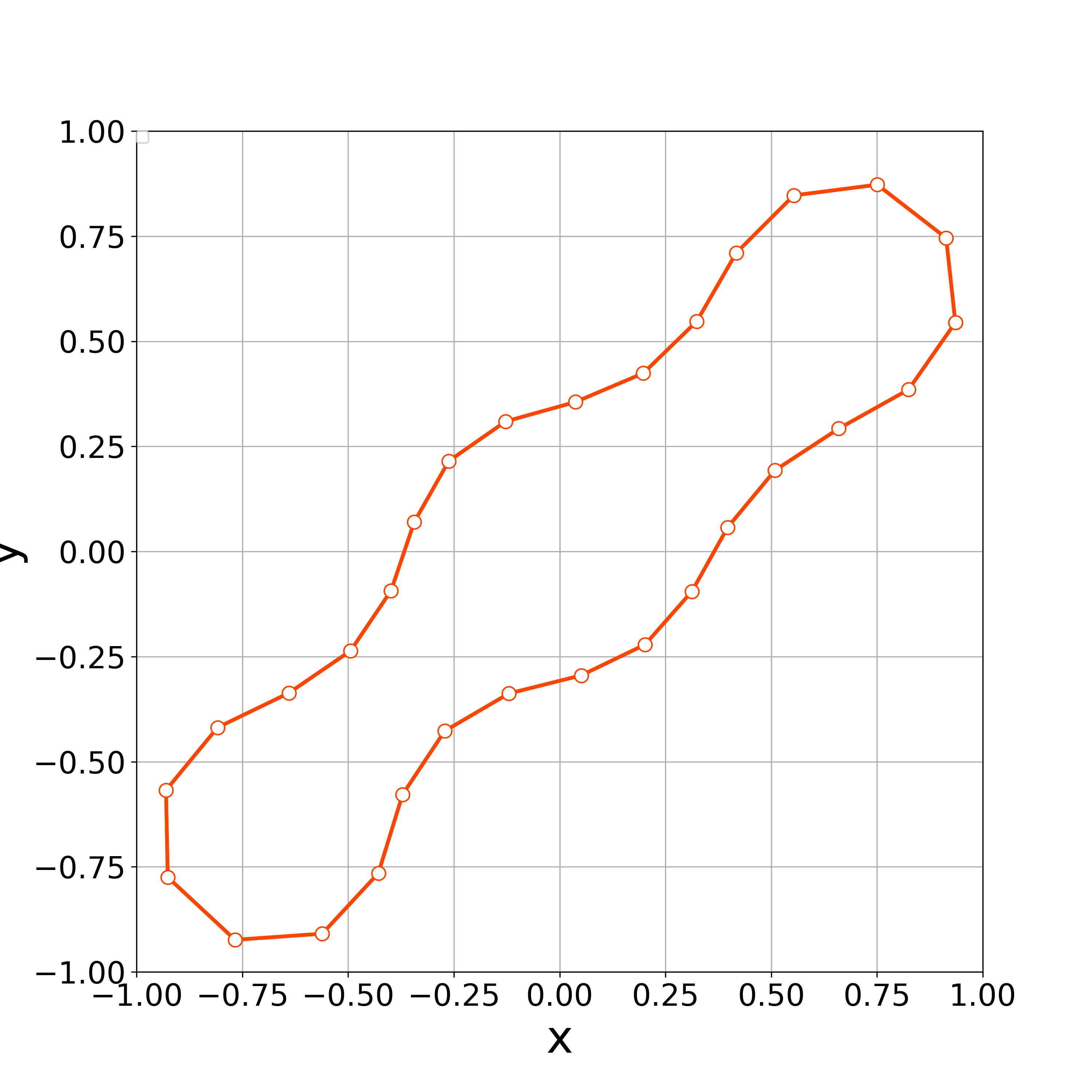}
   \end{center}
   \subcaption{$n=10.$}

  \end{minipage} 

  \begin{minipage}{0.5\hsize}
   \begin{center}
    \includegraphics[width=55mm]{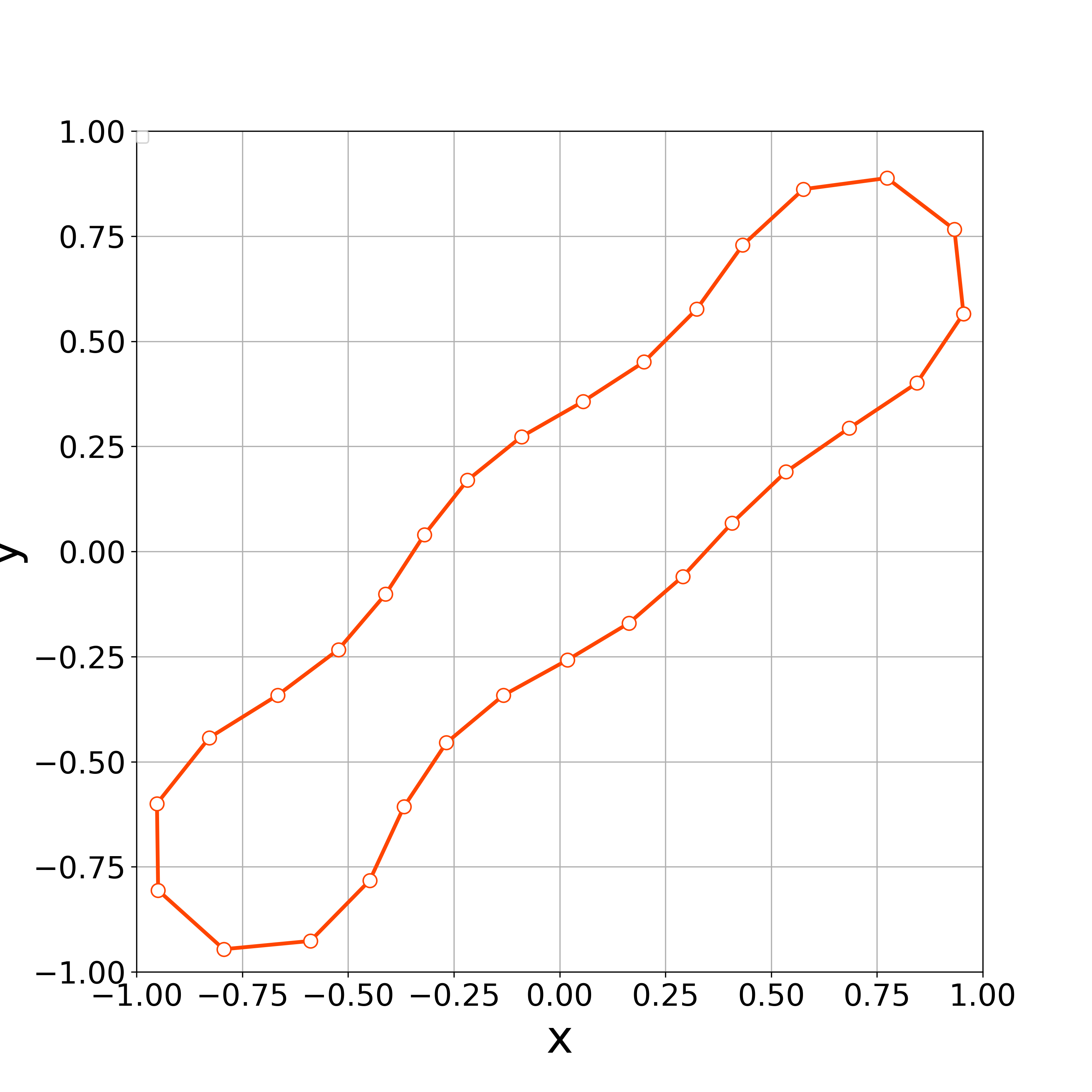}
   \end{center}
   \subcaption{$n=20.$}

  \end{minipage}
  \begin{minipage}{0.5\hsize}
   \begin{center}
    \includegraphics[width=55mm]{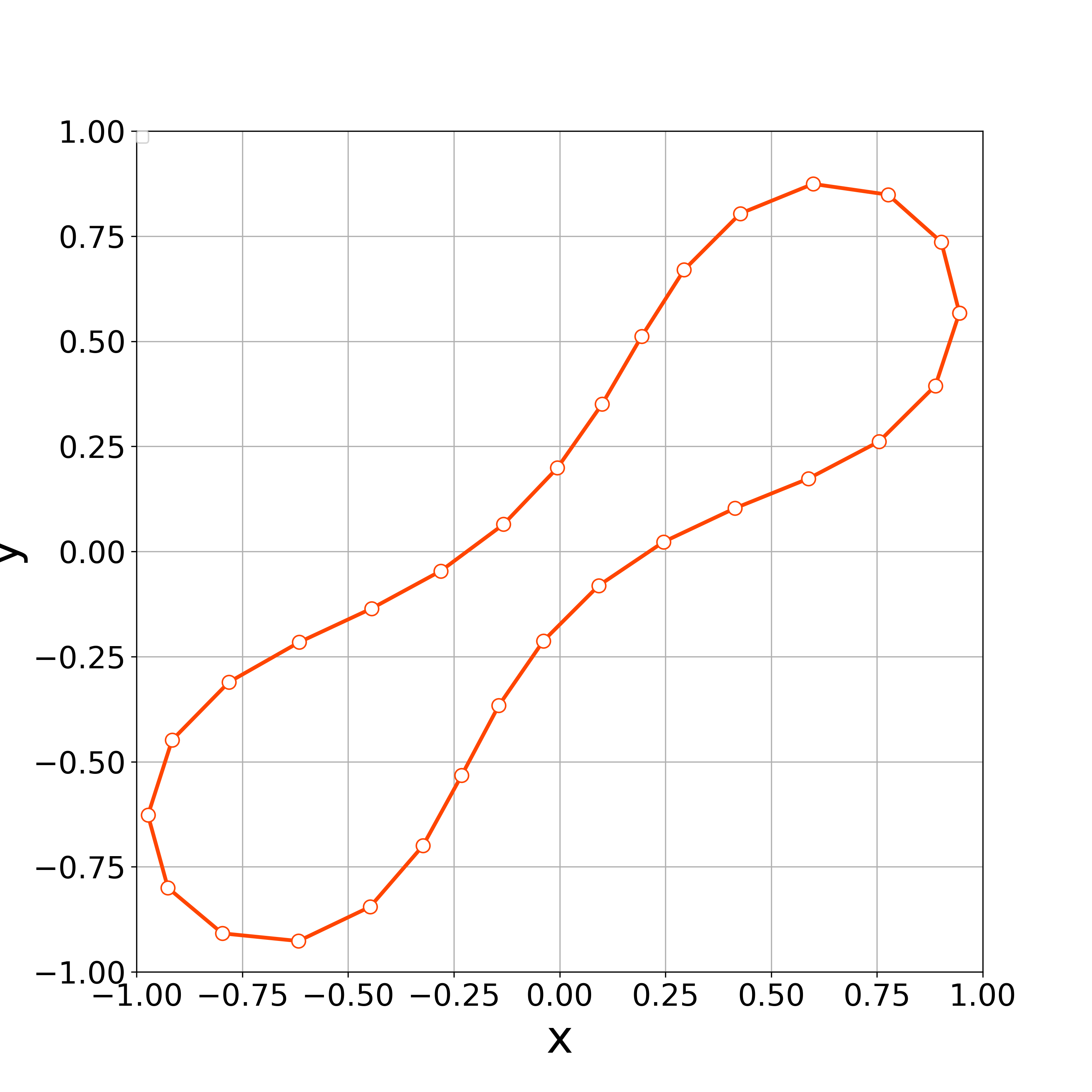}
   \end{center}
   \subcaption{$n=100.$}
  \end{minipage}
  \caption{The $n$th step of the polygonal curve with $\alpha=100$.}
  \label{fig:18}
 \end{figure}

\begin{figure}[H]
  \centering
  \includegraphics[width=55mm]{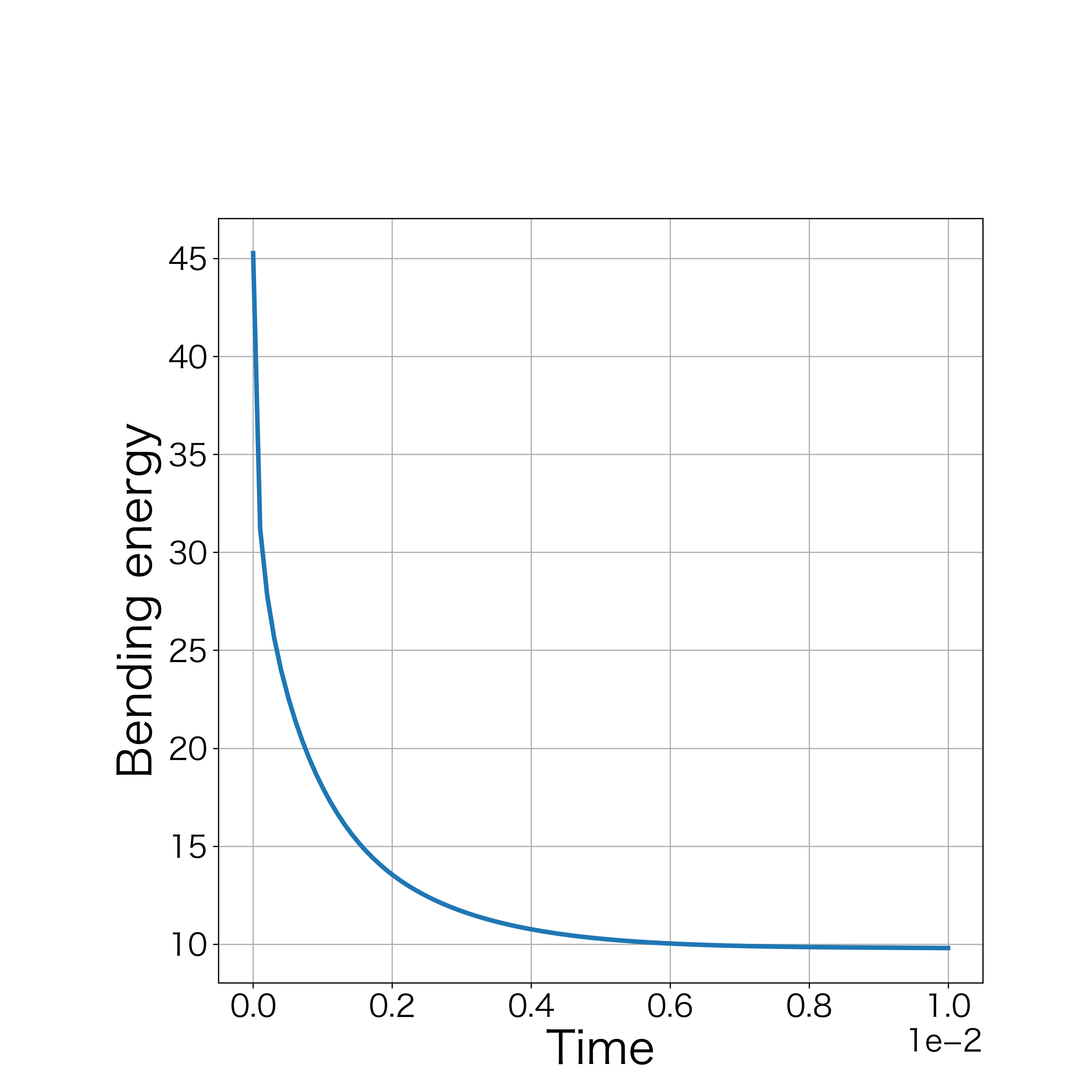}
  \caption{Time evolution of the bending energy ($\alpha=100$).}
    \label{fig:19}
   \end{figure}

\begin{figure}[H]
    \begin{minipage}{0.5\hsize}
     \begin{center}
      \includegraphics[width=55mm]{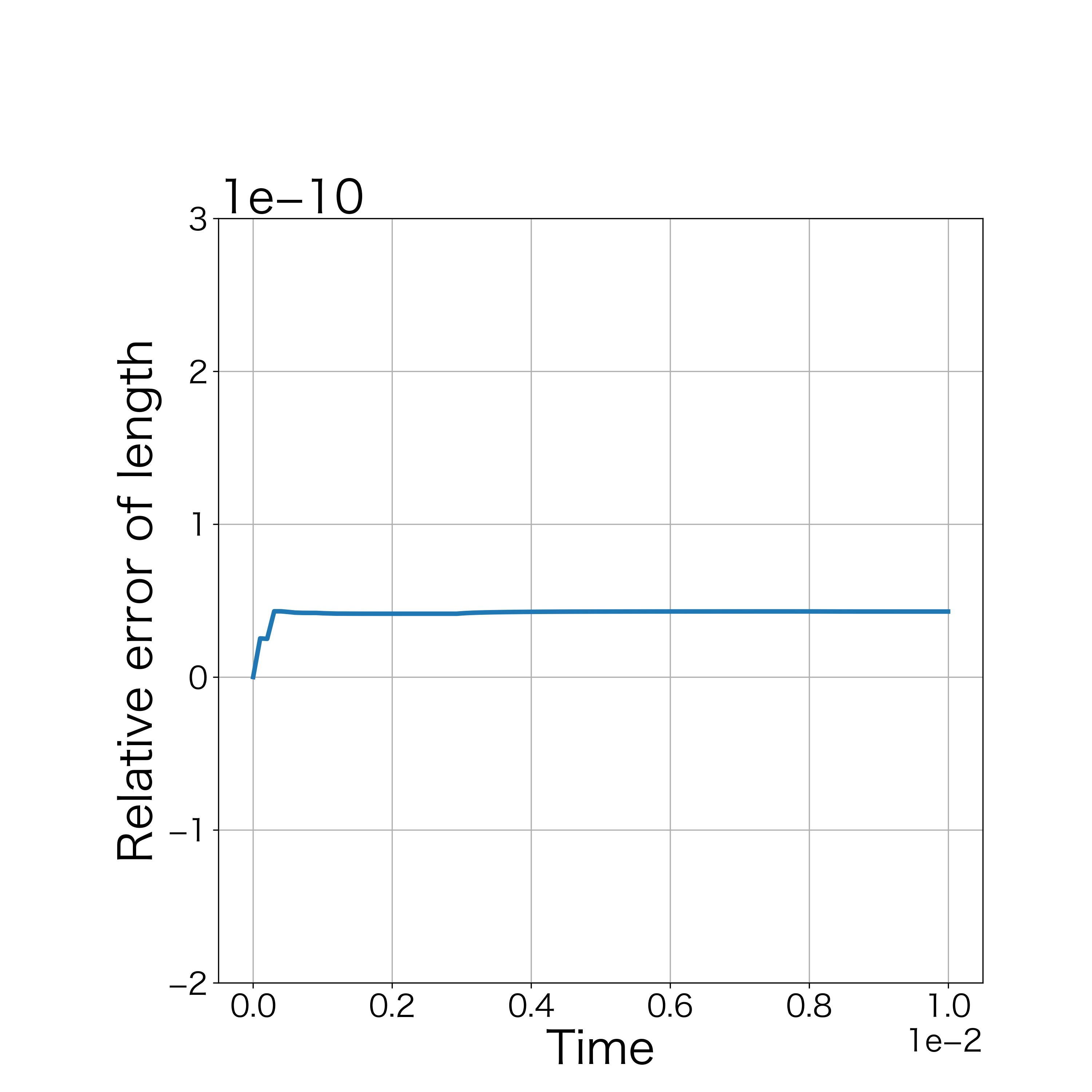}
     \end{center}
     \label{fig:20}
    \end{minipage}
    \begin{minipage}{0.5\hsize}
     \begin{center}
      \includegraphics[width=55mm]{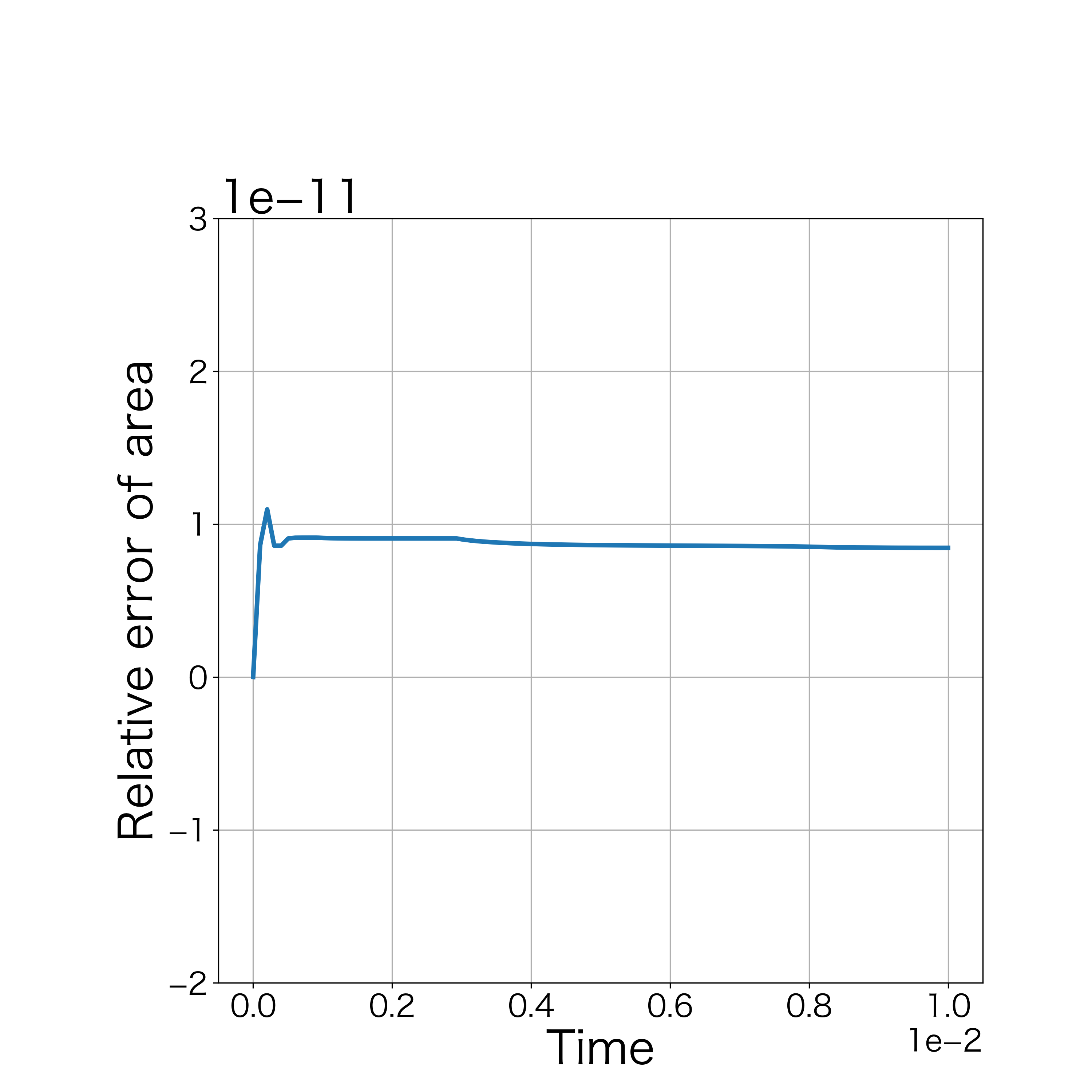}
     \end{center}
     \label{fig:21}
    \end{minipage}
    \caption{Relative error of the length and the enclosed area ($\alpha=100$).}
    \label{fig:20-21}
   \end{figure}

When $\alpha=100$, the numerical computations demonstrates that the nonlinear solver converges up to the 100th step. Fig.\ \ref{fig:18} demonstrates the initial polygonal curve and the numerical solutions at $n=10$, 20, and 100 steps.  Fig.\ \ref{fig:19}--\ref{fig:20-21} respectively demonstrates the time evolution of the bending energy, the relative error of the length and the enclosed area for each step. The discrete bending energy's dissipative nature and the conservation of the length and the enclosed area can be observed. Note that the relative error of the energy $E$ ($E=A_d,L_d$) at the $n$th step is defined by 
\begin{equation}
  \frac{E^{(n)}-E^{(0)}}{E^{(0)}}.
\end{equation}
\clearpage



\begin{figure}[H]
    \begin{minipage}{0.5\hsize}
     \begin{center}
      \includegraphics[width=55mm]{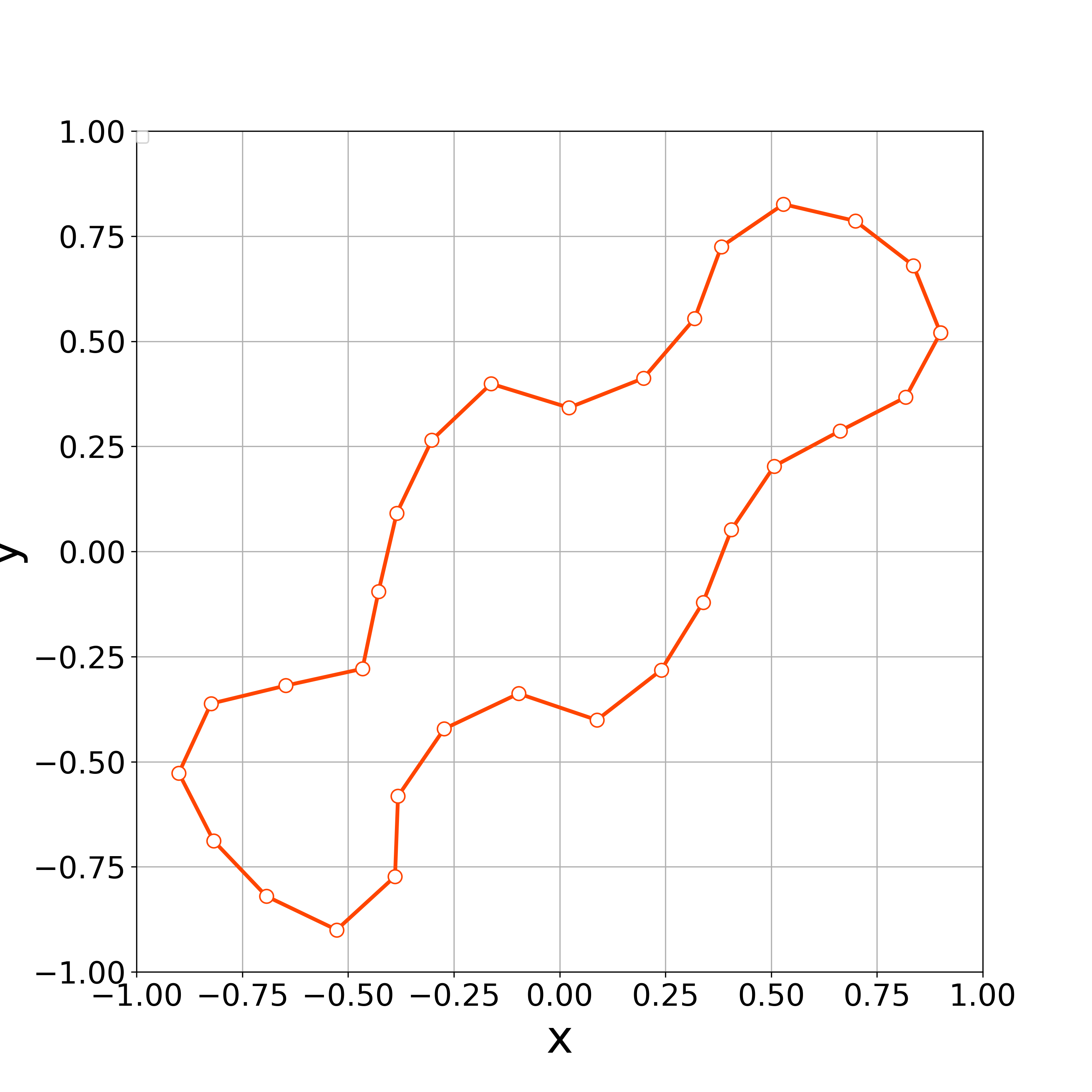}
     \end{center}
     \subcaption{$n=0.$}
  
    \end{minipage} 
    \begin{minipage}{0.5\hsize}
     \begin{center}
      \includegraphics[width=55mm]{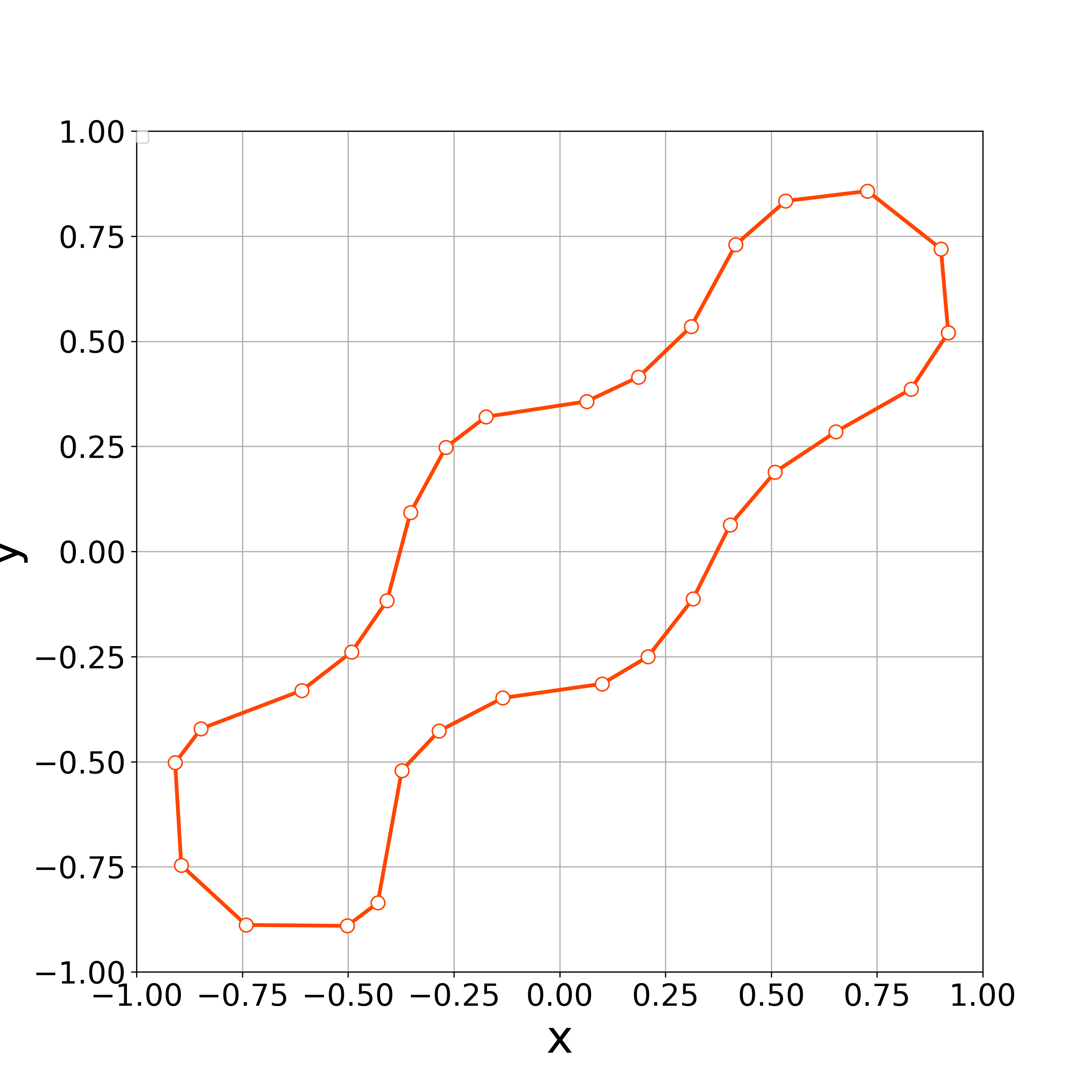}
     \end{center}
     \subcaption{$n=8.$}
  
    \end{minipage} 
  
    \begin{minipage}{0.5\hsize}
     \begin{center}
      \includegraphics[width=55mm]{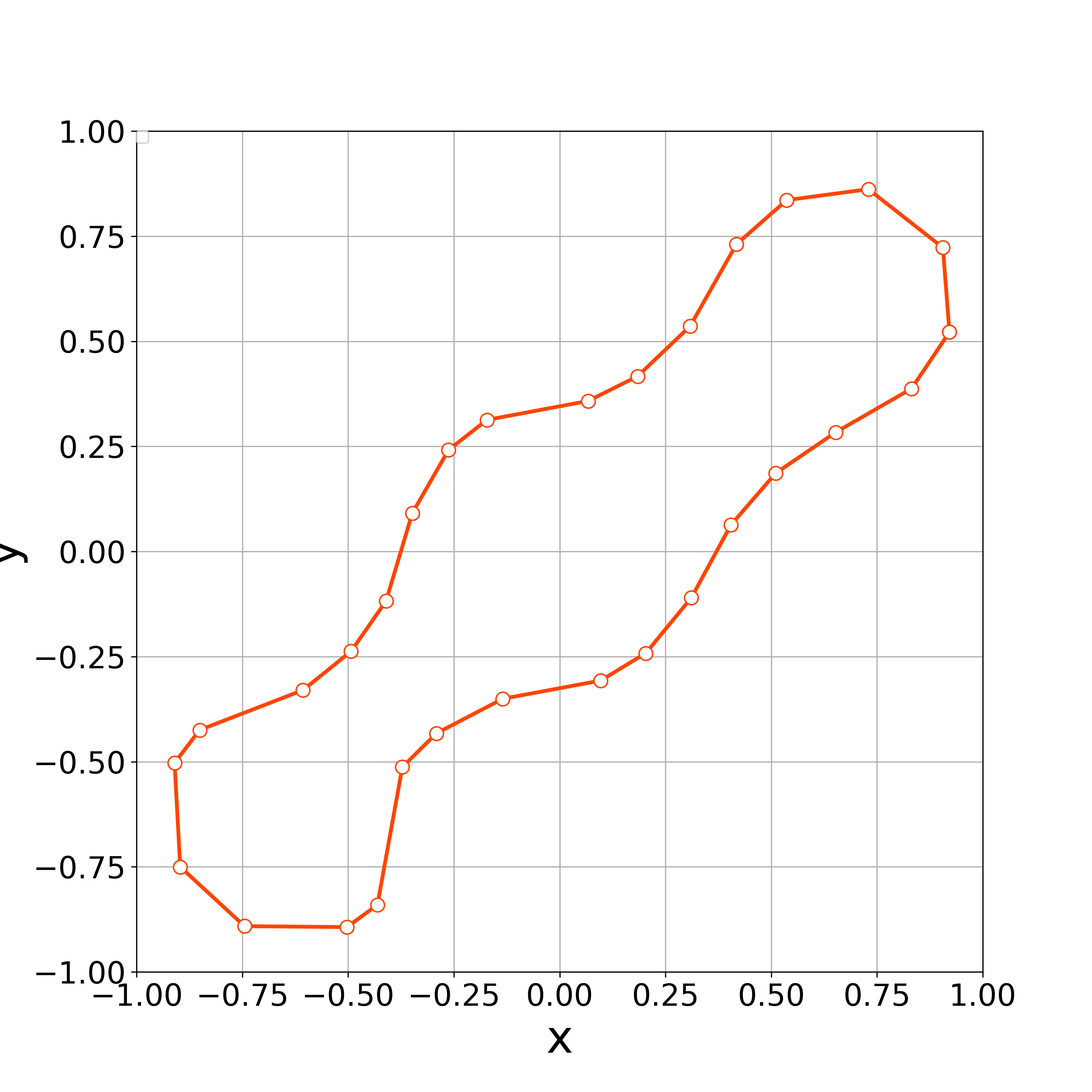}
     \end{center}
     \subcaption{$n=9.$}
  
    \end{minipage}
    \begin{minipage}{0.5\hsize}
     \begin{center}
      \includegraphics[width=55mm]{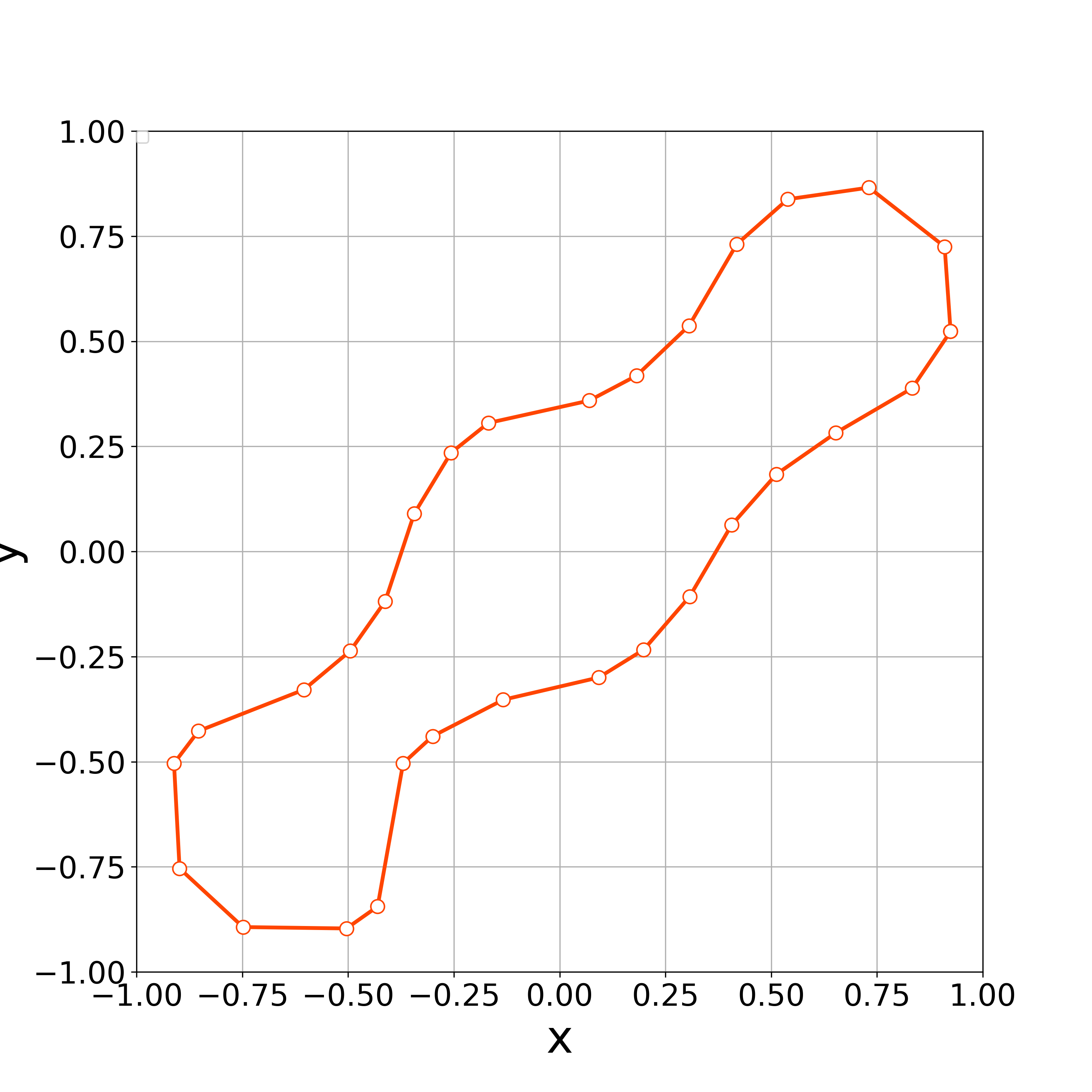}
     \end{center}
     \subcaption{$n=10.$}
    \end{minipage}
    \caption{The $n$th step of the polygonal curve with $\alpha=200$.}
    \label{fig:25}
   \end{figure}

   When $\alpha=200$, the vertices started to oscillate more tangentially from the 8th step, and the nonlinear solver stopped converging at the 10th step (Fig.\ \ref{fig:25}).

  Even in the case of the Helfrich flow, the numerical results suggest that the value of $\alpha$ is a crucial parameter in numerical computations and that it must be selected appropriately.
 
\subsubsection{Examination. 4}
\begin{figure}[H]
  \begin{minipage}{0.5\hsize}
   \begin{center}
    \includegraphics[width=55mm]{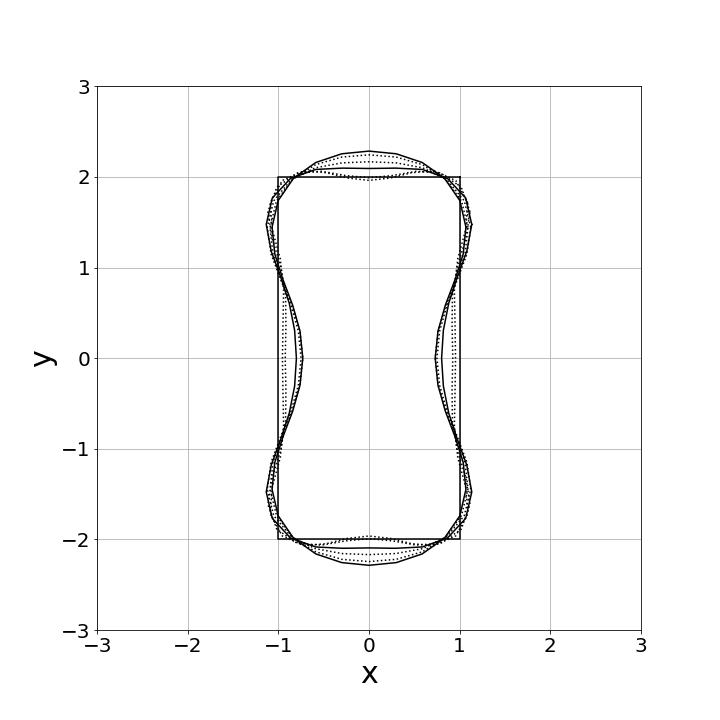}
   \end{center}
   \caption{Numerical solutions with $\alpha=50$ \\(initial rectangular curve and the numerical \\solutions at $n=10$, 50, 100, 300, 500 and\\ 1000th steps).}
   \label{fig:27}
  \end{minipage}
  \begin{minipage}{0.5\hsize}
   \begin{center}
    \includegraphics[width=55mm]{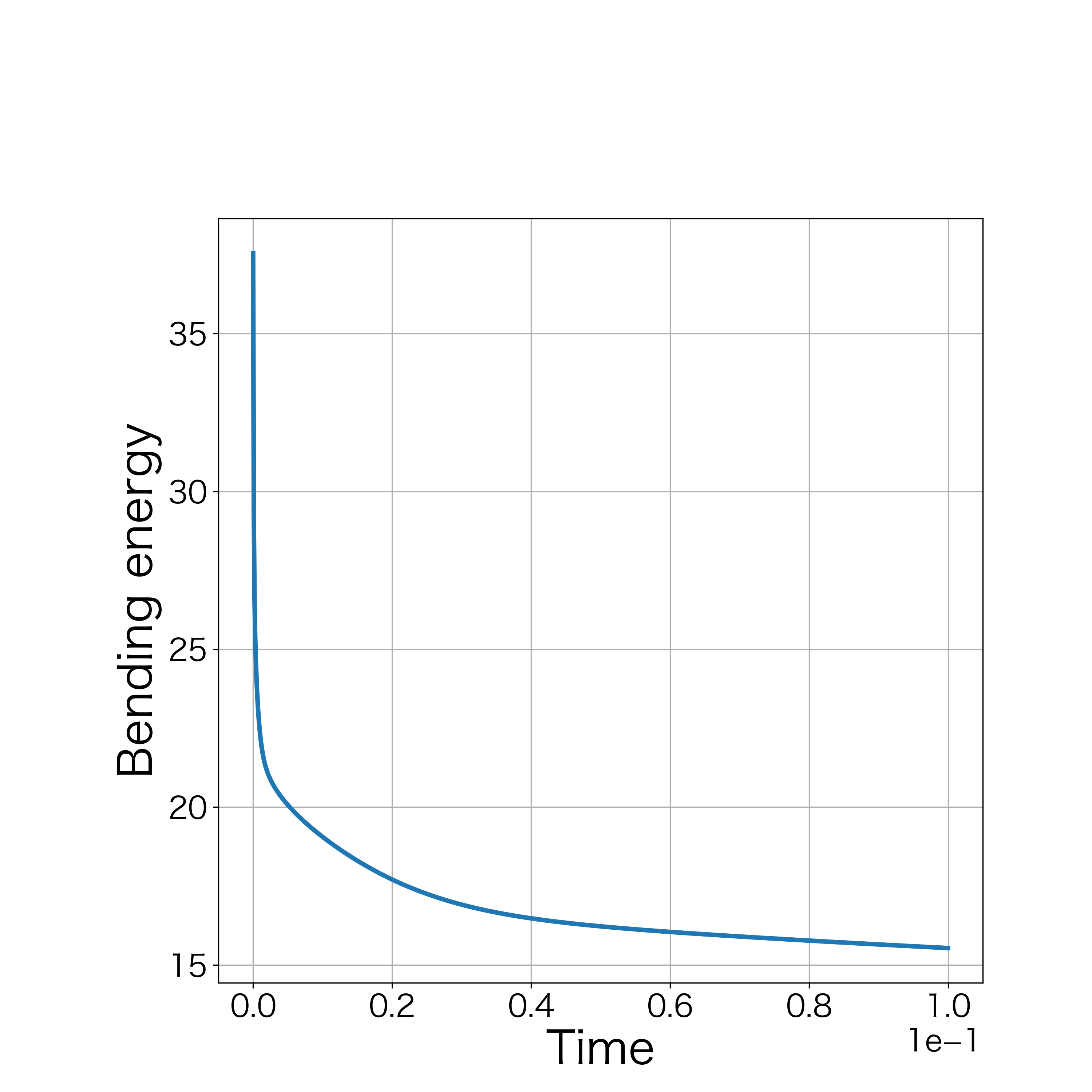}
   \end{center}
   \caption{The time evolution of the bending energy (rectangular curve).}
   
   \label{fig:28}
  \end{minipage}
\end{figure}

\begin{figure}[h]
  \begin{minipage}{0.5\hsize}
   \begin{center}
    \includegraphics[width=55mm]{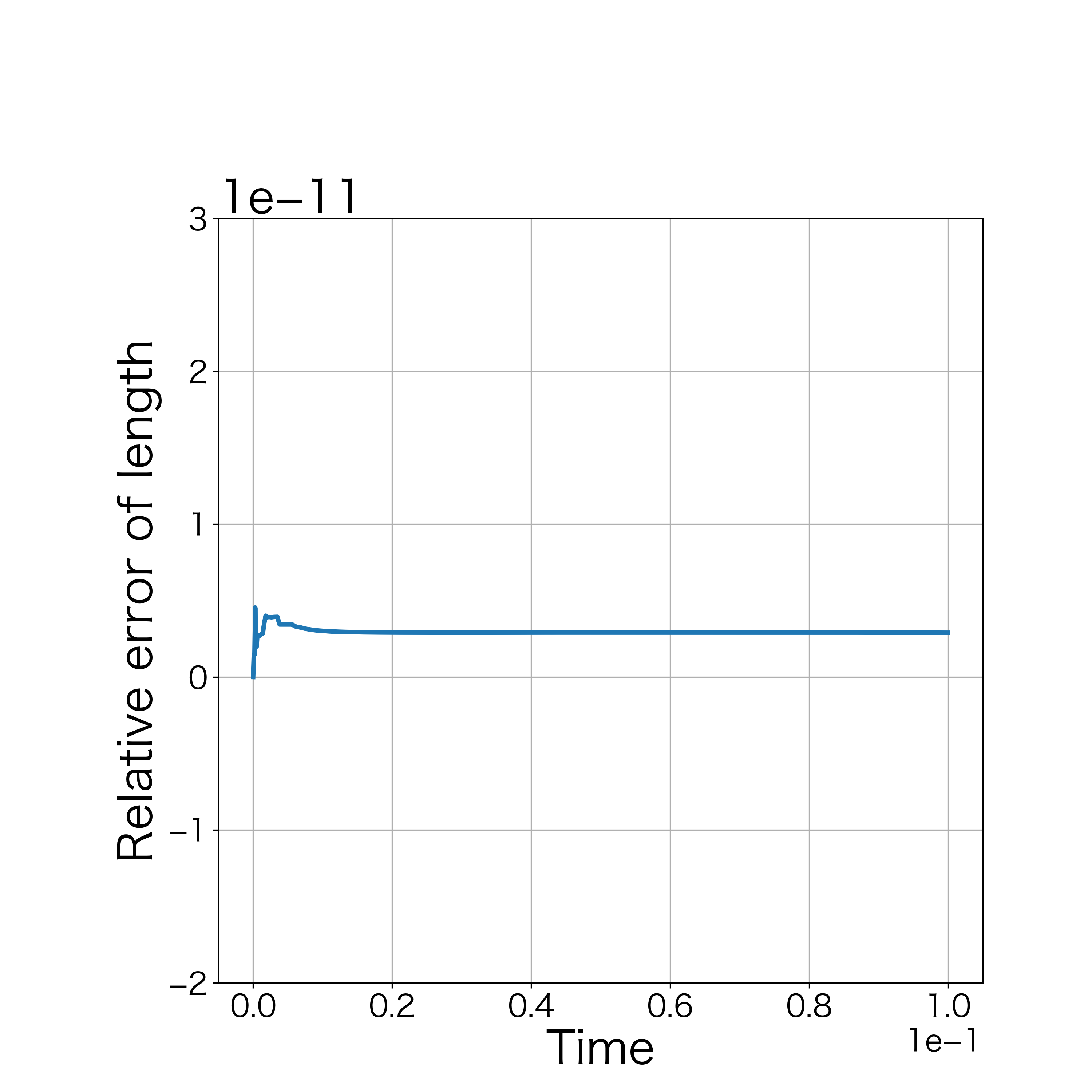}
   \end{center}
  \end{minipage}
  \begin{minipage}{0.5\hsize}
    \begin{center}
     \includegraphics[width=55mm]{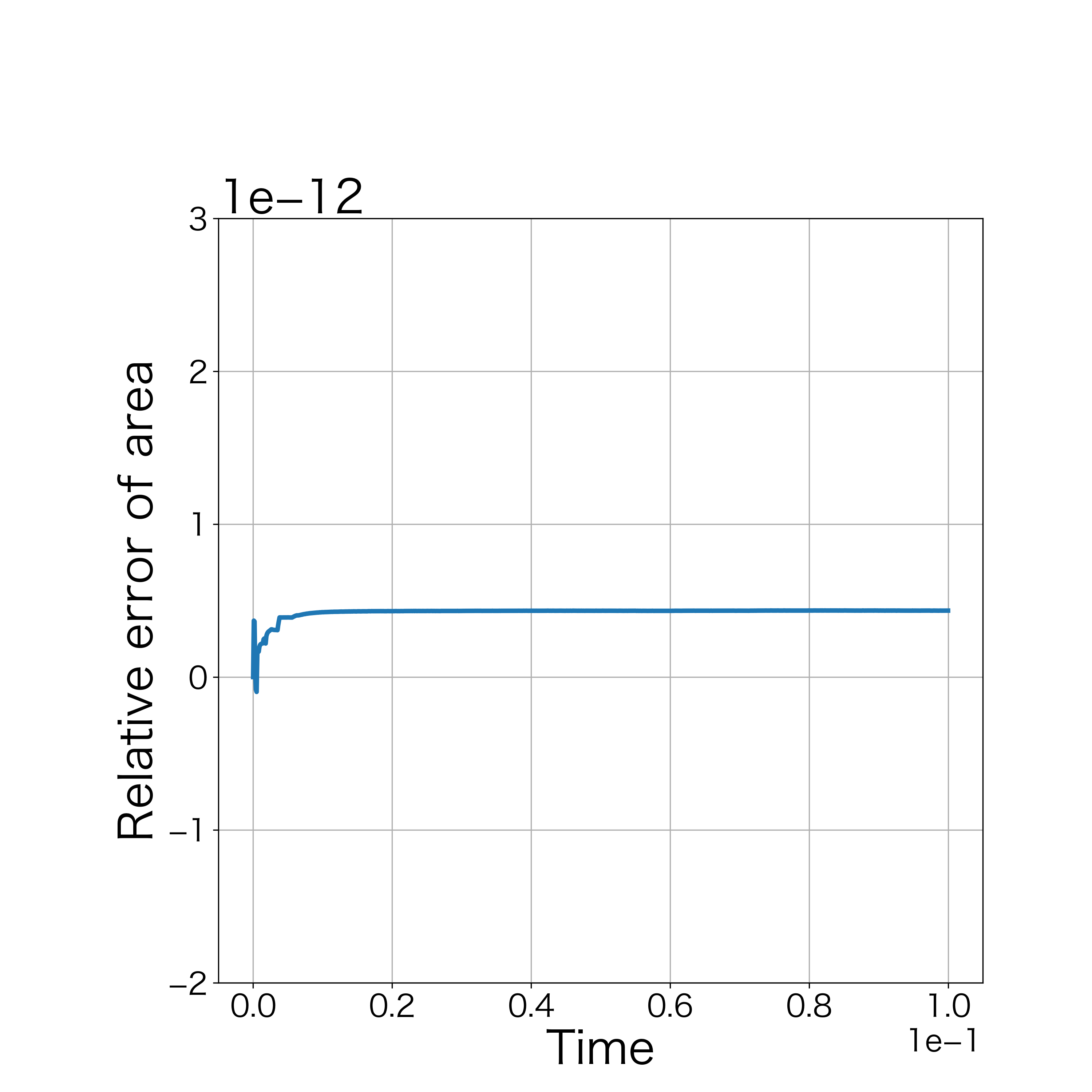}
    \end{center}
   \end{minipage}
   \caption{Relative error of the length and the enclosed area (rectangular curve).}
   \label{fig:29}
 \end{figure}

Numerical experiments are conducted with $N=40$ vertices and the rectangular initial curve illustrated in Fig.\ \ref{fig:26}. For the time step size $\Delta t=10^{-4}$, the constant $c_0=2$ and the tangential velocity parameter $\alpha=100$ are used. Numerical computations are stopped when the condition 
\begin{equation}
  \label{teishi}
  \frac{B^{(n)}-B^{(n+1)}}{B^{(n)}} < \epsilon
\end{equation}
is satisfied. In this examination, we set $\epsilon=10^{-5}$.
Numerical computation stopped at the 2567th step by condition \eqref{teishi}.
Fig.\ \ref{fig:27} demonstrates the initial polygonal curve and the numerical solutions at $n=50$, 100, 300, 500, and 1000 steps.  Fig.\ \ref{fig:28}--\ref{fig:29} respectively show the time evolution of the bending energy and the relative error of the length and the enclosed area for each step. The discrete bending energy's dissipative nature and the conservation of the length and the enclosed area can be seen.


\subsubsection{Examination. 5}
By \eqref{yo_gamma}, the value of $\gamma^{(n)}$ obtained by the scheme \eqref{hell_ski} is considered to approach 0 in the continuous limit. To confirm this numerically, two experiments are conducted: one with a fixed number of vertices and different time step sizes, and the other with a fixed time step size and different number of vertices. For each experiment, we check the value of $\gamma$ when the numerical experiment stops by condition \eqref{teishi}. The conditions changing time step size are as follows: the number of vertices $N=50$, the constant $c_0=2$, the tangential velocity parameter $\alpha=100$, and the stopping condition $\epsilon=10^{-5}$.
We use five time step sizes $\Delta t=10^{-6}$, $2.5\times10^{-6}$, $5\times10^{-6}$, $7.5\times10^{-6}$, and $10^{-5}$. The condition changing the number of vertices is conducted with the time step size $\Delta t=10^{-5}$, the constant $c_0=2$, the tangential velocity parameter $\alpha=100$, and the stopping condition $\epsilon=10^{-4}$. We use five initial polygonals with vertices $N=30$, 40, 50, 60, and 70. In these examination, we set initial polygonal curves by \eqref{kubire}.

\begin{figure}[H]
  \begin{minipage}{0.5\hsize}
   \begin{center}
    \includegraphics[width=55mm]{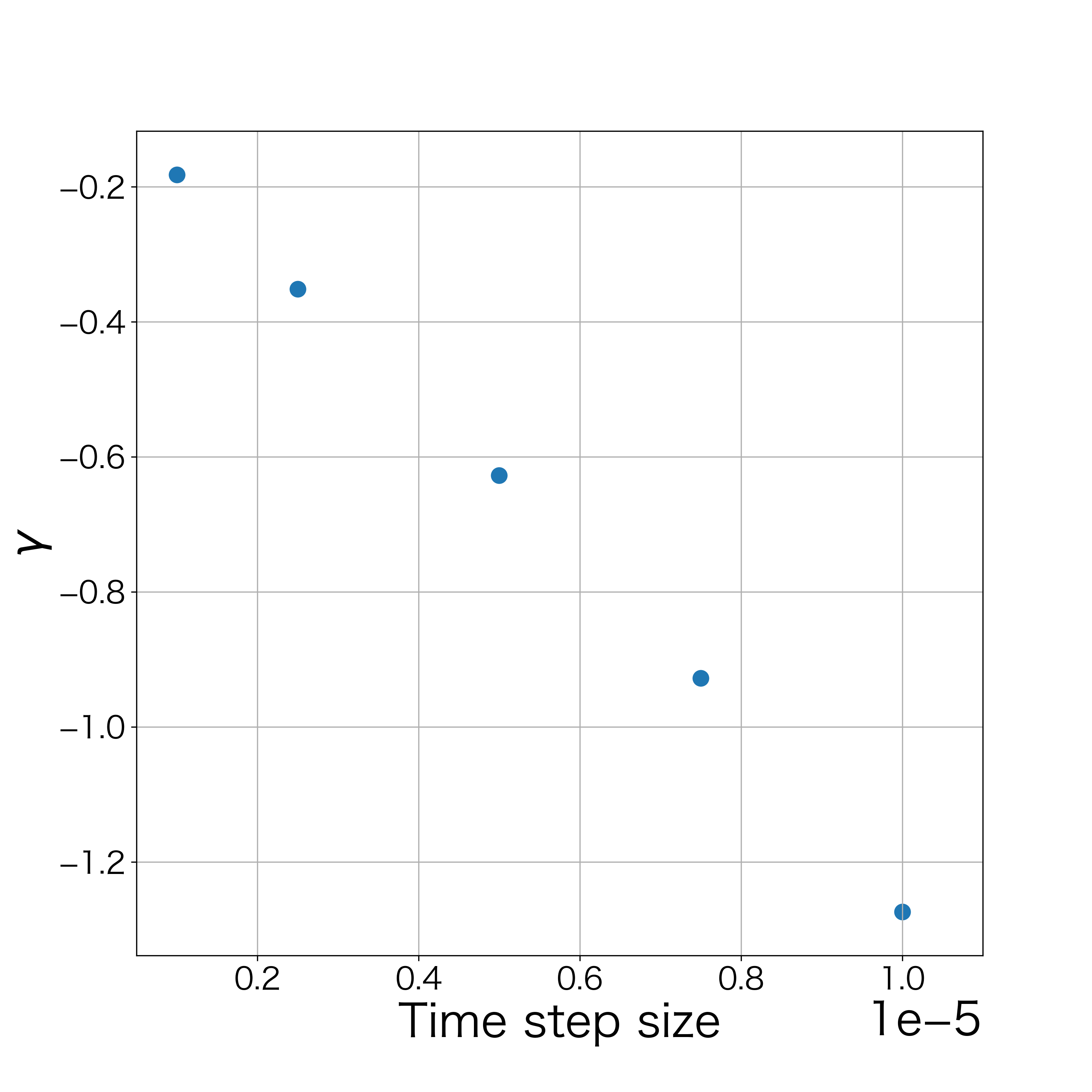}
   \end{center}
   \caption{Value of gamma when numerical\\ calculation stops (varying the time step size).}
   \label{fig:31}
  \end{minipage}
  \begin{minipage}{0.5\hsize}
    \begin{center}
     \includegraphics[width=55mm]{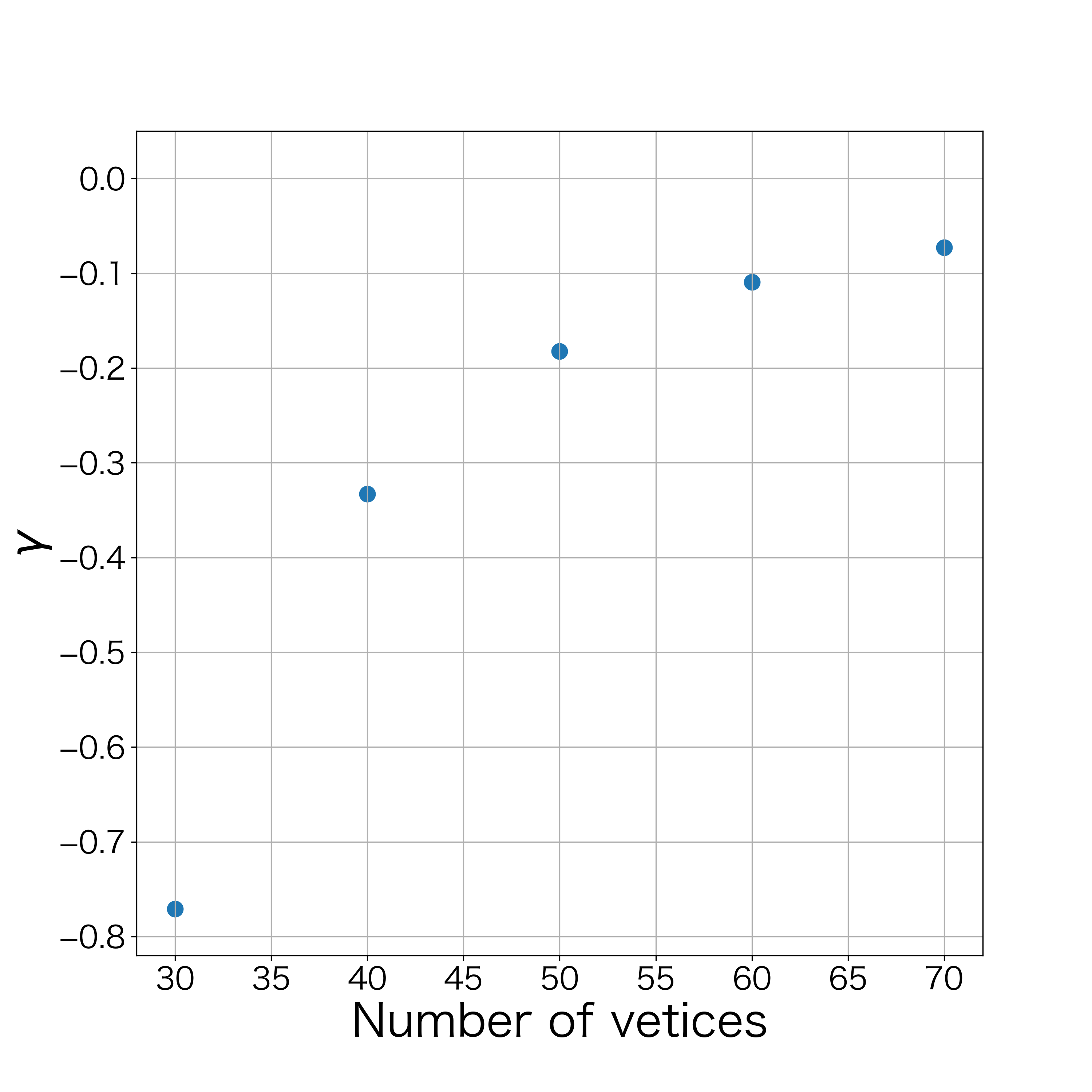}
    \end{center}
    \caption{Value of gamma when numerical \\calculation stops (varying the number of vertices).}
    \label{fig:32}
   \end{minipage}
 \end{figure}

In the experiments with various time step sizes, the value of $\gamma$ at the final step became smaller as the time step size decreased (Fig. \ref{fig:31}). In the experiments with a different number of vertices, the value of $\gamma$ at the final step also became smaller as the number of vertices increased (Fig.\ \ref{fig:32}). Therefore, in both experiments, the value of $\gamma$ may approach 0 in the continuous limit.

\section{Conclusion}
In this study, structure-preserving numerical methods for the Willmore and the Helfrich flows are developed. The energy structure was reproduced by discretizing the variational derivatives of the bending energy, the length, and the enclosed area based on DVDM and by determining Lagrange multipliers properly. Moreover, the tangential velocity by Deckelnick was appropriately introduced to prevent the concentration of vertices. Through numerical experiments, we can verify that the Willmore and the Helfrich flows can be computed correctly and our method reproduces the energy structures.
 
One of the problems of this method is that it needs the solution of nonlinear equations. Therefore, the well-posedness of the scheme is nontrivial and the computational cost is expensive. Thus, theoretical treatment and linearization of the scheme are future tasks. Additionally, the appropriate choice of $\alpha$ for tangential velocity is also a future task.




\section*{Acknowledgment}
This work was supported by JSPS KAKENHI Grant Numbers 19K14590 and 21K18301, Japan. The authors would like to thank Enago (\verb|www.enago.jp|) for the
English language review.
\appendix
\section{Appendix}
\subsection{Discrete variational derivative of the bending energy}
Let us calculate the discrete variational derivative of the bending energy. For simplicity, we denote $x=x^{(n+1)}$ and $\tilde{x}=x^{(n)}$ for $x=\boldsymbol{X}_i,\; k_i,\; r_i,$ etc. The bending energy of the $(n+1)$-th step can be transformed as 
\begin{equation}
  \begin{split}
  B&=\frac{1}{2}\sum_{i=1}^N(k_i-c_0)^2\hat{r_i}=\frac{1}{2}\sum_{i=1}^N(k_i-c_0)^2\frac{r_i+r_{i+1}}{2}\\
  &=\sum_{i=1}^N\frac{(k_i-c_0)^2+(k_{i-1}-c_0)^2}{4}r_i\\
  &=\sum_{i=1}^NC_ir_i,
  \end{split}
\end{equation}
where 
\begin{equation}
C_i=\frac{(k_i-c_0)^2+(k_{i-1}-c_0)^2}{4}.
\end{equation}
From the above form, the energy difference is 
\begin{equation}
  \begin{split}
   B-\tilde{B}=\sum_{i=1}^N\left(C_ir_i-\tilde{C}_i\tilde{r}_i\right)=\sum_{i=1}^N\frac{C_i+\tilde{C}_i}{2}(r_i-\tilde{r}_i)+\sum_{i=1}^N\frac{r_i+\tilde{r}_i}{2}(C_i-\tilde{C}_i),
  \end{split}
\end{equation}
where the formula 
\[
  ab-cd=\frac{a+c}{2}(b-d)+(a-c)\frac{b+d}{2}
\]
is used.
Now let 
\[A_1=\sum_{i=1}^N\frac{C_i+\tilde{C}_i}{2}(r_i-\tilde{r}_i),\quad A_2=\sum_{i=1}^N\frac{r_i+\tilde{r}_i}{2}(C_i-\tilde{C}_i).
\]
Then, we obtain 
\[
B-\tilde{B}=A_1+A_2.
\]
The first term $A_1$ can be transformed as 
\begin{equation}
  \begin{split}
  A_1&=\sum_{i=1}^N\frac{C_i+\tilde{C}_i}{2}(r_i-\tilde{r}_i)\\
  &=\sum_{i=1}^N\frac{C_i+\tilde{C}_i}{2}\frac{|\boldsymbol{X}_i-\boldsymbol{X}_{i-1}|^2-|\tilde{\boldsymbol{X}}_i-\tilde{\boldsymbol{X}}_{i-1}|^2}{|\boldsymbol{X}_i-\boldsymbol{X}_{i-1}|+
  |\tilde{\boldsymbol{X}}_{i}-\tilde{\boldsymbol{X}}_{i-1}|}\\
  &=\sum_{i=1}^N\frac{C_i+\tilde{C}_i}{2}\frac{r_i\boldsymbol{t}_i+\tilde{r}_i\tilde{\boldsymbol{t}}_i}{r_i+\tilde{r}_i}
  \cdot(\boldsymbol{X}_{i}-\tilde{\boldsymbol{X}}_{i}-\boldsymbol{X}_{i-1}+\tilde{\boldsymbol{X}}_{i-1})\\
  &=\sum_{i=1}^N \left(\frac{C_i+\tilde{C}_i}{2}\frac{r_i\boldsymbol{t}_i+\tilde{r}_i\tilde{\boldsymbol{t}}_i}{r_i+\tilde{r}_i}-
  \frac{C_{i+1}+\tilde{C}_{i+1}}{2}\frac{r_{i+1}\boldsymbol{t}_{i+1}+\tilde{r}_{i+1}\tilde{\boldsymbol{t}}_{i+1}}{r_{i+1}+\tilde{r}_{i+1}}\right)\cdot(\boldsymbol{X}_i-\tilde{\boldsymbol{X}_i})\\
  &=\sum_{i=1}^N\boldsymbol{D}_i\cdot(\boldsymbol{X}_i-\tilde{\boldsymbol{X}_i}),
  \end{split}
  \end{equation}
where 
\[\boldsymbol{D}_i=
\frac{C_i+\tilde{C}_i}{2}\frac{r_i\boldsymbol{t}_i+\tilde{r}_i\tilde{\boldsymbol{t}}_i}{r_i+\tilde{r}_i}-
\frac{C_{i+1}+\tilde{C}_{i+1}}{2}\frac{r_{i+1}\boldsymbol{t}_{i+1}+\tilde{r}_{i+1}\tilde{\boldsymbol{t}}_{i+1}}{r_{i+1}+\tilde{r}_{i+1}}.
\]
Subsequently, $A_2$ can be transformed as 
\begin{equation}
  \begin{split}
  A_2&=\sum_{i=1}^N\frac{r_i+\tilde{r}_i}{2}(C_i-\tilde{C}_i)\\
  &=\sum_{i=1}^N\frac{r_i+\tilde{r}_i}{2}\left\{\frac{(k_i-c_0)^2+(k_{i-1}-c_0)^2}{4}-\frac{(\tilde{k}_i-c_0)^2+(\tilde{k}_{i-1}-c_0)^2}{4}
  \right\}\\
  &=\sum_{i=1}^N\left(\frac{r_i+\tilde{r}_i}{8}+\frac{r_{i+1}+\tilde{r}_{i+1}}{8}\right)(k_i-c_0+\tilde{k}_i-c_0)(k_i-\tilde{k}_i)\\
  &=\sum_{i=1}^NE_i(k_i-\tilde{k}_i),
  \end{split}
  \end{equation}
  where 

\[
  E_i=\left(\frac{r_i+\tilde{r}_i}{8}+\frac{r_{i+1}+\tilde{r}_{i+1}}{8}\right)(k_i-c_0+\tilde{k}_i-c_0).
  \]

  Furthermore, $A_2$ can be transformed as 
  \begin{equation}
    \begin{split}
      \label{B}
    A_2&=\sum_{i=1}^NE_i(k_i-\tilde{k}_i)\\
    &=\sum_{i=1}^NE_i\left\{\frac{\mathrm{det}\left[\delta_i^{\left<1\right>}\boldsymbol{X}_i,\delta_i^{\left<2\right>}\boldsymbol{X}_i\right]}{g_i^3}-\frac{\mathrm{det}\left[\delta_i^{\left<1\right>}\tilde{\boldsymbol{X}}_i,\delta_i^{\left<2\right>}\tilde{\boldsymbol{X}}_i\right]}{\tilde{g}_i^3}\right\}\\
    &=\sum_{i=1}^NE_i\left(\frac{1}{g_i^3}+\frac{1}{\tilde{g}_i^3}\right)\frac{1}{2}\left(\mathrm{det}\left[\delta_i^{\left<1\right>}\boldsymbol{X}_i,\delta_i^{\left<2\right>}\boldsymbol{X}_i\right]-\mathrm{det}\left[\delta_i^{\left<1\right>}\tilde{\boldsymbol{X}}_i,\delta_i^{\left<2\right>}\tilde{\boldsymbol{X}}_i\right]\right)\\
    &+\sum_{i=1}^NE_i\frac{1}{2}\left(\mathrm{det}\left[\delta_i^{\left<1\right>}\boldsymbol{X}_i,\delta_i^{\left<2\right>}\boldsymbol{X}_i\right]+\mathrm{det}\left[\delta_i^{\left<1\right>}\tilde{\boldsymbol{X}}_i,\delta_i^{\left<2\right>}\tilde{\boldsymbol{X}}_i\right]\right)\left(\frac{1}{g_i^3}-\frac{1}{\tilde{g}_i^3}\right)\\
    &=\sum_{i=1}^NE_iF_i+\sum_{i=1}^NE_iG_i,
  \end{split}
  \end{equation}
where 
\[
  F_i=\left(\frac{1}{g_i^3}+\frac{1}{\tilde{g}_i^3}\right)\frac{1}{2}\left(\mathrm{det}\left[\delta_i^{\left<1\right>}\boldsymbol{X}_i,\delta_i^{\left<2\right>}\boldsymbol{X}_i\right]-\mathrm{det}\left[\delta_i^{\left<1\right>}\tilde{\boldsymbol{X}}_i,\delta_i^{\left<2\right>}\tilde{\boldsymbol{X}}_i\right]\right)
,\]
\[
G_i=\frac{1}{2}\left(\mathrm{det}\left[\delta_i^{\left<1\right>}\boldsymbol{X}_i,\delta_i^{\left<2\right>}\boldsymbol{X}_i\right]+\mathrm{det}\left[\delta_i^{\left<1\right>}\tilde{\boldsymbol{X}}_i,\delta_i^{\left<2\right>}\tilde{\boldsymbol{X}}_i\right]\right)\left(\frac{1}{g_i^3}-\frac{1}{\tilde{g}_i^3}\right).
\]
Transformations are then performed for $F_i$ and $G_i$.
Note that the component of $\boldsymbol{X}_i$ is 
\begin{equation}
  \boldsymbol{X}_i=
  \left(\begin{matrix}x_i\\y_i\end{matrix}\right).
\end{equation}
Firstly, $F_i$ can be transformed as 
\begin{equation}
  \begin{split}
F_i&=\frac{1}{2}\left(\frac{1}{g_i^3}+\frac{1}{\tilde{g}_i^3}\right)\left(\mathrm{det}\left[\delta_i^{\left<1\right>}\boldsymbol{X}_i,\delta_i^{\left<2\right>}\boldsymbol{X}_i\right]-\mathrm{det}\left[\delta_i^{\left<1\right>}\tilde{\boldsymbol{X}}_i,\delta_i^{\left<2\right>}\tilde{\boldsymbol{X}}_i\right]\right)\\
&=\frac{1}{2}\left(\frac{1}{g_i^3}+\frac{1}{\tilde{g}_i^3}\right)(\delta_i^{\left<1\right>}x_i \delta_i^{\left<2\right>}y_i-\delta_i^{\left<2\right>}x_i\delta_i^{\left<1\right>}y_i-
\delta_i^{\left<1\right>}\tilde{x}_i \delta_i^{\left<2\right>}\tilde{y}_i+\delta_i^{\left<2\right>}\tilde{x}_i\delta_i^{\left<1\right>}\tilde{y}_i)\\
&=\frac{1}{2}\left(\frac{1}{g_i^3}+\frac{1}{\tilde{g}_i^3}\right)
\left(\begin{matrix}\frac{\delta_i^{\left<2\right>}{y}_i+\delta_i^{\left<2\right>}\tilde{y}_i}{2}\\-\frac{\delta_i^{\left<2\right>}{x}_i+\delta_i^{\left<2\right>}\tilde{x}_i}{2}\end{matrix}\right)\cdot
\left(\begin{matrix}\delta_i^{\left<1\right>}{x}_i-\delta_i^{\left<1\right>}\tilde{x}_i  \\\delta_i^{\left<1\right>}{y}_i -\delta_i^{\left<1\right>}\tilde{y}_i \end{matrix}\right)\\
&+\frac{1}{2}\left(\frac{1}{g_i^3}+\frac{1}{\tilde{g}_i^3}\right)
\left(\begin{matrix}-\frac{\delta_i^{\left<1\right>}{y}_i+\delta_i^{\left<1\right>}\tilde{y}_i}{2}\\\frac{\delta_i^{\left<1\right>}{x}_i+\delta_i^{\left<1\right>}\tilde{x}_i}{2}\end{matrix}\right)\cdot
\left(\begin{matrix}\delta_i^{\left<2\right>}{x}_i-\delta_i^{\left<2\right>}\tilde{x}_i  \\\delta_i^{\left<2\right>}{y}_i -\delta_i^{\left<2\right>}\tilde{y}_i \end{matrix}\right)\\
& =\boldsymbol{H}_i\cdot \delta_i^{\left<1\right>}(\boldsymbol{X}_i-\tilde{\boldsymbol{X}}_i)+\boldsymbol{I}_i\cdot \delta_i^{\left<2\right>}(\boldsymbol{X}_i-\tilde{\boldsymbol{X}}_i)
,
  \end{split}
\end{equation}
where 
\[
  \boldsymbol{H}_i=\frac{1}{2}\left(\frac{1}{g_i^3}+\frac{1}{\tilde{g}_i^3}\right)
  \left(\begin{matrix}\frac{\delta_i^{\left<2\right>}{y}_i+\delta_i^{\left<2\right>}\tilde{y}_i}{2}\\-\frac{\delta_i^{\left<2\right>}{x}_i+\delta_i^{\left<2\right>}\tilde{x}_i}{2}\end{matrix}\right),
\]
\[
  \boldsymbol{I}_i=\frac{1}{2}\left(\frac{1}{g_i^3}+\frac{1}{\tilde{g}_i^3}\right)
  \left(\begin{matrix}-\frac{\delta_i^{\left<1\right>}{y}_i+\delta_i^{\left<1\right>}\tilde{y}_i}{2}\\\frac{\delta_i^{\left<1\right>}{x}_i+\delta_i^{\left<1\right>}\tilde{x}_i}{2}\end{matrix}\right).
\]
Subsequently, $G_i$ can be transformed as 
\begin{equation}
  \begin{split}
    &G_i=\frac{1}{2}\left(\mathrm{det}\left[\delta_i^{\left<1\right>}\boldsymbol{X}_i,\delta_i^{\left<2\right>}\boldsymbol{X}_i\right]+\mathrm{det}\left[\delta_i^{\left<1\right>}\tilde{\boldsymbol{X}}_i,\delta_i^{\left<2\right>}\tilde{\boldsymbol{X}}_i\right]\right)\left(\frac{1}{g_i^3}-\frac{1}{\tilde{g}_i^3}\right)\\
    &=-\frac{(g_i^2+g_i\tilde{g}_i+\tilde{g}_i^2)}{2g_i^3\tilde{g}_i^3}\left(\mathrm{det}\left[\delta_i^{\left<1\right>}\boldsymbol{X}_i,\delta_i^{\left<2\right>}\boldsymbol{X}_i\right]+\mathrm{det}\left[\delta_i^{\left<1\right>}\tilde{\boldsymbol{X}}_i,\delta_i^{\left<2\right>}\tilde{\boldsymbol{X}}_i\right]\right)({g_i}-\tilde{g}_i)\\
    &=G_{2,i}({g_i}-\tilde{g}_i),
  \end{split}
\end{equation}
where
\[
  G_{2,i}=-\frac{(g_i^2+g_i\tilde{g}_i+\tilde{g}_i^2)}{2g_i^3\tilde{g}_i^3}\left(\mathrm{det}\left[\delta_i^{\left<1\right>}\boldsymbol{X}_i,\delta_i^{\left<2\right>}\boldsymbol{X}_i\right]+\mathrm{det}\left[\delta_i^{\left<1\right>}\tilde{\boldsymbol{X}}_i,\delta_i^{\left<2\right>}\tilde{\boldsymbol{X}}_i\right]\right).
\]
Furthermore,
\begin{equation}
  \begin{split}
    G_i&=G_{2,i}({g_i}-\tilde{g}_i)\\
    &=\frac{G_{2,i}}{g_i+\tilde{g}_i}(|\delta_i^{\left<1\right>}{\boldsymbol{X}}_i|^2-|\delta_i^{\left<1\right>}\tilde{\boldsymbol{X}}_i|^2)\\
    &=\frac{G_{2,i}}{g_i+\tilde{g}_i}(\delta_i^{\left<1\right>}{\boldsymbol{X}}_i+\delta_i^{\left<1\right>}\tilde{\boldsymbol{X}}_i)\cdot(\delta_i^{\left<1\right>}{\boldsymbol{X}}_i-\delta_i^{\left<1\right>}\tilde{\boldsymbol{X}}_i)\\
    &=\frac{G_{2,i}}{g_i+\tilde{g}_i}(\delta_i^{\left<1\right>}{\boldsymbol{X}}_i+\delta_i^{\left<1\right>}\tilde{\boldsymbol{X}}_i)\cdot\delta_i^{\left<1\right>}({\boldsymbol{X}}_i-\tilde{\boldsymbol{X}}_i)\\
    &=\boldsymbol{L}_i\cdot\delta_i^{\left<1\right>}({\boldsymbol{X}}_i-\tilde{\boldsymbol{X}}_i),
  \end{split}
\end{equation}
where 
\[
  \boldsymbol{L}_i=\frac{G_{2,i}}{g_i+\tilde{g}_i}(\delta_i^{\left<1\right>}{\boldsymbol{X}}_i+\delta_i^{\left<1\right>}\tilde{\boldsymbol{X}}_i).
\]
Substituting $F_i$ and $G_i$ into \eqref{B} and using the  summation-by-parts formulas, we obtain
\begin{equation}
  \begin{split}
  A_2&=\sum_{i=1}^N\left\{E_i\boldsymbol{H}_i\cdot \delta_i^{\left<1\right>}(\boldsymbol{X}_i-\tilde{\boldsymbol{X}}_i)+E_i\boldsymbol{I}_i\cdot \delta_i^{\left<2\right>}(\boldsymbol{X}_i-\tilde{\boldsymbol{X}}_i)
+E_i\boldsymbol{L}_i\cdot\delta_i^{\left<1\right>}({\boldsymbol{X}}_i-\tilde{\boldsymbol{X}}_i)\right\} \\
&=\sum_{i=1}^N\left\{-\delta_i^{\left<1\right>}(E_i\boldsymbol{H}_i)+\delta_i^{\left<2\right>}(E_i\boldsymbol{I}_i)-\delta_i^{\left<1\right>}(E_i\boldsymbol{L}_i)\right\}\cdot(\boldsymbol{X}_i-\tilde{\boldsymbol{X}}_i).
\end{split}
\end{equation}
By transforming the equation so far, the difference of the bending energy can be transformed as 
\footnotesize
\begin{equation}
  \begin{split}
  &\frac{B^{(n+1)}-\tilde{B}^{(n)}}{\Delta t}\\
  &=\sum_{i=1}^N\frac{\boldsymbol{D}^{(n+1)}_i-\delta_i^{\left<1\right>}(E^{(n+1)}_i\boldsymbol{H}^{(n+1)}_i)+\delta_i^{\left<2\right>}(E^{(n+1)}_i\boldsymbol{I}^{(n+1)}_i)-\delta_i^{\left<1\right>}(E^{(n+1)}_i\boldsymbol{L}^{(n+1)}_i)}{\hat{r}^{(n+1/2)}_i}\cdot\frac{\boldsymbol{X}^{(n+1)}_i-\tilde{\boldsymbol{X}^{(n)}}_i}{\Delta t}\hat{r}^{(n+1/2)}_i. 
  \end{split}
\end{equation}
\normalsize
We can define $\boldsymbol{\delta}_dB$ by
\begin{equation}
  \boldsymbol{\delta}_dB_i^{(n+1)}=\frac{\boldsymbol{D}^{(n+1)}_i-\delta_i^{\left<1\right>}(E^{(n+1)}_i\boldsymbol{H}^{(n+1)}_i)+\delta_i^{\left<2\right>}(E^{(n+1)}_i\boldsymbol{I}^{(n+1)}_i)-\delta_i^{\left<1\right>}(E^{(n+1)}_i\boldsymbol{L}^{(n+1)}_i)}{\hat{r}^{(n+1/2)}_i}.
\end{equation}
\subsection{Discrete variational derivative of the length}
By calculating the time difference of the discrete length as well as the bending energy, we obtain
\footnotesize
\begin{equation}
  \begin{split}
    &\frac{L^{(n+1)}-L^{(n)}}{\Delta t}=\sum_{i=1}^N\left(r_i^{(n+1)}-r_i^{(n)}\right)\\
  &=\sum_{i=1}^N\frac{|\boldsymbol{X}^{(n+1)}_i-\boldsymbol{X}^{(n+1)}_{i-1}|^2-|\boldsymbol{X}_i^{(n)}-\boldsymbol{X}^{(n)}_{i-1}|^2}{|\boldsymbol{X}^{(n+1)}_i-\boldsymbol{X}^{(n+1)}_{i-1}|-
  |\boldsymbol{X}^{(n)}_{i}-\boldsymbol{X}^{(n)}_{i-1}|}\\
  &=\sum_{i=1}^N\frac{r_i^{(n+1)}\boldsymbol{t}^{(n+1)}_i+r^{(n)}_i\boldsymbol{t}_i^{(n)}}{r_i^{(n+1)}+r_i^{(n)}}
  \cdot(\boldsymbol{X}^{(n+1)}_{i}-\boldsymbol{X}^{(n)}_{i}-\boldsymbol{X}^{(n+1)}_{i-1}+\boldsymbol{X}^{(n)}_{i-1})\\
  &=\sum_{i=1}^N \left(\frac{r_i^{(n+1)}\boldsymbol{t}^{(n+1)}_i+r^{(n)}_i\boldsymbol{t}_i^{(n)}}{r_i^{(n+1)}+r_i^{(n)}}-
  \frac{r_{i+1}^{(n+1)}\boldsymbol{t}^{(n+1)}_{i+1}+r^{(n)}_{i+1}\boldsymbol{t}_{i+1}^{(n)}}{r_{i+1}^{(n+1)}+r_{i+1}^{(n)}}\right)\cdot(\boldsymbol{X}^{(n+1)}_i-\boldsymbol{X}_i^{(n)})\\
  &=\sum_{i=1}^N \frac{1}{\hat{r}_i^{(n+1/2)}}\left(\frac{r_i^{(n+1)}\boldsymbol{t}^{(n+1)}_i+r^{(n)}_i\boldsymbol{t}_i^{(n)}}{r_i^{(n+1)}+r_i^{(n)}}-
  \frac{r_{i+1}^{(n+1)}\boldsymbol{t}^{(n+1)}_{i+1}+r^{(n)}_{i+1}\boldsymbol{t}_{i+1}^{(n)}}{r_{i+1}^{(n+1)}+r_{i+1}^{(n)}}\right)\cdot(\boldsymbol{X}^{(n+1)}_i-\boldsymbol{X}_i^{(n)})\hat{r}_i^{(n+1/2)}.
  \end{split}
\end{equation}
\normalsize
Thus, the discrete variational derivative of the length is defined by
\begin{equation}
  \boldsymbol{\delta}_d L_i^{(n+1)}=\frac{1}{\hat{r}_i^{(n+1/2)}}\left(\frac{r_i^{(n+1)}\boldsymbol{t}^{(n+1)}_i+r^{(n)}_i\boldsymbol{t}_i^{(n)}}{r_i^{(n+1)}+r_i^{(n)}}-
  \frac{r_{i+1}^{(n+1)}\boldsymbol{t}^{(n+1)}_{i+1}+r^{(n)}_{i+1}\boldsymbol{t}_{i+1}^{(n)}}{r_{i+1}^{(n+1)}+r_{i+1}^{(n)}}\right).
\end{equation}
\subsection{Discrete variational derivative of the area}
Finally, we calculate the discrete variation of the area.
The energy change can be transformed as 
\begin{equation}
  \begin{split}
    &\frac{A^{(n+1)}-A^{(n)}}{\Delta t}=\frac{1}{2}\sum_{i=1}^N \left({\boldsymbol{X}^{(n+1)}_{i-1}}^{\perp}\cdot\boldsymbol{X}^{(n+1)}_{i}-{\boldsymbol{X}^{(n)}_{i-1}}^{\perp}\cdot\boldsymbol{X}^{(n)}_{i}\right)\\
    &=\frac{1}{2}\sum_{i=1}^N \left(
    x^{(n+1)}_{i-1}y^{(n+1)}_{i}-x^{(n+1)}_{i}y^{(n+1)}_{i-1}-x^{(n)}_{i-1}y^{(n)}_{i}+x^{(n)}_{i}y^{(n)}_{i-1}\right)\\
    &=\sum_{i=1}^N \left(\begin{matrix}
      \left(-y^{(n+1)}_{i-1} -y^{(n)}_{i-1}+ y^{(n+1)}_{i+1}+ y^{(n)}_{i+1}\right)/4 \\
      \left(x^{(n+1)}_{i-1} +x^{(n)}_{i-1}-x^{(n+1)}_{i+1}-x^{(n)}_{i+1}\right)/4 
      \end{matrix}\right)\cdot(\boldsymbol{X}^{(n+1)}_i-\boldsymbol{X}_i^{(n)})\\
      &=\sum_{i=1}^N\frac{1}{\hat{r}_i^{(n+1/2)}} \left(\begin{matrix}
        \left(-y^{(n+1)}_{i-1} -y^{(n)}_{i-1}+ y^{(n+1)}_{i+1}+ y^{(n)}_{i+1}\right)/4 \\
        \left(x^{(n+1)}_{i-1} +x^{(n)}_{i-1}-x^{(n+1)}_{i+1}-x^{(n)}_{i+1}\right)/4 
        \end{matrix}\right)\cdot(\boldsymbol{X}^{(n+1)}_i-\boldsymbol{X}_i^{(n)})\hat{r}_i^{(n+1/2)},
  \end{split}
\end{equation}
where $A^{\perp}$ denotes $\begin{pmatrix}
  -a_2\\
  a_1
\end{pmatrix}$
for 
$\boldsymbol{A}=\begin{pmatrix}
  a_1\\
  a_2
\end{pmatrix}$.
Thus, discrete variational derivative of the area is defined by
\begin{equation}
  \boldsymbol{\delta}_d A_i^{(n+1)}=\frac{1}{\hat{r}_i^{(n+1/2)}}
    \left(\begin{matrix}
    \left(-y^{(n+1)}_{i-1} -y^{(n)}_{i-1}+ y^{(n+1)}_{i+1}+ y^{(n)}_{i+1}\right)/4 \\
    \left(x^{(n+1)}_{i-1} +x^{(n)}_{i-1}-x^{(n+1)}_{i+1}-x^{(n)}_{i+1}\right)/4 
    \end{matrix}
    \right).
\end{equation}

\end{document}